\documentclass[11 pt, reqno]{amsart}
\usepackage[margin=1in]{geometry}
\usepackage{amssymb, amsmath, amsthm}
\usepackage{amstext, amsfonts}
\usepackage[T1]{fontenc}
\usepackage{fourier}

\usepackage[usenames, dvipsnames, svgnames, x11names, hyperref]{xcolor}
\usepackage[pagebackref, colorlinks=true, citecolor={blue!80!black}, linkcolor={blue!80!black}, urlcolor={blue!80!black}, linktoc=section]{hyperref}
\usepackage{wasysym}
\usepackage{setspace}
\usepackage{cancel}
\usepackage[all]{xy}
\usepackage{empheq}
\usepackage{extpfeil}
\usepackage{graphicx}
\usepackage{ifthen}
\usepackage{bbm}

\makeatletter
\newtheorem*{rep@theorem}{\rep@title}
\newcommand{\newreptheorem}[2]{%
\newenvironment{rep#1}[1]{%
 \def\rep@title{#2 \ref{##1}}%
 \begin{rep@theorem}}%
 {\end{rep@theorem}}}
\makeatother

\graphicspath{{images/}, {./../images/}, {new-images/}}

\numberwithin{equation}{section}

\newtheorem{theorem}{Theorem}[section]
\newtheorem{mainthm}{Theorem}

\newreptheorem{theorem}{Theorem}
\newtheorem{lemma}[theorem]{Lemma}
\newreptheorem{lemma}{Lemma}
\newtheorem{proposition}[theorem]{Proposition}
\newtheorem{corollary}[theorem]{Corollary}
\newtheorem{conjecture}[theorem]{Conjecture}
\newtheorem{question}[theorem]{Question}

\theoremstyle{definition}
\newtheorem{remark}[theorem]{Remark}
\newtheorem{definition}[theorem]{Definition}

\newtheorem{example}[theorem]{Example}

\def\beq{\begin{eqnarray*}}
\def\eeq{\end{eqnarray*}}
\def\Q{\mathbb{Q}}

\def\R{\mathbb{R}}
\def\Z{\mathbb{Z}}

\def\incl{\hookrightarrow}
\def\to{\rightarrow}

\def\eps{\varepsilon}

\def\dim{\mathrm{dim}\>}

\def\x{\times}
\def\d{\partial}

\def\phi{\varphi}

\def\Emb{\mathrm{Emb}}
\def\BrEmb{\mathrm{BrEmb}}

\def\Diff{\mathrm{Diff}}
\def\Link{\mathrm{Link}}
\def\Conf{\mathrm{Conf}}

\def\P{\mathcal{P}}

\def\K{\mathcal{K}}
\def\L{\mathcal{L}}

\def\Map{\mathrm{Map}}
\def\Aut{\mathrm{Aut}}

\def\im{\mathrm{im}}

\def\hofib{\mathrm{hofib}}

\def\scu{\sqcup}

\newcommand{\PE}[3][asdf]
{
\ifthenelse{\equal{#1}{asdf}}
{\mathcal{P}_{#2, \, #2}^{#3}}
{\mathcal{P}_{#1, \,  #2}^{#3}}
}

\newcommand{\LL}[3][asdf]
{
\ifthenelse{\equal{#1}{asdf}}
{\mathcal{L}_{#2, \, #2}^{#3}}
{\mathcal{L}_{#1, \, #2}^{#3}}
}

\def\del{\partial}

   \title{Graphing, homotopy groups of spheres, and spaces of long links and knots}
    \author{Robin Koytcheff}
    \email{koytcheff@louisiana.edu}
    \address{Department of Mathematics, University of Louisiana at Lafayette, Lafayette, LA, USA 70504}

\makeatletter
\@namedef{subjclassname@2020}{%
  \textup{2020} Mathematics Subject Classification}
\makeatother

\begin{document}

\begin{abstract}
We study homotopy groups of spaces of long links in Euclidean space of codimension at least three.  With multiple components, they admit split injections from homotopy groups of spheres.  We show that, up to knotting, these account for all the homotopy groups in a range which depends on the dimensions of the source manifolds and target manifold and which roughly generalizes the triple-point-free range for isotopy classes.  Just beyond this range, joining components sends both a parametrized long Borromean rings class and a Hopf fibration to a generator of the first nontrivial homotopy group of the space of long knots.  For spaces of equidimensional long links of most source dimensions, we describe generators for the homotopy group in this degree in terms of these Borromean rings and homotopy groups of spheres.  A key ingredient in most of our results is a graphing map which increases source and target dimensions by one.
\end{abstract}

\maketitle

\tableofcontents

\section{Introduction}

This paper concerns homotopy groups of spaces of links of various dimensions, where the source and target manifolds are either Euclidean spaces or spheres.  We focus mainly on spaces $\mathcal{L}_{p_1,\dots,p_m}^n:=\Emb_c \left( \coprod_{i=1}^m \R^{p_i}, \, \R^n\right)$ of long links, meaning smooth embeddings $\coprod_{i=1}^m \R^{p_i} \incl \R^n$ with fixed compact support.  For the most part, we will assume $p_1,\dots, p_m\leq n-3$.  
Extensive progress that has been made on the rational homotopy type of these spaces,
but less is known over the integers, and less is known about explicit ``geometric'' generators of these homotopy groups than about their ranks.

We study map 
$G:\Omega \mathcal{L}_{p_1,\dots,p_m}^n \to \mathcal{L}_{p_1+1,\dots,p_m+1}^{n+1}$ 
given by graphing, in the sense of the graph of a function.
It produces an embedding from a based loop of embeddings by incorporating the loop parameter into both the source and target; see Definition \ref{D:Graphing}.  It is most easily visualized when $n=2$ and $p_1 =\dots =p_m=0$.  Extending work of Budney \cite{Budney:Family} from one component to multiple components, we use it to determine the homotopy groups of $\mathcal{L}_{p_1,\dots,p_m}^n$ in a range, with the key ingredient being Goodwillie's connectivity estimates for pseudoisotopy embedding spaces \cite{Goodwillie:Thesis, Goodwillie:CalcI, Goodwillie-Klein:2015}.
Then we explicitly describe generators for these groups, up to describing those for homotopy groups of spheres.
Our main results are organized as follows:
\begin{itemize}
\item 
Theorem \ref{T:A} concerns injectivity on homotopy groups of graphing.  It is an easy generalization of a theorem on isotopy classes of links to higher homotopy groups of spaces of links.  We use it to prove other main results.
\item 
Theorems \ref{T:B} and \ref{T:C} calculate, modulo knotting, the homotopy groups $\pi_i$ of spaces $\mathcal{L}_{p_1,\dots,p_m}^n$ of long links, in roughly the ``metastable'' or ``triple-point-free'' range.  If $p_1=\dots=p_m=p$, this range is just below $i=2n-3p-3$; there is no knotting in this range (see Corollary \ref{C:DirectSum}); and all the classes are in the image of the map induced by iterated graphing from spheres.  
Theorem \ref{T:D} extends this calculation to $i=2n-3p-3$ for $p_1=\dots=p_m=p$ with mild lower bounds on $p$.
\item 
Theorem \ref{T:E} gives explicit generators for the group $\pi_{2n-3p-3}\mathcal{L}_{p,\dots, p}^n$ calculated in Theorem \ref{T:D}.  They are described for $p\geq 3$ and any number $m$ of components, as well as for $p=2$ and 2 components.  
Modulo torsion, we obtain generators for $p \geq1$ and any $m$.
For the space $\mathcal{K}_p^n:=\Emb_c(\R^p, \R^n)$ of long knots,
 generators of the previously known group $\pi_{2n-3p-3}\mathcal{K}_p^n$ are described in terms of  $2$- and $3$-component links and homotopy groups of spheres.  Ultimately, all the generators in Theorem \ref{T:E} come from either $\pi_{2n-2p-3} S^{n-p-1}$ or an analogue of the Borromean rings.
\end{itemize}
Each batch of results is proven using a different method.
Before describing them in more detail, we survey some earlier related work.

\subsection{Previous related results}
In all of the following previous results, the codimensions of the embeddings are assumed to be at least 3.  An early antecedent of our approach here is the work of Zeeman \cite{Zeeman:LinkingSpheres}, who established an injection of $\pi_p(S^{n-p-1})$ into the set $\pi_0 \Emb(S^p \scu S^p, \,  S^n)$ of isotopy classes of spherical links.  Shortly afterwards, Haefliger \cite{Haefliger:1966CMH} determined this set in a range of values of $p,q$, and $n$ and showed that it is an abelian group.  This range was improved slightly by M.~Skopenkov \cite{Skopenkov:2009}.  Their result applies in the 2-metastable range. This is roughly the quadruple-point-free range, meaning that if $p=q$, then roughly $p < 3n/4$.  More precisely their range is $3n -2p-2q \geq  6$.  Crowley, Ferry, and M.~Skopenkov \cite{Crowley-Ferry-Skopenkov} computed rational isotopy classes of spherical links with an arbitrary number of components.  
Songhafouo Tsopm\'en\'e and Turchin \cite[Theorem 3.2]{Songhafouo-Turchin:Forum} described rational isotopy classes of long links in terms of trivalent trees.  
They conjectured an extension to all rational homotopy groups of spaces of long links in terms of the homology of a graph complex, which was proven in the extensive work of Fresse, Turchin, and Willwacher \cite[Section 5.5]{Fresse-Turchin-Willwacher:Emb}.

Haefliger \cite{Haefliger:1966Annals} showed that the group $\pi_0 \Emb(S^p, \,  S^n)$ 
(which is isomorphic to $\pi_0\K_p^n$) 
 is trivial in the metastable range, which is roughly the triple-point-free range $p<2n/3$.  More precisely, there is no knotting if $2n -3p \geq  4$.  He also computed this group for the first family of dimensions $(p,n)$ where they are nontrivial, namely where $2n-3p=3$.  Using iterates of the graphing map, Budney \cite{Budney:Family} generalized these results from isotopy classes to families of knots, specifically from $\pi_0\K_p^n$ to $\pi_i\K_p^n$ where $0 \leq i \leq 2n-3p-3$.  He showed that $\pi_i \K_p^n=0$ if $i\leq 2n-3p-4$ and calculated $\pi_{2n-3p-3}\K_p^n$.

\subsection{Main results and organization}
Our results use the same iterated graphing construction, adapted easily from long knots to long links.  
To get nontrivial classes for knots from iterated graphing, one must start with at least 1-dimensional knots, whereas in our setting of links, we can start with links where one or both components have dimension 0.  The space of such links has a subspace homeomorphic to a sphere, and most of our results involve classes from homotopy groups of spheres.  
Though some of our results hold in codimension less than 3, they say nothing new in these cases.

Our first main result gives nontrivial homotopy classes in spaces of 2-component links from homotopy groups of spheres with little restriction on the dimensions involved:  

\begin{mainthm}
\label{T:A}
If  $\ 0 \leq p \leq q \leq n-1$ and $i\geq 0$, then $\pi_i \Emb_c(\R^p \scu \R^q, \, \R^n)$ contains a direct summand of $\pi_{i+p}S^{n-q-1}$.
An inclusion of this summand is given by composing maps induced by a based homotopy equivalence 
$S^{n-q-1} \to \R^{n-p} - \R^{q-p}$,
the inclusion 
$\Emb_c(\ast, \R^{n-p} - \R^{q-p}) \incl \Emb_c(\ast \scu \R^{q-p}, \, \R^{n-p})$, 
and the $p$-fold graphing map 
$G^p: \Omega^p \Emb_c(\ast \scu \R^{q-p}, \, \R^{n-p}) \to \Emb_c(\R^p \scu \R^q, \, \R^n)$.
\end{mainthm}

Theorem \ref{T:A} appears in the main body as Theorem \ref{T:Lambda}.  In the Appendix, we prove Theorem \ref{T:LambdaSpherical} and Theorem \ref{T:Alpha}, which are analogues of it for spherical links and for link maps (i.e., smooth maps such that the images of the components are disjoint).
Putting $i=0$ gives the above-mentioned result of 
Zeeman.
The case $i=0$ and $p+q=n-1$ corresponds to classes dual to the generalized linking number.
See Example \ref{Ex:LkNumExamples}.  
Theorem \ref{T:A} clearly yields many nontrivial torsion classes in spaces of links.  
The proof of Theorem \ref{T:A} relies on showing that graphing is the inclusion of a retract up to homotopy.

\smallskip

In our second set of results, we determine certain homotopy groups.
The first among them mutually extends to $\pi_i\Emb_c(\R^p \scu \R^q, \, \R^n)$ both Budney's result 
on $\pi_i\Emb_c(\R^p,\, \R^n)$ and the result of Haefliger 
and M.~Skopenkov 
on $\pi_0\Emb (S^p \scu S^q, \, \R^n)$.  Indeed, Lemma \ref{L:SphericalLinks} identifies spherical isotopy classes with long isotopy classes.
The next theorem applies in a range that generalizes the triple-point-free range for isotopy classes to $i$-parameter families.  
There are however many nontrivial groups for links in this range, in contrast to the setting of knots.  
Recall that $\LL[p]{q}{n}:=\Emb_c(\R^p \scu \R^q, \, \R^n)$ and $\K_p^n :=\Emb_c(\R^p, \, \R^n)$.  Theorem \ref{T:B} appears as Theorem \ref{T:pi_iLinksIso} in the main body of the paper.

\begin{mainthm}
\label{T:B}
If $1 \leq p \leq q \leq n-3$ and $i \leq 2n - p - 2q -4$, then in the sequence of maps 
\[
\pi_{i+p} \LL[0]{q-p}{n-p} 
\overset{G_*}{\longrightarrow}
\pi_{i+p-1} \LL[1]{q-p+1}{n-p+1} \overset{G_*}{\longrightarrow}
 \dots 
 \overset{G_*}{\longrightarrow}  
 \pi_i \LL[p]{q}{n} \overset{G_*}{\longrightarrow} \dots 
 \overset{G_*}{\longrightarrow} \pi_0\LL[i+p]{i+q}{i+n}
\]
induces by graphing, each map is an isomorphism, except possibly the first.
The first map is always a surjection, and it is an isomorphism if $i \leq 2n-p-2q-5$ or $p=q$.
\end{mainthm}

Graphing preserves the quantity $2n-p-2q-4-i$, that is, replacing $i,p,q,$ and $n$ by the corresponding four parameters in any term in the sequence gives the same number and thus preserves the assumed inequality involving them.
Corollary \ref{C:DirectSum} describes $\pi_i\LL[p]{q}{n}$ as $\pi_{i+p}S^{n-q-1} \oplus \pi_0 \K_{i+q}^{i+n}$.
If $n-q\geq 2$, then we can identify $\pi_0 \K_{i+q}^{i+n}$ with isotopy classes of spherical knots in the sphere or  Euclidean space.  This group vanishes if $2n-3q \geq 4$, hence it vanishes in the assumed range if $p=q$.
It is also known if $2n-3q$ is $3$ or $2$, by work of Haefliger \cite{Haefliger:1966Annals} and Milgram \cite{Milgram:1972}.  The latter further identifies the 2-primary part of this group for some smaller values of this quantity.  See the exposition by A.~Skopenkov \cite{Skopenkov:Knots}. 
The proof of Theorem \ref{T:B} uses Goodwillie's connectivity results on pseudoisotopy embedding spaces, much like the result of Budney's that it generalizes.  
The iterated graphing in Theorems \ref{T:A} and \ref{T:B} becomes simpler when $p=q$, in which case the source space $\LL[0]{q-p}{n-p}$ is just the configuration space of two points in $\R^{n-p}$.  


Theorem \ref{T:C} appears as Theorem \ref{T:ImagesOfiS} in the main body.
It says that for spaces of links with $m$ components, in a range analogous to that in Theorem \ref{T:B}, all elements of those homotopy groups come from links with at most 2 components.  
More specifically, for $i\leq 2n - p_1 - p_{m-1} - p_m - 4$, it allows us to describe all classes in $\pi_i \mathcal{L}_{p_1, \dots, p_m}^n$ as in Corollary \ref{C:DirectSum}.
Its proof uses the Hilton--Milnor theorem, restriction fibrations, and a homotopy retract as in the proof of Theorem \ref{T:A}

\begin{mainthm}
\label{T:C}
Suppose that $0 \leq \ell \leq m$, $1 \leq p_1 \leq \dots \leq p_m \leq n-3$, and 
$0 \leq i < 1-p_1+\sum_{k=m-\ell+1}^m (n-p_k-2)$.
Then every class in $\pi_i \mathcal{L}_{p_1, \dots, p_m}^n$ is in 
$\displaystyle \sum_{S\subset \{1,\dots,m\}, |S|\leq \ell} \im (\iota_S)_*$, where $\iota_S$ is the inclusion of the subspace of links with components labeled by a subset $S$ of $\{1,\dots,m\}$.
\end{mainthm}


Next, methods like those used to prove Theorem \ref{T:B} yield a calculation of homotopy groups of equidimensional 2-component links of dimension at least 2 in the degree where graphing from spheres ceases to be surjective.
Theorem \ref{T:D} is a slightly abbreviated version of Theorem \ref{T:GraphingDim>1}.

\begin{mainthm}
\label{T:D}
Suppose $1 \leq p \leq n-3$.  
\begin{itemize}
\item[(a)]
For 2-component links, consider the sequence of maps induced by graphing:
\begin{equation*}
\pi_{2n-3p-3}\LL{p}{n} \to \pi_{2n-3p-4}\LL{p+1}{n+1} \to \dots \to \pi_0 \LL{2n-2p-3}{3n-3p-3}.
\end{equation*}
If $p\geq 2$, then all the maps are isomorphisms, and these groups are isomorphic to 
\begin{align*}
\begin{array}{ll}
\Z^3 \oplus \pi_{2n-2p-3} S^{n-p-1} & \text{ if $n-p$ is odd}\\
(\Z/2)^3 \oplus \pi_{2n-2p-3} S^{n-p-1} & \text{ if $n-p$ is even}. 
\end{array}
\end{align*}
If $p=1$, then the first map is surjective, and the remaining maps are isomorphisms.  
\item[(b)]
For 3-component links, consider the sequence of maps induced by graphing:
\begin{equation*}
\pi_{2n-3p-3}\mathcal{L}_{3\cdot p}^{n} \to \pi_{2n-3p-4}\mathcal{L}_{3\cdot(p+1)}^{n+1} \to \dots \to \pi_0 \mathcal{L}_{3\cdot(2n-2p-3)}^{3n-3p-3}.
\end{equation*}
If $p\geq 3$, then all the maps are isomorphisms, and these groups are isomorphic to 
\begin{align*}
\begin{array}{ll}
\Z^7 \oplus \left(\pi_{2n-2p-3} S^{n-p-1}\right)^3 & \text{ if $n-p$ is odd}\\
\Z \oplus (\Z/2)^6 \oplus \left(\pi_{2n-2p-3} S^{n-p-1}\right)^3 & \text{ if $n-p$ is even}. 
\end{array}
\end{align*}
If $p=2$, then the first map is surjective, and the remaining maps are isomorphisms.  
\end{itemize}
\end{mainthm}


Our last main result gives explicit generators of these groups for equidimensional long links.  It connects spaces of 2- and 3-component pure braids to  spaces of long knots and 2-component long links, in the homotopy group just outside the ranges of Theorems \ref{T:B} and \ref{T:C}, i.e., the group $\pi_{2n-3p-3}$, which appeared in Theorem \ref{T:D}.  

\begin{figure}[h!]
\includegraphics[scale=0.25]{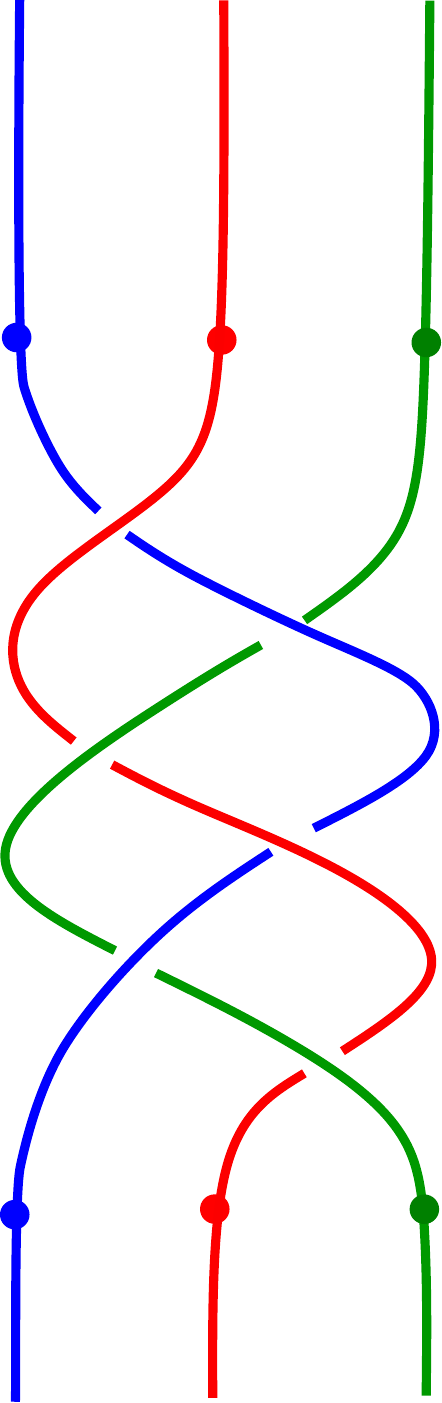}
\qquad
\raisebox{6pc}{$\xmapsto{J}$}
\qquad
\includegraphics[scale=0.25]{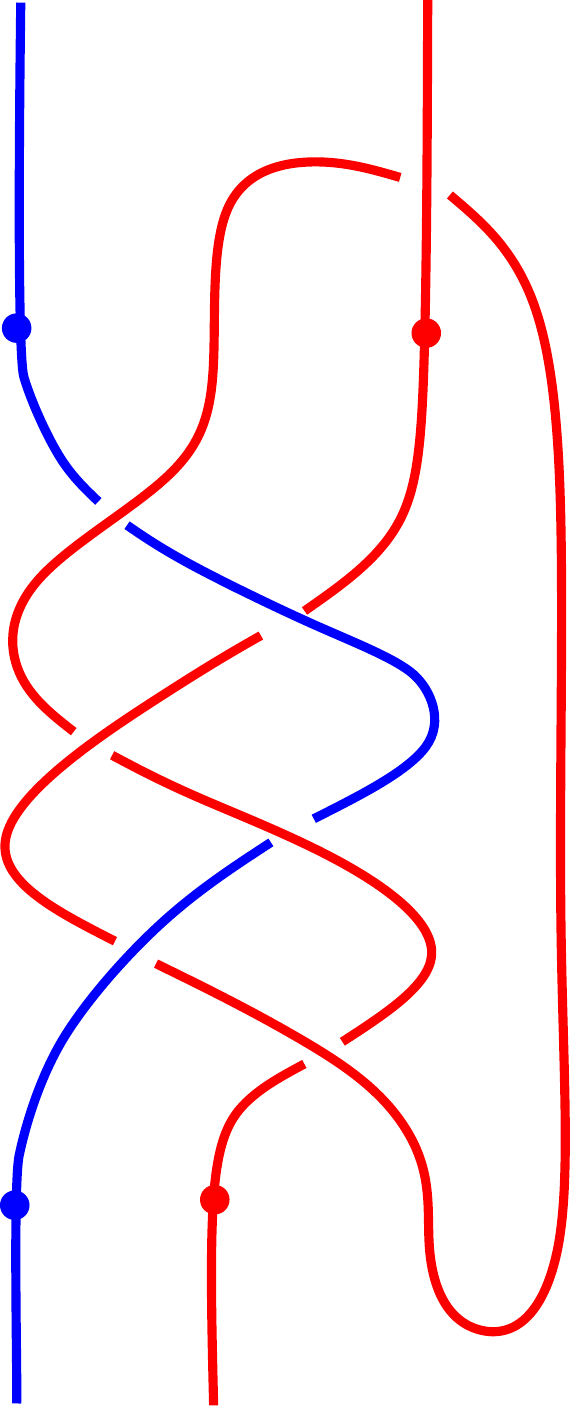}
\qquad
\raisebox{6pc}{$\xmapsto{J}$}
\qquad
\includegraphics[scale=0.25]{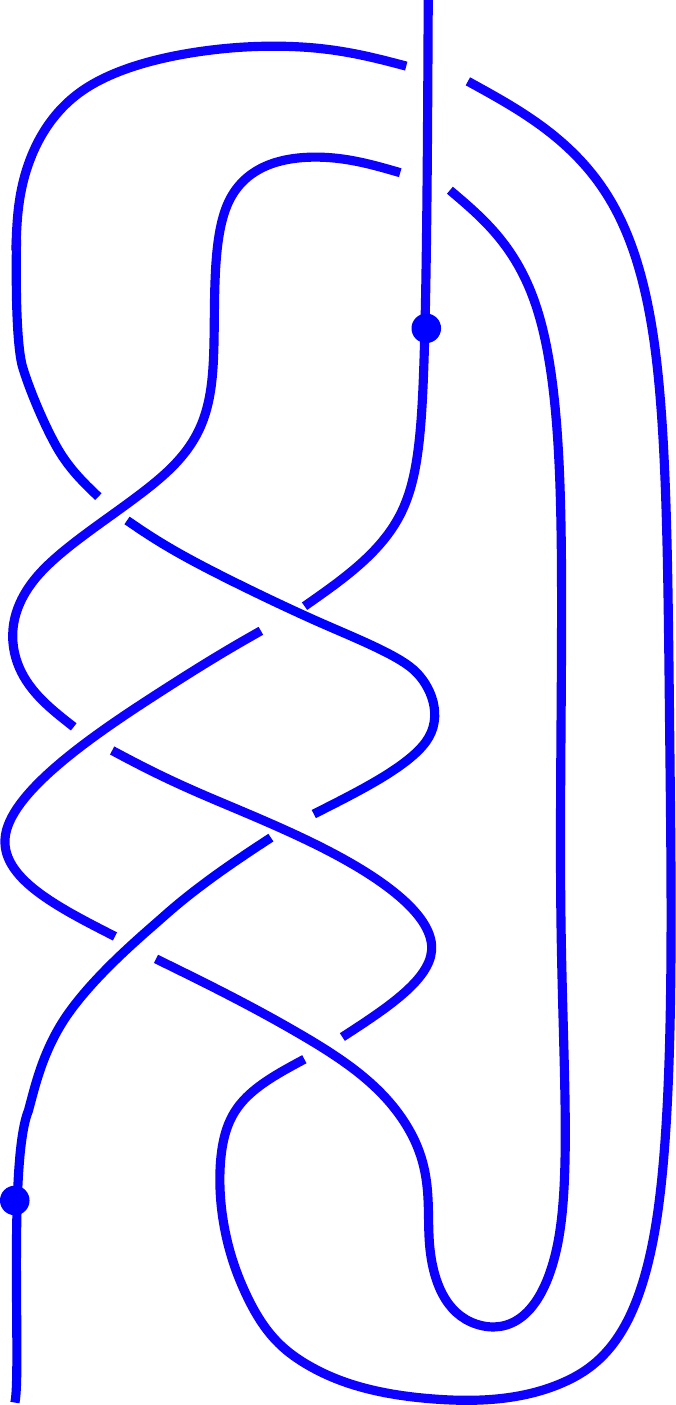}
\caption{Picture of the long Borromean rings (i.e., a pure braid commutator) $[b_{21},b_{31}]$ and the results of joining its components in the classical setting where $p=1$ and $n=3$.}
\label{F:LongBorr}
\end{figure}

For (spherical) 1-dimensional links in $\R^3$, 
joining all three components of the Borromean rings yields a trefoil knot, and joining just two of the three components yields the Whitehead link.  
(Figure \ref{F:LongBorr} shows long links which have the same Vassiliev invariants of order at most $2$ as long versions of these links.)
A generalization to isotopy classes of higher-dimensional spherical links is also known \cite{Skopenkov:HighCodimLinks}.
Theorem \ref{T:E} generalizes this fact to higher homotopy groups.
Certain classes in it can be viewed as analogues of the trefoil, Borromean rings, and Whitehead link.
Indeed, Theorem \ref{T:Borr} says that the ``parametrized long Borromean rings'' maps via graphing and closure to the high-dimensional spherical Borromean rings defined by Haefliger \cite[Section 4.1]{Haefliger:1962Annals}.  
(Using Theorem \ref{T:Borr} and Theorem \ref{T:E}, we also establish in Corollary \ref{C:HaefligerTrefoilEvenCodim} that the Haefliger trefoil generates $\pi_0\K_{2k-1}^{3k}$ for $k$ odd, the analogue of a fact proven by Haefliger for $k$ even.)

For equidimensional links, the domain of the graphing map $G^p:\mathcal{L}_{m\cdot 0}^{n-p}\to \mathcal{L}_{m\cdot p}^n$ is the configuration space $\Conf(m,\R^{n-p})$, denoted in this way in Theorem \ref{T:E}.
For $m=2$, we pre-compose by a homotopy equivalence $S^{n-p-1} \to \Conf(2,\R^{n-p})$.  
Theorem \ref{T:E} involves a map $J$ (defined in Definition \ref{D:Join}) that joins components.  
More specifically, it joins the last two components, and when $p=1$, component $m-1$ is traversed before component $m$.  Of course $J$ can be iterated.  
Theorem \ref{T:E} appears as Theorem \ref{T:Knots} in the main body.

\begin{mainthm}
\label{T:E}
Suppose $1 \leq p\leq n-3$ and $2n-3p-3\geq 0$.
\begin{itemize}
\item[(a)]
If $n-p$ is odd, then the map $\pi_{2n-2p-3} S^{n-p-1}\to \pi_{2n-3p-3}\K_p^n (\cong \Z)$ given by $p$-fold graphing followed by joining the two link components sends the Whitehead square $[\mathbbm{1}_{n-p-1}, \mathbbm{1}_{n-p-1}]$ of the identity to twice a generator.
Thus if $n-p=3$, $5$, or $9$, it sends the Hopf fibration to a generator. 
\item[(b)]
The map $\pi_{2n-3p-3}\Omega^p\Conf(3,\R^{n-p}) \to \pi_{2n-3p-3}\K_p^n (\cong \Z$ or $\Z/2$) induced by the composition of $p$-fold graphing followed by joining the three components together maps the ``parametrized long Borromean rings'' class $[b_{21}, b_{31}]$ to a generator.
\item[(c)] 
For $p=1$, there is a basis for $\pi_{2n-6}\LL{1}{n}$ modulo torsion, 
consisting of the images of a generator of $\pi_{2n-6}\K_1^n$ under the inclusions $\iota_1,\iota_2: \K_1^n \to \LL{1}{n}$;
the result of graphing and then joining two components of $[b_{21},b_{31}]$; and
for $n-p$ odd, the result of graphing $[\mathbbm{1}_{n-2}, \mathbbm{1}_{n-2}]$.

If $p\geq 2$, then $\pi_{2n-3p-3}\LL{p}{n}$ is generated by the two inclusions of a generator of $\pi_{2n-3p-3}\K_p^n$; the result of graphing and then joining two components of $[b_{21},b_{31}]$; and the image of a generating set of $\pi_{2n-2p-3}S^{n-p-1}$ under graphing.

If $p\geq 3$ and $m\geq 3$, then $\pi_{2n-3p-3}\mathcal{L}_{m\cdot p}^{n}$ is generated by the $m$ inclusions of a generator of $\pi_{2n-3p-3}\K_p^n$; the result of graphing and then joining two components of $[b_{21},b_{31}]$ for every pair of components $(i,j)$ with $1\leq i < j \leq m$; the image under graphing of a generating set of $\pi_{2n-2p-3}S^{n-p-1}$ for every $(i,j)$ with $1\leq i < j \leq m$; and the result of graphing $[b_{21},b_{31}]$ for every $(i,j,k)$ with $1 \leq i < j < k \leq m$.
\end{itemize}
\end{mainthm}

Recall that the first nontrivial homotopy group of $\K_p^n$ occurs in dimension $2n-3p-3$.
This is also the lowest dimension in which a class not detected by the first two stages of the Goodwillie--Weiss Taylor tower can appear.  Equivalently, it is the lowest dimension of classes dual to analogues of type-2 Vassiliev invariants.  
Remark \ref{R:ClassicalAnalogue} is a detailed discussion of this analogy in terms of the construction used in the proof.
The branching of part (c) into several cases according to the value of $p$ is not ideal, and 
Conjecture \ref{C:ExtendThmF} is that the statement for $p\geq 3$ extends to $p\geq 1$.
By recent results on rational homotopy of these spaces \cite{Fresse-Turchin-Willwacher:Emb}, the only hurdle is ruling out any other torsion.

Theorem \ref{T:E} is proven via cohomology classes from configuration space integrals for long $1$-dimensional links in $\R^n$ \cite{Cattaneo:2002, KMV:2013} and dual homology classes from resolving singular links.  
These integrals are indexed by graphs as in formulas \eqref{Eq:KappaCocycle} through \eqref{Eq:LkNum3}, and the singular links we use are illustrated in Figures \ref{F:Singular3Braid}, \ref{F:Singular2Braid}, \ref{F:Singular2Link}, \ref{F:SingularHopfBraid}, \ref{F:KnotFromHopf}, and \ref{F:SingularKnot}.  
We represent all the homotopy classes in the Theorem by such resolutions and their images under graphing.
In previous joint work of ours, we related the restriction to $1$-dimensional pure braids in $\R^n$ of these integrals to Chen's integrals for based loop spaces \cite{KKV:2020}.  We then determined the values of Chen integrals on iterated pure braid commutators (i.e.~Whitehead products) \cite{KKV:Primitive}.  
These two results of ours help us identify $G_*[b_{21},b_{31}]$ and $G_*[\mathbbm{1}_{n-p-1}, \mathbbm{1}_{n-p-1}]$ with resolved singular links.
We also rely on a calculation of $\pi_{2n-6} \LL{1}{n} \otimes \Q$ \cite{Arone-Turchin:Htpy, Songhafouo-Turchin:HHA} to understand 2-component link classes that come from neither braiding nor knotting.
We  use Budney's work \cite{Budney:Family} to obtain a generator of $\pi_{2n-6}\K_1^n$ and importantly to bootstrap from 
$1$-dimensional links to $p$-dimensional links.

Theorem \ref{T:E} suggests the possibility of obtaining arbitrary classes of long links from pure braids, generalizing known results on obtaining isotopy classes and $n$-equivalence classes of knots via pure braids; see Questions \ref{Q:PB} and \ref{Q:CablesSplitLinks}.

\medskip

The paper is organized as follows.
In Section \ref{S:Definitions}, we define various spaces of links and various maps between them.
Readers familiar with the subject matter may wish to initially skip this Section and refer back to it as needed.
The remaining Sections are mostly independent of each other, except that results from previous Sections are used, as noted.
In Section \ref{S:Injectivity}, we prove Theorem \ref{T:Lambda} (Theorem \ref{T:A}) about the injectivity of homotopy groups of spheres into homotopy groups of spaces of long 2-component links.
Section \ref{S:Bijectivity} contains Theorem \ref{T:pi_iLinksIso} (Theorem \ref{T:B}) and Theorem \ref{T:ImagesOfiS} (Theorem \ref{T:C}), which determine homotopy groups of spaces of long links, up to knotting, in a certain range.  Their proofs use Theorem \ref{T:A}.
In Section \ref{S:Knots}, we prove Theorem \ref{T:Knots} (Theorem \ref{T:E}) on nontrivial classes in spaces of long knots and links from classes in spaces of pure braids, including classes from spheres.  
Its proof uses Theorems \ref{T:A}, \ref{T:C}, and \ref{T:D}.
We conclude this Section with 
conjectures and questions to explore in future work.
In Appendix \ref{S:Appendix}, we prove the injectivity of graphing for spaces of 2-component spherical links and 2-component link maps.

\subsection{Acknowledgments}
This work was supported by the Louisiana Board of Regents Support Fund, contract number LEQSF(2019-22)-RD-A-22 and by the National Science Foundation, Award No.~DMS-2405370.
We thank Mark Behrens, Ryan Budney, Neeti Gauniyal, Danica Kosanovi\'c, and Victor Turchin for useful conversations.  We thank the referee for a careful reading and helpful comments.

\section{Key definitions}
\label{S:Definitions}

In Section \ref{S:SpacesDef}, we define a handful of spaces of smooth embeddings or smooth maps.  
In Sections \ref{S:MapsDef}
and \ref{S:RhoIota}, we define a number maps involving these spaces.  We begin with some basic conventions and notation.

For any $k\geq 0$, let $D^k$ denote the closed $k$-dimensional unit disk in $\R^k$.
Let $I:=D^1=[-1,1]$.
For a smooth manifold $X$ with a basepoint $x$, we write $\Omega^k X$ for the space of smooth based $k$-fold loops in $X$.  We take these to be smooth maps $\R^k \to X$ which are constant at $x$ outside $I^k$.  
We may sometimes allow the role of $I^k$ to be played by $D^k$, or by a product of disks whose dimensions sum to $k$, via a homeomorphism that is a diffeomorphism between the interiors.  
One can define $\pi_k X$ as $\pi_0 \Omega^k X$, and  $\pi_i \Omega^j X \cong \pi_{i+j} X$ for any non-negative integers $i$ and $j$.
We may also represent a $k$-fold loop in $X$ by a based map $S^k\to X$.
We write $\pi^s_i$ for the $i$-th stable homotopy group of spheres, i.e., $\pi^s_i:=\underset{k\to \infty}{\mathrm{colim}}\ \pi_{k+i} S^k$.

\subsection{Spaces of embeddings, link maps, and pseudoisotopy embeddings}
\label{S:SpacesDef}

For smooth manifolds $P$ and $N$, we write $\Emb(P,N)$ for the space of smooth embeddings of $P$ into $N$.  
Let $P_1, \dots, P_m$ be the path components of $P$.
A \textbf{link map} is a smooth map $(f_1, \dots, f_m): \coprod_1^m P_i\to N$ such that the images of the $f_i$ are pairwise disjoint.  We denote the space of such {link maps} by $\Link\left( \coprod_1^m P_i, N \right)$.
We will usually use the term {\bf links} for embeddings $\coprod_1^m P_i\to N$, especially when $m>1$, though we may occasionally use it for link maps when there is no risk of confusion.

If the components of $P$ have basepoints $b_1, \dots, b_m$ and $y_1, \dots, y_m$ are fixed distinct points in $N$, we define $\Emb_*(P,N)$ as the space of {\bf based embeddings}, meaning embeddings $f:P\to N$ such that $f(b_1)=y_1, \dots, f(b_m)=y_m$.
We write $\Emb_c(P,N)$ for a space of {compactly supported embeddings} of $P$ into $N$, which we will define more precisely in Definition \ref{D:LongLinks}.
We use the notation $\Link_*(P,N)$ and $\Link_c(P,N)$ for similarly defined spaces of based link maps and compactly supported link maps respectively.
For spaces $X$ and $Y$, we write $\Map_*(X,Y)$ for the space of based continuous maps $X \to Y$.

All spaces of smooth maps are equipped with the Whitney $C^\infty$-topology, while spaces of continuous maps have the compact-open topology.   
In this paper, $N$ will often be a sphere or Euclidean space, while $P$ will often be a disjoint union of finitely many spheres or Euclidean spaces.

\begin{definition}
\label{D:LongLinks}
Let $m$ and $n$ be positive integers and $p_1, \dots, p_m$ be nonnegative integers less than $n$.
Let $t_1^*, \dots, t_m^*$ be points with $-1< t_1^* < t_2^* < \dots < t_m^* <1$ and with equal gaps between each successive pair; that is, $t_i^*=-1+2i/(m+1)$.
We define the space $\Emb_c\left( \R^{p_1} \scu \dots \scu \R^{p_m}, \, \R^n \right)$ of \textbf{long links} in $\R^n$ as follows.
An element $f$ of this embedding space is required to map $\coprod_1^m (-1,1)^{p_i}$ into $(-1,1)^n$.
Outside of $\coprod_1^m (-1,1)^{p_i}$, $f$ must agree with the embedding 
\begin{align*}
e=(e_1, \dots, e_m): \coprod_1^m \R^{p_i} &\incl \R^n \\
e_{i}: (t_1, \dots, t_{p_i}) &\mapsto (t_i^*, 0, \dots, 0, t_1, \dots, t_{p_i}).
\end{align*}
We take $e$ to be the basepoint of $\Emb_c \left(\coprod_1^m \R^{p_i}, \, \R^n \right)$.  
A {\bf long knot} is a long link with one component.  
We define the space $\Link_c \left(\coprod_1^m \R^{p_i}, \, \R^n \right)$ of {\bf long link maps} similarly, just replacing embeddings by link maps.
\end{definition}

We sometimes abbreviate $\mathcal{L}_{p_1,\dots,p_m}^n:=\Emb_c \left(\coprod_1^m \R^{p_i}, \, \R^n \right)$ and write $\mathcal{L}_{m\cdot p}^n$ when $p_1=\dots=p_m=p$.  Many authors use the term ``string links'' instead of ``long links.''
Our convention of using the last rather than first $p_i$ coordinates of $\R^n$ in Definition \ref{D:LongLinks} is 
chosen to match
our conventions in Definition \ref{D:Graphing} of the graphing map.  Any other choice of standard long link $e$ produces a space of long links homeomorphic to the one in Definition \ref{D:LongLinks}.  Above each $e_i$ depends on $m$, which is suppressed from the notation when there is no risk of confusion.

If $p_1=\dots=p_m=0$ in Definition \ref{D:LongLinks}, one obtains the {\bf configuration space} of $m$ points in $\R^n$, which we also denote $\Conf(m, \R^n)$.  
Applying Definition \ref{D:LongLinks}, we obtain $((t_1^*,0, \dots,0), \dots, (t_m^*, 0,\dots,0))$ as the basepoint of $\Conf(m, \R^n)$.  
For a finite set $S$, we will also write $\Conf(S, \R^n)$ for $\Emb(S, \R^n)$.  
We call an element of the space $\Omega^p \Conf(m,\R^{n-p})$ a $p$-dimensional, $m$-component {\bf pure braid} in $\R^n$.  
Using the graphing map (see Definition \ref{D:Graphing} below), we may sometimes view such an element as lying in $\mathcal{L}_{m \cdot p}^n$.

For any $i$ and $j$ with $1 \leq i \neq j \leq m$, let $b_{ij} \in \pi_{n-1} \Conf(m, \R^n)$ be the cycle obtained from the map 
\begin{align*}
S^{n-1} &\to \Conf(m, \R^n)\\
v &\mapsto (x_1, \dots, x_{i-1}, x_j + \eps v, x_{i+1}, \dots, x_m)
\end{align*}
where $\eps< \min_i (t_i^* - t_{i-1}^*)$.  More precisely, we get a based map by conjugating by a path which takes the $i$-th point from its position at the basepoint to the image of the sphere above.  For $m=2$, the map $b_{12}$ (or $b_{21}$) is a homotopy equivalence.

\begin{definition}
\label{D:PEDef}
If $m\geq 1$ and $1 \leq p_1, \dots, p_m < n$, 
we define the space $\mathrm{PEmb}\left( \coprod_{i=1}^m \R^{p_i}, \, \R^n \right)$ of \textbf{pseudoisotopy embeddings} as the subspace of embeddings $f=(f_1, \dots, f_m): \coprod_{i=1}^m \R^{p_i} \incl \R^n$ such that 
\begin{itemize}
\item 
$f$ agrees with the standard embedding $e=(e_1, \dots, e_m)$ outside of $\coprod_{i=1}^m I^{p_i-1} \x [-1, \infty)$ and
\item 
there is a long link 
$g=(g_1, \dots, g_m) \in \LL[p_1-1,\dots]{p_m-1}{n-1}$
such that if $t_{p_i} \in [1, \infty)$, then \\
$f_i(t_1, \dots, t_{p_i})=(g_i(t_1,\dots,t_{p_i-1}), t_{p_i})$.
\end{itemize}
We will often abbreviate this space by $\P_{p_1, \dots, p_m}^n:=\mathrm{PEmb}\left( \coprod_{i=1}^m \R^{p_i}, \, \R^n \right)$.
\end{definition}

Some authors write $\mathrm{PE}(I^{p-1}, I^{n-1})$ or $\mathrm{CE}(I^{p-1}, I^{n-1})$ to denote a similarly defined space of embeddings $I^p \incl I^n$.  This space is homeomorphic to $\P_{p}^{n}$ via restriction and extension maps.  
One can think of a pseudoisotopy embedding as an embedding that is fixed on the bottom and sides of the domain cubes and takes the top faces of the domain into the top face of the codomain.

\subsection{Stacking, graphing, closure, and joining maps}
\label{S:MapsDef}

We now define various maps between spaces of links, starting with a multiplication on spaces of long links.  
\begin{definition}
\label{D:Stacking}
Let $m\geq 1$ and $1 \leq p_1,\dots,p_m < n$.  
There is a product on $\Emb_c \left( \coprod_1^m \R^{p_i}, \, \R^n \right)$ called {\bf stacking} which sends $(f,g) \in \left(\Emb_c \left( \coprod_1^m \R^{p_i}, \, \R^n \right)\right)^2$ to the map $f \# g\in \Emb_c \left( \coprod_1^m \R^{p_i}, \, \R^n \right)$ defined on $\coprod_1^m [-1,1]^{p_i}$ by 
\[
(f\# g)_i(t_1, \dots, t_{p_i}) := 
\left\{
\begin{array}{ll}
(\mathbbm{1}^{n-1} \x \rho_-^{-1}) \circ f_i \circ (\mathbbm{1}^{p_i-1} \x \rho_-)
 & 
 \text{ on } [-1,1]^{p_i-1} \x [-1,0] \\
(\mathbbm{1}^{n-1} \x \rho_+^{-1}) \circ f_i \circ (\mathbbm{1}^{p_i-1} \x \rho_+)
&
 \text{ on } [-1,1]^{p_i-1} \x [0,1]
\end{array}
\right.
\]
 for each $i=1,\dots,m$, where $\rho_{\pm}(t):=2t\mp1$ and where $\mathbbm{1}^k$ denotes the identity map on $\R^k$.  
 \end{definition}
This map gives the space of long links the structure of a homotopy-associative H-space.  
An analogous operation gives a homotopy associative H-space structure on $\Link_c \left( \coprod_1^m \R^{p_i}, \, \R^n\right)$.  If all the codimensions are at least 3, then $\pi_0$ of either of these spaces is an abelian group under stacking \cite{Haefliger:1962Top}.

\smallskip

The next map appears in most of our main results.  
In contrast with the usual conventions for graphs, we put the range coordinates of $f$ before its domain coordinates because for $n=2$ and $p_1=\dots=p_m=0$, we view braids as flowing vertically rather than horizontally.

\begin{definition}
\label{D:Graphing}
For any $m\geq 1$, $0 \leq p_1,\dots, p_m < n$, define the \textbf{graphing map}
\begin{align*}
G: \Omega \Emb_c \left( \coprod_1^m \R^{p_i}, \, \R^n \right) 
& \to \Emb_c \left( \coprod_1^m \R^{p_i+1}, \, \R^{n+1}\right) \\
(t \mapsto f(t)=(f_1(t), \dots, f_m(t))) &\mapsto (G(f)_1, \dots, G(f)_m)
\end{align*}
by 
\[
G(f)_i(t_1, \dots, t_{p_i+1}):=(f_i(t_{p_i+1})(t_1, \dots, t_{p_i}), \, t_{p_i+1})
\]
for $i=1,\dots, m$.  
\end{definition}

Iterating such graphing maps gives rise to maps of the form
\[
G^k: \Omega^k \Emb_c \left( \coprod_1^m \R^{p_i}, \, \R^n \right) 
 \to \Emb_c \left( \coprod_1^m \R^{p_i+k}, \, \R^{n+k}\right).
\]
More precisely, $G^k$ is the composition of the maps induced on various iterated loop spaces by various graphing maps, all of which by abuse of notation we will denote $G$.

\smallskip

The next map lets us pass from long links to spherical (based) links in a Euclidean space.
It will be used to relate isotopy classes of the two types in Lemma \ref{L:SphericalLinks} and to prove variants of the injectivity of graphing in Appendix \ref{S:Appendix}.
It is roughly given by fixing disjoint disks in the complement of $[-1,1]^n$ and gluing these ``closing disks'' to the parts of the long link in $[-1,1]^n$.

First, for any $k\geq 1$ and any subspace $S\subset \R^k$, let $rS$ be the result of scaling $S$ by $r$.
In particular, $[-1,1]^k$ is a cube inscribed in $\sqrt{k}S^{k-1}$.

\begin{definition}
\label{D:Closure}
For $m\geq 1$ and $1 \leq p_1, \dots, p_m \leq n-2$, define the {\bf closure map} 
\begin{align}
\label{Eq:Closure}
\begin{split}
\widehat{\cdot} \ : \Emb_c \left(\coprod_1^m \R^{p_i}, \, \R^n\right)
&\to \Emb_* \left(\coprod_1^m S^{p_i}, \, \R^n\right) \\
f=(f_1, \dots, f_m) &\mapsto \widehat{f}=(\widehat{f}_1, \dots, \widehat{f}_m)
\end{split}
\end{align}
where for each $i=1,\dots, m$, $\widehat{f}_i$ is given by 
\begin{itemize}
\item 
fixing a homeomorphism $S^{p_i} \to  \sqrt{n}D^{p_i} \cup \sqrt{n}D^{p_i}$, 
\item 
taking the union of $f_i|_{\sqrt{n}D^{p_i}}$ with an embedding $g_i$ of another copy of $\sqrt{n}D^{p_i}$ whose image is the hemisphere 
\[
\left\{(t_i^*,0,\dots,0,t_{n-p_i},\dots,t_n): t_{n-p_i}^2+\dots+t_n^2=n\right\}, 
\]
\item 
and smoothing the resulting injection $S^{p_i}\incl \R^n$ in a fixed small neighborhood of the intersection of the two disks 
$\sqrt{n}D^{p_i}$ using a fixed partition of unity.
\end{itemize}
The codomain of the map $\ \widehat{\cdot} \ $ is the space of based embeddings where we take the basepoint $b_i$ in each copy of $S^{p_i}$ to be any point in the image of $g_i$, and we take the image of $b_i$ in $\R^n$ to be the base value $y_i$.  
\end{definition}

Each $\widehat{f}_i$ has image in $t_i^* \x 0^{n-p_i-2} \x \R^{p_i+1} \subset \R^n$, so the assumption that $n-p_i\geq 2$ for all $i$ guarantees that their images are disjoint.  
Moreover, the closure $\widehat{e}$ of the standard long link $e$ is a trivial link, meaning that its components bound disjoint $(p_i+1)$-dimensional disks in $\R^n$.
We take $\widehat{e}$  to be the basepoint of $\Emb_* \left(\coprod_{i=1}^m S^{p_i}, \, \R^n\right)$.  We can also use it as the basepoint of the space $\Emb \left(\coprod_{i=1}^m S^{p_i}, \, \R^n\right)$ of unbased embeddings.
With different choices of embeddings $g_i$, one could define a closure map where one allows $n-p_i=1$, but it cannot take a standard long link to a trivial link if $n-p_i=1$ for three or more values of $i$.

We will use similarly defined closure maps denoted by the same symbol:
\begin{align*}
\begin{split}
\widehat{\cdot} \ : \Link_c \left(\coprod_{i=1}^m \R^{p_i}, \, \R^n \right)
& \to \Link_* \left(\coprod_{i=1}^m S^{p_i}, \, \R^n \right),
\end{split}
\begin{split}
\widehat{\cdot} \ : \Emb_c \left(\coprod_1^m \R^{p_i}, \, \R^n \right)
& \to \Emb \left(\coprod_1^m S^{p_i}, \, S^n \right) 
\end{split}
\end{align*}
The second one is obtained by post-composing the map \eqref{Eq:Closure} by the map induced by a fixed inclusion $\R^n \incl S^n$, 
and forgetting that the resulting embeddings are based.

\smallskip

The upcoming definition of the joining map is lengthy, but the payoff is the facilitation of compatibility with graphing.  
That compatibility is crucial in proving Theorem \ref{T:E}.
The idea is indicated in Figure \ref{F:LongBorr}, though without the modifications needed to obtain a based map of spaces of links. 
Roughly, the joining map closes the $m$-th component and then connect sums the $(m-1)$-th component, which could be done with a tube $S^{p-1} \x D^1$.  
Instead, we will essentially use a tubular neighborhood of a path in $\R^n$, and apply a flow of $\R^n$ along this path. 
The tube implicitly lies in the boundary of the tubular neighborhood.

We begin the precise definition by fixing auxiliary data for each $m$ and each codimension $n-p \geq 2$.
Recall that $t_1^*, \dots, t_m^*$ are the first coordinates of the components $e_1, \dots, e_m$ of the standard link $e$.
\begin{itemize}
\item 
Fix $\eps<1/m$.  Find a real number $k> 2$ and a smooth embedding $\gamma: [0,1] \to \R^{n-p+1}$ such that 
\begin{itemize}
\item 
$\gamma(0)=P:=(t_m^*, 0,\dots,0)$; 
$\gamma(1)=Q:=(t_{m-1}^*, 0, \dots, 0, k)$;  $\gamma'(0), \gamma'(1) \perp 0^{n-p} \x \R$; and
\item 
$\gamma$ has a tubular neighborhood $\eta$ of radius $1+\eps$ 
whose interior is contained in
\[
[-1,1] \x [-2k,2k]^{n-p} - \left( [-1,t_{m-1}^*+\eps] \x [-1,1]^{n-p} \cup e\left(\coprod_1^{m-1} \R^p\right) \right).
\]
\end{itemize}
\item
Fix a framing $\eta \xrightarrow{\cong} \gamma \x \R^{n-p}$ which restricts to orientation-preserving affine-linear maps
\begin{align*} 
{(t_m^*,0,\dots,0) \x [-1,1]} &\to P \x 0^{n-p-1} \x [-1,1] \text{ and} \\
 (t_{m-1}^*, 0,\dots,0)\x[k-1,k+1] &\to Q \x 0^{n-p-1} \x [-1,1].
 \end{align*}
 This is possible because up to homotopy, such a framing is given by a path in $O(n-p)$ with prescribed last vectors at the endpoints, and we assume $n-p\geq 2$.  
 \item 
 Let $H \subset \eta$ be the preimage of $\gamma \x 0^{n-p-1} \x [-1,1]$ under the framing.
\item 
Fix a diffeotopy $\phi_s: \R^{n-p+1} \to \R^{n-p+1}$, $s\in[0,1]$, such that 
\begin{itemize}
\item 
$\phi_0=\mathbb{1}^{n-p+1}$, 
\item 
$\phi_1(P)=Q$, and 
for each $s>0$, $\phi_s|_H$ is defined by flowing along $\gamma$, using the framing of $\eta$, 
\item 
for each $s\in [0,1]$, $\phi_s$ is supported in the set of points with distance at most $\eps$ from $H$, and 
\item 
for each $s\in [0,1]$ and points $x$ within $\eps$ of $H$, $\phi_s(x)$ is given by an interpolation (via a partition of unity)  between flowing along $\gamma$ and the identity map.  
\end{itemize}
See \cite[Chapter 8.1]{Hirsch:1976} for constructions of maps similar to but more general than $\phi_s$, which flows along an {embedded} path.
\end{itemize}

Because of the condition on the framing, $\phi_1$ takes 
${(t_m^*,0,\dots,0) \x [-1,1]}$ onto $(t_{m-1}^*, 0,\dots,0)\x[k-1,k+1]$ 
by the unique orientation-preserving affine linear map between them.  
Since $\phi_s$ flows along $\gamma$ and 
since $\mathrm{int}(\eta) \cap e(\coprod_1^{m-1} \R^p) = \varnothing$, 
we have $\phi_1((t_{m-1}^*, 1] \x [-1,1]^{n-1}) \cap e(\coprod_1^{m-1} \R^p) = \phi_1( (t_m^*,0,\dots,0) \x [-1,1]^p)$.

Recall that $\mathcal{L}_{m\cdot p}^n$ stands for  
$\mathcal{L}_{p_1,\dots, p_m}^n = \Emb_c\left(\coprod_{i=1}^m \R^{p_i} \, \R^n \right)$ 
where each $p_i=p$.  
To disambiguate between standard links of different numbers of components, we denote the $m$-component standard link by $e^m=(e^m_1, \dots, e^m_m)$ and the first coordinates of $e^m_1, \dots, e^m_m$ by $(t^m_1)^*, \dots, (t^m_m)^*$ for the rest of this Section.

\begin{definition}
\label{D:Join}
For $m\geq 2$ and $1\leq p \leq n-2$, we define the {\bf joining map}
\[
J: \mathcal{L}_{m\cdot p}^n \to \mathcal{L}_{(m-1)\cdot p}^n
\]
on a link $f=(f_1,\dots,f_m)$ in $\mathcal{L}_{m\cdot p}^n$ in two steps, using the diffeotopy $\phi_s$ and the number $k$ fixed above:
\begin{enumerate}
\item 
First define an embedding $\coprod_1^{m-1} \R^p \to \R^n$ as 
\[
(\phi_1 \x \mathbb{1}^{p-1}) \circ f  \ \ \cup \ \ (\phi_1 \x \mathbb{1}^{p-1}) \circ f \circ T \ \  \cup \ \ e^m
\]
on $A \cup B \cup C$, where 
\[
A:=\coprod_1^{m-1}[-1,1]^p , \quad  
B:=\coprod_1^{m-2} \varnothing \ \sqcup \ \left([k-1,k+1]\x[-1,1]^{p-1}\right),
\]
\[
C:=\coprod_1^{m-2} (\R^p - [-1,1]^p)\sqcup (\R^p - ([-1,1]^p \cup [k-1,k+1]\x[-1,1]^{p-1})),
\]
and $T$ is the affine-linear map taking $B$ onto $[-1,1]^p$ in the $m$-th component.
The maps and all their derivatives agree on the intersections because of the behavior of long link components outside $[-1,1]^p$ and because of the behavior of 
$\phi_1$ on $(t_m^*,0,\dots,0) \x [-1,1]^p$.
The above-mentioned property of $\phi_1((t_{m-1}^*, 1] \x [-1,1]^{n-1})$ guarantees that we get an embedding.
\item
Pre-compose the embedding $\coprod_1^{m-1} \R^p \to \R^n$ from step (1) by the map $2k \x \mathbb{1}^{p-1}: \R^p \to \R^p$ in each summand.
Post-compose it by the map $(S, 1/(2k), \dots, 1/(2k)) \x \mathbb{1}^{p-1}:\R^n \to \R^n$, where $S: \R \to \R$ is the affine-linear map that sends $(t^m_1)^*, \dots, (t^m_{m-1})^*$ to $(t^{m-1}_1)^*, \dots, (t^{m-1}_{m-1})^*$.
\end{enumerate}
\end{definition}

Step (2) ensures that outside of $\coprod_1^{m-1} [-1,1]^p$, $J(f)$ agrees with $e^{m-1}$, so $J(f)$ is a long link.
Regarding basepoints, applying step (1) to $f=e^m$ yields a link where every component is affine-linear and has derivative $(0^{n-p}, \mathbb{1}^p): \R^p \to \R^n$.  The effect of step (2) when $f=e^m$ is only to alter the first coordinate in each component, thus yielding $e^{m-1}$.  Hence $J$ preserves basepoints, i.e., $J((e^m_1, \dots, e^m_m))=(e^{m-1}_1, \dots, e^{m-1}_{m-1})$.  
As with graphing, we denote by abuse of notation all of the joining maps for various $m$, $n$, and $p$ by the same symbol $J$.
A generalization of $J$ to joining different pairs of components is possible but not needed for our purposes.  We now show that graphing commutes with joining components:

\begin{proposition}
\label{P:JG=GJ}
Suppose $m\geq 2$ and $1 \leq p \leq n-2$.
The composites $J \circ G$ and $G\circ \Omega J$ are homotopic as maps 
$\Omega \mathcal{L}_{m\cdot p}^n \to \mathcal{L}_{(m-1)\cdot (p+1)}^{n+1}$.
\end{proposition}

\begin{proof}
The key point is that $J$ is defined by maps that use only the first coordinate of the domain $\coprod_1^m \R^p$ and the first $n-p+1$ coordinates of the codomain $\R^n$.
Let $\omega \in \Omega \mathcal{L}_{m\cdot p}^n$.  
Applying either of the two maps to $\omega$ gives an embedding $\coprod_1^{m-1} \R^{p+1} \to \R^{n+1}$.   
We view $\omega$ as a family of embeddings $\omega_t:\coprod_1^m \R^p \x t \to \R^n \x t$, $t\in[-1,1]$, and hence $\Omega J (\omega)$ as a family of embeddings $\coprod_1^{m-1} \R^p \x t \to \R^n \x t$, $t\in [-1,1]$.
In these terms, $G \circ \Omega J (\omega): \coprod_1^{m-1} \R^{p+1} \to \R^{n+1}$ is obtained from $\Omega J (\omega)$ by first taking the union of the latter family over both the domain and codomain and then extending by the standard link $(e_1,\dots,e_{m-1})$ for $t_{p+1} \notin [-1,1]$.
Thus we must essentially compare the union of the family $\Omega J(\omega)$ with the embedding $\coprod_1^{m-1} \R^p\x[-1,1] \to \R^n \x [-1,1]$ obtained by restricting the domain and codomain of $J \circ G(\omega)$.

Let $A$, $B$, and $C$ be the subspaces in the definition of joining $p$-dimensional links in $\R^n$.  Then each of $\Omega J (\omega)$ and $J \circ G(\omega)$  is obtained by gluing embeddings of $A \x [-1,1]$, $B \x [-1,1]$, and $C \x [-1,1]$.  On $C \x [-1,1]$, both are given by the restriction of the standard $(p+1)$-dimensional $(m-1)$-component long link, since the combined effect of the pre- and post-compositions by the affine-linear maps of the domain and codomain leaves the standard long link fixed.  
Recall that $\phi_1$ is a diffeomorphism of $\R^{n-p+1}$.
On $A\x[-1,1]$, $J \circ G(\omega)$ is given by 
\begin{equation}
\label{Eq:JG}
((S, 1/(2k), \dots, 1/(2k)) \x \mathbb{1}^{p}) \  \circ \  ((\phi_1 \x \mathbb{1}^{p}) \circ G(\omega))  \ \circ \ (2k \x \mathbb{1}^{p}).
\end{equation}
On each slice $A\x t$, $\Omega J (\omega)$ is given by 
\begin{equation}
\label{Eq:GJ}
((S, 1/(2k), \dots, 1/(2k)) \x \mathbb{1}^{p-1}) \  \circ \ ((\phi_1 \x \mathbb{1}^{p-1}) \circ \omega_t) \ \circ \ (2k \x \mathbb{1}^{p-1}).
\end{equation}
Applying $G$ to \eqref{Eq:GJ} means taking the union over $t\in [-1,1]$, which has the effect of replacing $\omega_t$ by $G(\omega)$ and replacing each instance of $\mathbb{1}^{p-1}$ by $\mathbb{1}^p$, thus yielding \eqref{Eq:JG}.  
On $B\x[-1,1]$, we have a similar comparison of two expressions, except that the second factors are $(\phi_1 \x \mathbb{1}^{p}) \circ G(\omega)\circ T$ and $(\phi_1 \x \mathbb{1}^{p-1}) \circ \omega_t \circ T$ respectively.
\end{proof}

We have actually shown that $J\circ G = G \circ \Omega J$, but in Section \ref{S:KnotsThm} we will consider an alternative description of $J$ that may agree only up to homotopy with Definition \ref{D:Join}.

\begin{remark} 
\label{R:ClassesOfJoin}
One can show that $J$ depends only on the isotopy $h_s \x \mathbb{1}^{p-1}$ of the $m$-th copy of $[-1,1]^p$ into the $(m-1)$-th component, where $h_s$ is the restriction of $\phi_s$ to $(t^m_m)^*\x 0^{n-p-1} \x [-1,1]$.  
If we generalize the construction from $h_s \x \mathbb{1}^{p-1}$ to any isotopy of $[-1,1]^p$ in $\R^n$ coming from an embedded path $\gamma$ in $\R^n$, then the space of possible choices in this construction is homotopy equivalent to 
\[
(\R^p - [-1,1]^p)) \ltimes 
\left( \Emb_c\left(\R, \R^n - \left([-1,1]^n \cup e\left(\coprod_1^{m-1} \R^p\right)\right)\right) \x \Omega V_{p}(\R^{n-1}) \right)
\]
 where $B \ltimes F$ is imprecise notation for a bundle with base $B$ and fiber $F$, and where $V_p(\R^{n-1})$ is the Stiefel manifold of $p$-frames in $\R^{n-1}$.  
If we restrict to paths $\gamma$ lying in $\R^{n-p+1}\x 0^{p-1}$ and isotopies $h_s \x \mathbb{1}^{p-1}$ (as we do to ensure compatibility with graphing), then the space of possible choices is as above but with every instance of $n$ replaced by $n-p+1$ and every instance of $p$ replaced by $1$.  In either setting, the possible homotopy classes of $J$ depend only on $\pi_0$ of this space.
\end{remark}

\subsection{Restriction and inclusion maps}
\label{S:RhoIota}

We now define restriction and inclusion maps between spaces of links of different numbers of components.  
Let $S \subset \{1,\dots,m\}$.
Our link components are always ordered, and accordingly we view $S$ as an ordered set using the standard order on $\{1,\dots,m\}$.
Let $t_i^*\in (-1,1)$, $i\in\{1,\dots, m\}$ and $s_i^*\in (-1,1)$, $i\in S$ be the basepoints' first coordinates for $\mathcal{L}_{p_1, \dots, p_m}^n$ and $\mathcal{L}_{p_i : i\in S}^n$ respectively. 
Fix a diffeomorphism 
$g:\R^2\to\R^2$ isotopic to the identity 
such that 
$g([-1,1]^2) \subset [-1,1]^2$ and
$g(t_i^*,0)=(s_i^*,0)$ for all $i\in S$.
Fix a diffeomorphism 
$h:\R^2\to\R^2$ isotopic to the identity 
such that 
$h([-1,1]^2) \subset [-1,1]^2$,
$h(s_i^*,0)=(t_i^*,0)$ for all $i\in S$, 
and $(t_j^*,0)\notin h([-1,1]^2)$ for all $j\notin S$. 
If $S$ is a consecutive subset of $\{1,\dots, m\}$, then each of $g$ and $h$ can be taken to be the product of an affine-linear map $\R \to \R$ with $\mathbb{1}: \R \to \R$.  
It is convenient now to extend Definition \ref{D:LongLinks} to allow $m=0$, i.e.~an empty list $(\ )$ of source dimensions, in which case $\mathcal{L}^n_{(\ )}:=\{\ast\}$ (where the point may be viewed as the embedding of the empty set).

\begin{definition}
\label{D:Restriction}
Let $m\geq 1$ and $1 \leq p_1, \dots, p_m \leq n-1$.
Let $S \subset \{1,\dots, m\}$.
 Define the {\bf restriction} $\rho_S$ as the map 
\begin{equation}
\label{Eq:LinksToKnotsFibn}
\Emb_c\left( \coprod_{p_i : i \in \{1,\dots,m\} - S} \R^{p_i}, \, \R^n - \coprod_{i \in S} \R^{p_i} \right) \to \mathcal{L}_{p_1, \dots, p_m}^n \overset{\rho_S}{\longrightarrow} 
\mathcal{L}_{p_i: i \in S}^n
\end{equation}
which first restricts a link to the components indexed by $S$ and then applies $g \x \mathbb{1}^{n-2}$, where the fiber is taken over the standard long link $e=(e_i)_{i \in S}$.  
Above $\R^n - \coprod_{i\in S} \R^{p_i}$ is shorthand for $\R^n - \coprod_{i\in S} e_i(\R^{p_i})$.
\end{definition}
The map $\rho_S$ is a fibration by work of Palais or Lima \cite{Palais:LocalTriv, Lima:1964}.

\begin{definition}
\label{D:Inclusion}
Let $m\geq 1$ and $1 \leq p_1, \dots, p_m \leq n-2$.
Let $S \subset\{1,\dots, m\}$.
Define the {\bf inclusion} 
\[
\iota_S: \mathcal{L}_{p_i : i \in S}^n \to \mathcal{L}_{p_1, \dots, p_m}^n
\]
by setting $\iota_S(f)$ to be the standard embedding $(e_i)_{i \in \{1,\dots,m\} - S}$ together with 
$(h\x \mathbb{1}^{n-2})\circ f$.
\end{definition}

We continue to suppress the ambient set containing $S$ from the notation for $\rho_S$ and $\iota_S$.  Though this set varies below between $\{1,\dots,m\}$ and subsets thereof, it should be clear from the context.

\begin{proposition}
\label{P:iota-section-up-to-htpy}
Let $m\geq 1$ and $1 \leq p_1, \dots, p_m \leq n-2$.
Let $S,T \subset\{1,\dots, m\}$. 
\begin{itemize}
\item[(a)]
The map $\iota_S$ is a section of $\rho_S$ up to homotopy.
\item[(b)]
The composite $\rho_T \circ \iota_S$ is 
is homotopic to the composite
$\mathcal{L}_{p_i : i\in S}^n 
\xrightarrow{\rho_{S\cap T}}
\mathcal{L}_{p_i : i\in S \cap T}^n 
\xrightarrow{\iota_{S\cap T}} 
\mathcal{L}_{p_i : i\in T}^n$.
 Thus the induced maps on homotopy groups satisfy $\im(\rho_T \circ \iota_S)_* \subset \im(\iota_{S\cap T})_*$.
\item[(c)]
If $T\subset S$, the composite of the inclusions associated to $T\subset S$ and $S \subset \{1,\dots,m\}$ is homotopic to the inclusion associated to $T\subset \{1,\dots,m\}$.  An analogous statement holds for the restrictions.
 Thus $\im (\iota_T)_* \subset \im (\iota_S)_*$ and $\ker (\rho_S)_* \subset \ker (\rho_T)_*$.
\end{itemize}
\end{proposition}

\begin{proof}
For part (a),  $\rho_S \circ \iota_S$ takes a long link $f$ to $(g \circ h\x \mathbb{1}^{n-2})\circ f$.  By assumption $g \circ h:\R^2 \to \R^2$ is isotopic to $\mathbb{1}^2$, and an isotopy to $\mathbb{1}^2$ produces a homotopy from $\rho_S \circ \iota_S$ to the identity map on $\mathcal{L}_{p_i:i\in S}$.  
Similar homotopies yield the statements about maps in parts (b) and (c), from which the statements about the induced maps are immediate.
\end{proof}

\begin{corollary}
\label{C:IntKerPlusSumIm}
Let $m\geq 1$, $1 \leq p_1, \dots, p_m \leq n-3$, and $i \geq 0$.
Let $\mathcal{S}$ be any set of subsets of $\{1,\dots,m\}$.
Then 
\begin{equation}
\label{Eq:IntKerPlusSumIm}
\pi_i \mathcal{L}_{p_1, \dots, p_m}^n \cong
\bigcap_{S\in \mathcal{S}}
\ker (\rho_S)_* \oplus 
\sum_{S\in \mathcal{S}}
\im (\iota_S)_*.
\end{equation}
\end{corollary}

\begin{proof}
The left-hand side is an abelian group for all $i \geq 0$ because every $p_j\leq n-3$.
By Proposition \ref{P:iota-section-up-to-htpy} (a), for any $S\subset \{1,\dots, m\}$, $\iota_S$ gives a splitting of the long exact sequence in homotopy of the fibration $\rho_S$ and hence $\pi_i \mathcal{L}_{p_1, \dots, p_m}^n \cong \ker (\rho_S)_* \oplus \im (\iota_S)_*$.
Repeatedly applying this decomposition for every $S \in \mathcal{S}$ shows that $\pi_i \mathcal{L}_{p_1, \dots, p_m}^n $ is the sum of the intersection and the sum on the right-hand side.  One can show that this sum is direct by using induction on the cardinality of $\mathcal{S}$ and all three parts of Proposition \ref{P:iota-section-up-to-htpy}.
\end{proof}

\begin{corollary}
\label{C:2CompRhoIota}
For $m=2$, the restriction fibration 
\begin{equation}
\label{Eq:LinksToKnotsFibn2Comp}
\Emb_c(\R^p, \, \R^n - \R^q) 
\xrightarrow{\epsilon}
 \LL[p]{q}{n} \overset{\rho_2}{\longrightarrow} \K_{q}^n
\end{equation}
induces an isomorphism 
\begin{equation}
\label{Eq:pi_iLpqNDecomp1}
\pi_i \LL[p]{q}{n} \cong  \pi_i \Emb_c(\R^p, \, \R^n - \R^{q}) \oplus \pi_i \K_q^n
\end{equation}
where 
the inclusion of $\pi_i \Emb_c(\R^p, \, \R^n - \R^{q})$ is induced by the inclusion $\epsilon$ of the fiber of $\rho_2$ and 
the inclusion of $\pi_i \K_q^n$ is induced by the section $\iota_2$ of $\rho_2$.  In addition, 
\begin{equation}
\label{Eq:pi_iLpqNDecomp2}
\pi_i \LL[p]{q}{n} \cong (\ker (\rho_1)_* \cap \ker (\rho_2)_*) \oplus \pi_i \K_p^n \oplus \pi_i \K_q^n
\end{equation}
where the inclusions of the last two summands are induced by $\iota_1$ and $\iota_2$.
\end{corollary}

\begin{proof}
The decomposition \eqref{Eq:pi_iLpqNDecomp1}
follows from
 Corollary \ref{C:2CompRhoIota} with $\mathcal{S}=\{\{2\}\}$. 
We obtain the decomposition \eqref{Eq:pi_iLpqNDecomp2} from Corollary \ref{C:2CompRhoIota} with 
$\mathcal{S}=\{\{1\},\{2\}\}$, using the fact 
that $\im (\iota_1)_* \cap \im (\iota_2)_* =0$.
The latter fact holds because $(\rho_j\circ \iota_j)_*$ is the identity and $(\rho_j \circ \iota_k)_*=0$ if $j\neq k$, by Proposition \ref{P:iota-section-up-to-htpy} (a) and (b) respectively.
\end{proof}

\medskip

The graphing, restriction, and inclusion maps are maps of H-spaces, whereas the maps that join components are not.  
Graphing commutes with both restriction and inclusion.
Analogues of the stacking, joining, restriction, and inclusion maps for closed links appear in work of Haefliger \cite{Haefliger:1962Top}, where they are respectively called addition, contraction, projection, and inclusion.

\section{Injectivity of graphing for spaces of 2-component long links}
\label{S:Injectivity}

We now prove Theorem \ref{T:Lambda} (a.k.a.~Theorem \ref{T:A}),  the injectivity of graphing into spaces of 2-component links, which produces elements of homotopy groups in arbitrarily high degrees.  

Recall that there is a diffeomorphism 
\begin{align}
\label{Eq:ComplementDiffeo}
\begin{split}
h = h_q^n:\R^n - 0^{n-q} \x\R^q = (\R^{n-q} - 0^{n-q}) \x \R^q
&\to S^{n-q-1} \x \R^{q+1} \\
(t_1, \dots, t_n) &\mapsto \left( \frac{(t_1, \dots, t_{n-q})}{|(t_1, \dots, t_{n-q})|}, \ \ln|(t_1, \dots, t_{n-q})|, \ t_{n-q+1}, \dots, t_n \right).
\end{split}
 \end{align}
We can use it to define a homotopy equivalence 
\begin{align}
\label{Eq:ComplementHtpyEqv}
\begin{split}
S^{n-q-1} &\to \R^n - 0^{n-q}\x \R^q \\
s &\mapsto h^{-1}(s, 0^{q+1})
\end{split}
\end{align}
We can similarly define a based homotopy equivalence $S^{n-q-1} \to \R^n - e_2(\R^q)$, where $e_1(0)$ is the basepoint in the codomain,  by post-composing by the appropriate affine-linear map in the first coordinate.
By taking the one-point compactification of $\R^n$, we obtain from $h$ a diffeomorphism
\begin{equation}
\label{Eq:SpheresComplementDiffeo}
S^n - S^q \to S^{n-q-1} \x \R^{q+1}.
\end{equation}

We will now show there is a homotopy retraction from a space of long links to a wedge of spheres.  This next lemma can be viewed as an adaptation to long links of the $\lambda$-invariant given in Definition \ref{D:Lambda} below.

\begin{lemma}
\label{L:HtpyRetract}
Let $m\geq 2$ and $1 \leq p_1 \leq \dots \leq p_m \leq n-1$.
Then 
$\Omega^{p_1}\left(\bigvee_{k=2}^m S^{n-p_k-1}\right)$ 
is a retract up to homotopy of $\Emb_c\left(\R^{p_1}, \ \R^n -  \coprod_{k=2}^m \R^{p_k} \right)$, and   
a section of it is given by $\Omega^{p_1}j$ followed by the graphing map $G^{p_1}$, where 
$j:\bigvee_{k=2}^m S^{n-p_k-1} \to \R^n -  \coprod_{k=2}^m \R^{p_k}$ is a homotopy equivalence.
\end{lemma}

\begin{proof}
First, there is a map 
\begin{equation}
\label{Eq:EmbToMap}
r:\Emb_c\left(\R^{p_1}, \ \R^n -  \coprod_{k=2}^m \R^{p_k} \right) \to
\Omega^{p_1}\left(\R^{n-p_1} -  \coprod_{k=2}^m \R^{p_k-p_1} \right)
\end{equation}
given by viewing an embedding as a smooth map and projecting onto the first $n-p_1$ coordinates.  
As in Definition \ref{D:Restriction}, $\R^{p_k}$ and $\R^{p_k - p_1}$ are shorthand for their images under the standard embeddings $e_k$.
We claim that graphing gives a section of $r$.  Indeed, by a slight abuse of notation, consider the map 
\[
G^{p_1}:
\Omega^{p_1}\left(\R^{n-p_1} -  \coprod_{k=2}^m \R^{p_k-p_1} \right) \to 
\Emb_c\left(\R^{p_1}, \ \R^n -  \coprod_{k=2}^m \R^{p_k} \right)
\]
given by first viewing a point in $\Omega^{p_1}X$ as a $p_1$-fold loop in $\Emb_c(\{\ast\}, X)$ and then applying the $p_1$-fold graphing map.
Then $r \circ G^{p_1}$ is the identity:
\[
f \xmapsto{G} \Bigl( (t_1, \dots, t_{p_1}) \mapsto \bigl(f(t_{1}, \dots, t_{p_1}), t_{n-p_1+1}, \dots, t_{n} \bigr) \in \R^{n} \Bigr) 
\xmapsto{r} 
\bigl((t_1, \dots, t_{p_1}) \mapsto f(t_1, \dots, t_{p_1})\bigr).
\]
Thus 
$\Omega^{p_1}\left(\R^{n-p_1} -  \coprod_{k=2}^m \R^{p_k-p_1} \right)$ is a homotopy retract of 
$\Emb_c\left(\R^{p_1}, \ \R^n -  \coprod_{k=2}^m \R^{p_k} \right)$.  
A based homotopy equivalence 
$j:\bigvee_{k=2}^m S^{n-p_k-1} \xrightarrow{\simeq} \R^{n-p_1} -  \coprod_{k=2}^m \R^{p_k-p_1}$, where $e_1(0)$ is the basepoint in the codomain, can be obtained along similar lines to formula \eqref{Eq:ComplementHtpyEqv}, though an explicit formula is not as easily obtained for $m>2$ as for $m=2$.
Pre-composing $G^{p_1}$ by $\Omega^{p_1}(j)$ yields the desired section, since post-composing $r$ by the $p_1$-fold looping of the homotopy inverse of $j$ gives its one-sided inverse.
\end{proof}

\begin{theorem}
\label{T:Lambda}
 If  $\ 0 \leq p \leq q \leq n-1$ and $i\geq 0$, then $\pi_i \Emb_c(\R^p \scu \R^q, \, \R^n)$ contains a direct summand of $\pi_{i+p}S^{n-q-1}$.
An inclusion of it is given by composing maps induced by a based homotopy equivalence 
$j: S^{n-q-1} \to \R^{n-p} - e_2(\R^{q-p})$, 
the inclusion 
$\epsilon: \Emb_c(\ast, \R^{n-p} - e_2(\R^{q-p})) \incl \Emb_c(\ast \scu \R^{q-p}, \, \R^{n-p})$, 
and the $p$-fold graphing map $G^p$:
$\Omega^p \Emb_c(\ast \scu \R^{q-p}, \, \R^{n-p}) 
\to
 \Emb_c(\R^p \scu \R^q, \, \R^n)$.
\end{theorem}

Here the basepoint for $\R^{n-p} - e_2(\R^{n-q})$ is defined to be $(t_1^*,0,\dots,0)$.

\begin{proof}
There is a commutative square 
\begin{equation}
\label{Eq:FiberOfRestriction}
\xymatrix{
\Omega^p \Emb_c(\ast, \, \R^{n-p} - e_2(\R^{n-q})) \ar@{^(->}[r]^-\epsilon \ar[d]_-{G^p} & \Omega^p\Emb_c(\ast \scu \R^{q-p}, \, \R^{n-p})  \ar[d]^-{G^p} \\
\Emb_c(\R^p, \, \R^n - e_2(\R^q)) \ar@{^(->}[r]^-\epsilon & \Emb_c(\R^p \scu \R^q, \, \R^n)
}
\end{equation}
so we may consider $\epsilon \circ G^p \circ \Omega^p j$
instead of $G^p \circ \epsilon \circ \Omega^p j$.
By Corollary \ref{C:2CompRhoIota}, for each row of the square, the inclusion $\epsilon$ of the fiber of the restriction $\rho_2$ induces an inclusion of the homotopy groups of the left-hand side as a direct summand of those of the right-hand side.  The theorem follows from the fact that by Lemma \ref{L:HtpyRetract}, $G^p \circ \Omega^p j$
induces the inclusion of a direct summand on homotopy groups.
\end{proof}

\begin{example}[Generalized linking numbers]
\label{Ex:LkNumExamples}
One of the simplest cases of Theorem \ref{T:Lambda} and \ref{T:LambdaSpherical} is when $p=q=1$, $n=3$, and $i=0$, where the theorem reduces to the injection of $\pi_1(S^1) \cong \Z$, the group of 2-strand classical pure braids, into the monoid of isotopy classes of classical 2-component long, or closed, links.  
More generally, if $i+p+q=n-1$, then a generator of the resulting copy of $\Z$ in 
$\pi_{i}\Emb_c (\R^p \sqcup \R^q \, , \R^{i+p+q+1})$ is dual to the linking number of manifolds of dimensions say $i+p$ and $q$ in $\R^{i+p+q+1}$, by combining the domain of an element of $\pi_i$ and the source manifold $\R^p$ into an embedded $(i+p)$-dimensional manifold.
\end{example}

\section{Bijectivity of graphing in a range}
\label{S:Bijectivity}

We will now give a complete calculation, at least up to knotting, of homotopy groups of spaces of long links in a certain range.  
We will consider spaces of long (knots and) links, such as $\K_p^n$ and $\LL[p]{q}{n}$, as well as spaces of pseudoisotopy embeddings, such as $\mathcal{P}_p^n$ and $\PE[p]{q}{n}$.
The key ingredients, given in Section \ref{S:PELemmas}, are a fibration sequence involving spaces of long links and a space of pseudoisotopy embeddings, together with a lower bound on the connectivity of the latter space.  
In Section \ref{S:BijThm}, we prove Theorem \ref{T:pi_iLinksIso} (a.k.a.~Theorem \ref{T:B}), establishing bijectivity in a range for graphing of 2-component links. 
In this range, we will see in Theorem \ref{T:ImagesOfiS} (a.k.a.~Theorem \ref{T:C}) that this determines the homotopy groups of any space of long links.
The main result of Section \ref{S:SphericalLinks} is Theorem \ref{T:GraphingDim>1}, which extends the bijectivity result to one degree higher for equidimensional links by starting with 2-dimensional rather than 0-dimensional links.  This requires the identification in Lemma \ref{L:SphericalLinks} of isotopy classes between spherical links and long links.

\subsection{Lemmas on pseudoisotopy embedding spaces}
\label{S:PELemmas}

Restriction of a 2-component pseudoisotopy embedding $\R^p \scu \R^q \incl \R^n$ to the hyperplanes $\{(t_1,\dots,t_p)\in \R^p: t_p=1\} \scu \{(t_1, \dots, t_q)\in \R^q: t_q=1\}$ gives the following fibration:
\[
\xymatrix{
\LL[p]{q}{n} \ar[r]^-i & \PE[p]{q}{n} \ar[r]^-\rho & \LL[p-1]{q-1}{n-1}
}
\]
In turn, this leads to the sequence
\[
\xymatrix{
\Omega \LL[p-1]{q-1}{n-1}  \ar[r]^-\del & \LL[p]{q}{n} \ar[r]^i & \PE[p]{q}{n} 
}
\]
which up to homotopy is a fibration.
To be somewhat explicit, we review the general construction of the connecting map $\del$ for a fibration in this special case.  It comes from the following homotopy-commutative diagram:
\begin{equation}
\label{E:Hofiber}
\xymatrix{
\Omega \LL[p-1]{q-1}{n-1} \ar@{-->}[r]^-\del & \LL[p]{q}{n} \ar[r]^-i  & \PE[p]{q}{n} \ar@{=}[d]\\
\hofib(i) \ar[r]^-j \ar[u]^-p_-\simeq & \widetilde{\mathcal{L}}_{p,q}^n \ar[u]^-\pi_-{\simeq} \ar[r]^-{k} & \PE[p]{q}{n}
}
\end{equation}
The space $\widetilde{\mathcal{L}}_{p,q}^n$ is defined by 
\[
\widetilde{\mathcal{L}}_{p,q}^n := \left\{ (f, \gamma) : f \in \LL[p]{q}{n}, \ \gamma: [-1,1] \to \PE[p]{q}{n}, \ \gamma(1)=i(f) \right\}.
\]
The homotopy equivalence $\pi: \widetilde{\mathcal{L}}_{p,q}^n \to \LL[p]{q}{n}$ is the projection $\pi:(f, \gamma) \mapsto f$.
The map $k:\widetilde{\mathcal{L}}_{p,q}^n\to \PE[p]{q}{n}$ is given by $k:(f, \gamma) \mapsto \gamma(-1)$, and it is homotopic to $i \circ \pi$.
One defines $\hofib(i)$ as its fiber, i.e., 
\[
\hofib(i):= \left\{ (f, \gamma) : f \in \LL[p]{q}{n}, \ \gamma: [-1,1] \to \PE[p]{q}{n}, \ \gamma(-1)=e, \ \gamma(1)=i(f) \right\}.
\]
The homotopy equivalence $p:\hofib(i) \to \Omega \LL[p-1]{q-1}{n-1}$ is given by $(f, \gamma) \mapsto \rho \circ \gamma$, and for $h$ homotopy inverse to $p$, we can define $\del:= \pi \circ j \circ h$.

\begin{lemma}
\label{L:GraphingPEFibn}
The graphing map $G: \Omega \LL[p-1]{q-1}{n-1} \to \LL[p]{q}{n}$ agrees with $\del$.
\end{lemma}

\begin{proof}
We will specify a map $h: \Omega \LL[p-1]{q-1}{n-1} \to \hofib(i)$ homotopy inverse to $p$ so that with $\del=\pi \circ j \circ h$, we get $G = \del$.  
We thus need to construct out of a loop of long links a path of pseudoisotopy embeddings in one dimension higher, which starts at the standard embedding $e$ and ends at a long link.
Write a map $f \in \Omega \LL[p-1]{q-1}{n-1}$ as 
\begin{align*}
f=(f_1,f_2) : \R \x (\R^{p-1} \scu \R^{q-1}) &\to \R^{n-1}  \\ 
f_1:(s,(t_1, \dots, t_p)) & \mapsto f_1(s)(t_1, \dots, t_{p-1})\\
f_2:(s,(t_1, \dots, t_q)) & \mapsto f_2(s)(t_1, \dots, t_{q-1})
\end{align*}
by identifying $\R \x (\R^{p-1} \scu \R^{q-1})$ with $(\R \x \R^{p-1}) \scu (\R \x \R^{q-1})$.
We define $h$ by the formula 
\begin{align*}
h(f) = (h(f)_1, \, h(f)_2) : [-1,1] \x (\R^p \sqcup \R^q) &\to \R^n \\
h(f)_1 : (s, (t_1, \dots, t_p)) & \mapsto 
\left(
f_1(g(s,t_p))
(t_1, \dots, t_{p-1}), \, t_p \right) \\
h(f)_2 : (s, (t_1, \dots, t_q)) & \mapsto 
\left(
f_2(g(s,t_q))
(t_1, \dots, t_{q-1}), \, t_q \right) 
\end{align*}
where $g(s,t)=(s+1)\rho(t)-1$ and $\rho:\R\to \R$ is a smooth cutoff function that is $-1$ for $t\leq -1/2$ and $1$ for $t \geq 1/2$.  
We conclude the proof by noting that the function $g: \R^2\to \R$ satisfies the following properties:
\begin{itemize}
\item $g$ is a smooth function of $t$ for all $s\in [-1,1]$,
\item $g(-1,t) \leq -1$, 
so that $h(f)$ starts at $e$, 
\item $g(s,t) \leq -1$ if $t \leq -1$, so that at every time $s$, $h(f)$ is a pseudoisotopy embedding (standard on the bottom face of $I^p \scu I^q$ and all slices of $(I^{p-1} \scu I^{q-1}) \x \R$ below it),
\item $g(1,t) \geq 1$ if $t \geq 1$, so that $h(f)$ ends at a long knot,
\item $g(s,t)$ is independent of $t$ for all $t \geq 1$, so that at every time $s$, $h(f)$ is a pseudoisotopy embedding (given by the same long knot on the top face of $I^p \scu I^q$ and all slices of $(I^{p-1} \scu I^{q-1}) \x \R$ above it), and
\item $g(s,1)=s$, so that $\pi \circ j\circ h = G$.
\end{itemize}
\end{proof}

\begin{lemma}
\label{L:Connectivity}
If $1 \leq p\leq q \leq n-3$, then the space $\PE[p]{q}{n}$ is $(2n - 2q - 5)$-connected.
\end{lemma}
\begin{proof}
Consider the square 
\begin{equation}
\label{Eq:PESquare}
\xymatrix{
\PE[p]{q}{n} \ar[r] \ar[d] & \P_p^n \ar[d] \\
\P_q^n  \ar[r] & \ast
}
\end{equation}
where the maps with domain $\PE[p]{q}{n}$ are obtained by restricting an embedding to each of the two components.  This square is $(2n-p-q-4)$-cartesian.  This fact is the improved version of Morlet's disjunction lemma 
due to Goodwillie \cite[p.~6]{Goodwillie:Thesis}, who further generalized this result \cite[p.~12, Theorem D]{Goodwillie:Thesis} \cite[Lemma 7.4]{Goodwillie-Klein:2015}.  
One can see that the statement about the square \eqref{Eq:PESquare} is equivalent to the improved version of Morlet's disjunction lemma 
by considering the connectivity of the map between the fibers of the rows or columns. 
This property of the square means the map from $\PE[p]{q}{n}$ to the homotopy limit of the rest of the square is $(2n-p-q-4)$-connected.  
This homotopy limit consists of a point in $\P_p^n$, a point in $\P_q^n$, and a path $\Delta^1 \to \ast$ joining their images, hence it is homeomorphic to $\P_p^n \x \P_q^n$.
Another result of Goodwillie implies that the space $\P_p^n$ is $(2n-2p-5)$-connected \cite[pp.~9-10]{Goodwillie:Thesis}  \cite[Lemma 3.16]{Goodwillie:CalcI}.
Thus $\PE[p]{q}{n}$ has connectivity at least $\min\{2n-p-q-5, \ 2n-2p-5, \ 2n-2q-5\} = 2n-2q-5$.  
\end{proof}

\begin{lemma}
\label{L:3DPECube}
If $1 \leq p \leq n-3$, then the space $\P_{p,p,p}^n$ is $(2n-2p-6)$-connected.  
\end{lemma}
\begin{proof}
Consider the 3-cube of spaces 
\[
\xymatrix@R1pc{
& \P_{p,p}^n \ar[r] \ar[dr] & \P_p^n\ar[dr]&  \\
\P_{p,p,p}^n \ar[ur]\ar[r]\ar[dr] &  \P_{p,p}^n\ar[ur]\ar[dr]& \P_p^n \ar[r] & \ast \\
& \P_{p,p}^n\ar[r] \ar[ur] & \P_p^n \ar[ur] & 
}
\]
This cube is $(3n-3p-6)$-cartesian \cite[Lemma 7.4]{Goodwillie-Klein:2015}, i.e., the map from $\P_{p,p,p}^n$ to the homotopy limit $X$ of the rest of the cube is $(3n-3p-6)$-connected.  The space $X$ consists of three points in $\PE{p}{n}$, three paths $\Delta^1 \to \P_p^n$ joining pairs of images of these points, and a map $\Delta^2 \to \ast$ whose three faces are the images of the three paths.  The last piece of data is superfluous, so $X$ fibers over $(\P_{p,p}^n)^3$ with fiber $(\Omega \P_p^n)^3$.  Since both $\P_p^n$ and $\P_{p,p}^n$ are $(2n-2p-5)$-connected (by \cite[Lemma 3.16]{Goodwillie:CalcI} and Lemma \ref{L:Connectivity}), $X$ is $(2n-2p-6)$-connected, hence so is $\P_{p,p,p}^n$.
\end{proof}

\subsection{Bijectivity of graphing and homotopy groups of spaces of long links in a range}
\label{S:BijThm}
We now calculate homotopy groups of 2-component long links, at least up to knotting, in a certain range.  

\begin{theorem}
\label{T:pi_iLinksIso}
If $1 \leq p \leq q \leq n-3$ and $0 \leq i \leq 2n - p - 2q - 4$, then in the sequence of maps 
\begin{equation}
\label{Eq:GraphingComposite}
\pi_{i+p} \LL[0]{q-p}{n-p} 
\overset{G_*}{\longrightarrow}
\pi_{i+p-1} \LL[1]{q-p+1}{n-p+1} \overset{G_*}{\longrightarrow}
 \dots \overset{G_*}{\longrightarrow}  
 \pi_i \LL[p]{q}{n} \overset{G_*}{\longrightarrow} \dots 
 \overset{G_*}{\longrightarrow} \pi_0\LL[i+p]{i+q}{i+n}
\end{equation}
induced by graphing, each map except possibly the first is an isomorphism.
The first map is always a surjection, and it is an isomorphism if $i \leq 2n-p-2q-5$ or $p=q$.
\end{theorem}

The inequality involving $i$, $p$, $q$, and $n$ is precisely the condition that $\pi_{i+p}S^{n-q-1}$ lies within the stable range of the homotopy groups of $S^{n-q-1}$.

\begin{proof}
The maps of spaces underlying the sequence \eqref{Eq:GraphingComposite} are
\begin{equation}
\label{Eq:GraphingSequenceSpaces}
\Omega^{i+p}\LL[0]{q-p}{n-p} \to \Omega^{i+p-1}\LL[1]{q-p+1}{n-p+1} \to 
\dots \to
\Omega^i \LL[p]{q}{n}
\dots 
\to \Omega \LL[i+p-1]{i+q-1}{i+n-1} \to \LL[i+p]{i+q}{i+n}.
\end{equation}
We are interested in the maps they induce on $\pi_0$.
By Lemma \ref{L:GraphingPEFibn}, each such graphing map $G$ fits into a fibration 
\begin{equation}
\label{Eq:GraphingPEFibn}
\Omega^{j+1} \LL[i+p-j-1]{i+q-j-1}{i+n-j-1} 
\xrightarrow{G}
 \Omega^j \LL[i+p-j]{i+q-j}{i+n-j} \to \Omega^j \PE[i+p-j]{i+q-j}{i+n-j}
\end{equation}
where $0 \leq j \leq i+p-1$.  
By Lemma \ref{L:Connectivity}, the connectivity $c$ of the base space in this case satisfies
\[
c\geq (2(i+n-j)-2(i+q-j)-5) - j  = 2n-2q-j-5 \geq 2n - i - p - 2q - 4 
\]
If $i\leq 2n-p-2q-5$, then $c\geq 1$.  Thus each map in the sequence \eqref{Eq:GraphingSequenceSpaces} is $1$-connected and hence induces isomorphisms on $\pi_0$.  
If $i=2n-p-2q-4$, then $j=i+p-1$ yields $c\geq 0$, while
 $j \leq i+p-2$ yields $c \geq 1$.
 Hence in this case, the first map in the sequence \eqref{Eq:GraphingSequenceSpaces} induces a surjection in $\pi_0$ and the remaining maps induce isomorphisms on $\pi_0$.
Finally, if $p=q$, then the first instance of $G_*$ is the graphing map from
$\pi_{i+p}\LL{0}{n-p} \cong \pi_{i+p}S^{n-p-1}$ to $\pi_{i+p-1}\LL{1}{n-p+1}$, which by Theorem \ref{T:Lambda} is also injective.
\end{proof}

We will see in Remark \ref{R:SharpnessBijGraphing} that if $p=q$, the range of values for $i$ in Theorem \ref{T:pi_iLinksIso} is sharp, even if one ignores knotting.  That is, the map $\pi_{2n-3p-3}\Omega^pS^{n-p-1} \to \pi_{2n-3p-3}\LL{p}{n}$ does not surject even onto the factor that does not come from knotting.
However, the sequence of graphing maps involving $\pi_{2n-3p-3}\LL{p}{n}$ consists of isomorphisms if one starts with 2-dimensional instead of 0-dimensional links, as we will see in Theorem \ref{T:GraphingDim>1}.

\begin{corollary}
\label{C:DirectSum}
If $1 \leq p \leq q \leq n-3$ and $0 \leq i \leq 2n - p - 2q - 4$, then
\begin{align}
\label{Eq:pi_iLinksDirectSum}
\pi_i \LL[p]{q}{n}&\cong  \pi_{i+p}S^{n-q-1} \oplus \pi_0 \K_{i+q}^{i+n}
\end{align}
where the inclusion of the first summand is given by
the composite 
\begin{equation}
\label{Eq:InclusionSphereSummand}
\pi_{i+p} S^{n-q-1} \
{{\overset{j_*}{\cong}}} \
\pi_{i+p} (\R^{n-p} - \R^{q-p})  
\xhookrightarrow{\epsilon_*}
\pi_{i+p}\LL[0]{q-p}{n-p} 
{\overset{G^p_*}{\longrightarrow}}
 \pi_i \LL[p]{q}{n}
\end{equation}
with $j$ and $\epsilon$ as in Theorem \ref{T:Lambda}, and where the projection onto the second summand is given by the composite
\begin{equation}
\label{Eq:ProjectionKnotsSummand}
 \pi_i \LL[p]{q}{n} \cong 
 \pi_0 \Omega^i\LL[p]{q}{n}
{\overset{G^i_*}{\longrightarrow}}
 \pi_0\LL[i+p]{i+q}{i+n}
 \xrightarrow{(\rho_2)_*}
 \pi_0\K_{i+q}^{i+n}.
\end{equation}
For the case where $p=q$, we have $\pi_i \LL[p]{p}{n} \cong  \pi_{i+p}S^{n-p-1}$ for all $i \leq 2n - 3p - 4$.
\end{corollary}


\begin{proof}
By Corollary \ref{C:2CompRhoIota}, the restriction fibration $\rho_2$ and the inclusion $\epsilon$ of its fiber yield a direct sum decomposition of $\pi_0$ of each space in the sequence \eqref{Eq:GraphingSequenceSpaces}, starting with 
 \begin{align}
 \label{Eq:pi_iLinksDirectSum1}
 \begin{split}
  \pi_{i+p} \LL[0]{q-p}{n-p} &\cong \pi_{i+p} \Emb_c(\{\ast\}, \, \R^{n-p} - \R^{q-p})  \oplus \pi_{i+p} \K_{q-p}^{n-p}\\
    &\cong \pi_{i+p}  S^{n-q-1}  \oplus \pi_{i+p} \K_{q-p}^{n-p}
 \end{split}
 \end{align}
 and ending with 
 \begin{equation}
\label{Eq:pi_iLinksDirectSum2}
\pi_0 \LL[i+p]{i+q}{i+n} \cong
\pi_0 \Emb_c(\R^{i+p}, \, \R^{i+n} - \R^{i+q}) \oplus \pi_0 \K_{i+q}^{i+n}.
\end{equation}
Each graphing map $G_*$ preserves the direct sum decomposition and hence can be written as $G_*=G'_*\oplus G''_*$.  
We now apply Theorem \ref{T:pi_iLinksIso}.  
If $i \leq 2n-p-2q-5$, each instance of $G'_*$ and $G''_*$ is an isomorphism, so we can decompose $\pi_i \LL[p]{q}{n}$ as the first summand of \eqref{Eq:pi_iLinksDirectSum1} plus the second summand of \eqref{Eq:pi_iLinksDirectSum2}.
If $i=2n-p-2q-4$, the last $p$ instances of $G''_*$ are still isomorphisms because $p\geq 1$, and all but the first instance of $G'_*$ are isomorphisms.  The first instance of $G'_*$ is surjective, and by Theorem \ref{T:Lambda}, it is also injective, so we obtain the direct-sum decomposition in this case too.  
The claim about the maps in \eqref{Eq:InclusionSphereSummand} and \eqref{Eq:ProjectionKnotsSummand} is immediate from this direct-sum decomposition.   

The simplification of the group for $p=q$ holds because $\pi_i\K_p^n=0$ if $i\leq 2n-3p-4$ \cite[Proposition 3.9 (2)]{Budney:Family}.  
(This fact, together with the assumption on $i$, also explains the absence of a summand of $\pi_i\K_p^n$ in the decomposition \eqref{Eq:pi_iLinksDirectSum}).
\end{proof}

One can easily extend Corollary \ref{C:DirectSum} to the case where $0 = p \leq q \leq n-2$, while still assuming $0 \leq i \leq 2n-p-2q - 4$.  In that case, the splitting \eqref{Eq:LinksToKnotsFibn2Comp} gives $\pi_i \LL[0]{q}{n} \cong \pi_{i} S^{n-q-1} \oplus \pi_i \K_{q}^{n}$.

Here is a variation on Theorem \ref{T:pi_iLinksIso} and Corollary \ref{C:DirectSum}.  It allows a higher upper bound on $i$ at the cost of not addressing the graphing map on the factor that comes from knotting.  The proof similarly relies on a disjunction result.  The upper bound on $i$ here and the one in Corollary \ref{C:DirectSum} coincide for $p=q$.

\begin{proposition}[Improved range of graphing isomorphisms on the linking summands]
\label{P:KW-GK}
If $i \leq 2n - 2p - q - 4$, then $\pi_i \LL[p]{q}{n} \cong \pi_{i+p} S^{n-q-1} \oplus \pi_i \K_q^n$, where an inclusion of the first summand is given by a composite involving the graphing map $G^p_*$, as in \eqref{Eq:InclusionSphereSummand}.
\end{proposition}

\begin{proof}
A result of Klein and Williams \cite[Theorem 11.1]{Klein-Williams:2007}, proven in further generality by Goodwillie and Klein \cite[Theorem E]{Goodwillie-Klein:2015}, states that the square
\[
\xymatrix{
\Emb(P, N-Q) \ar[r]\ar[d] & \Emb(P,N)\ar[d]\\
\Map(P, N-Q) \ar[r] & \Map(P,N)
}
\]
is $(2n - 2p - q - 3)$-cartesian.  We apply it to long links (or maps) of $P \scu Q = \R^p \scu \R^q$ in $N=\R^n$, where $\Map(P,N) \simeq \ast$.  Proceeding as in the proof of Lemma \ref{L:Connectivity}, we get that for $i\leq 2n-2p-q-4$
\[
\pi_i \Emb_c(\R^p, \, \R^n - \R^q) \cong 
\pi_i S^{n-q-1} \oplus \pi_i \K_p^n \cong \pi_i S^{n-q-1}
\]
where $\pi_i \K_p^n=0$
because $i \leq 2n-3p-4$.
An application of the splitting \eqref{Eq:pi_iLpqNDecomp1} yields the desired direct sum.  
The left-hand map in the square above is the map $r$ from formula \eqref{Eq:EmbToMap}, of which $G^p$ is a section.  This verifies the claim about the inclusion of the first summand.
\end{proof}

\begin{remark}[Isotopy classes of links]
Lemma \ref{L:SphericalLinks} below implies that for $i=0$, Theorem \ref{T:pi_iLinksIso} reduces to the result of Haefliger \cite{Haefliger:1966CMH} and M.~Skopenkov \cite{Skopenkov:2009}.  Their result applies in a larger range, namely $3n - 2p -2q -6\geq 0$ (roughly, the quadruple-point-free range), in which the group is a direct sum of four factors.  There is one factor for knotting of each of the components, a third factor $\pi_p S^{n-q-1}$, and a fourth factor $\pi_{p+q+2-n} V_M(\R^{M+n-p-1})$ where $V_M(\R^{M+k})$ is the Stiefel manifold of $M$-frames in $\R^{M+k}$ and $M$ is large enough for this group to be stable (or more precisely, $M\geq 2p+q+5-2n$).  In the range $2n-p-2q-4\geq 0$ of Theorem \ref{T:pi_iLinksIso}, there is no knotting of the $\R^p$ component, and the homotopy group of the Stiefel manifold is zero.  
\end{remark}


Putting $\ell=2$ in the next result shows that Theorem \ref{T:pi_iLinksIso} describes the homotopy groups of any space of long links with $m$ components in a certain range.  If $p_1=\cdots=p_m=p$, this range coincides with the range in Theorem \ref{T:pi_iLinksIso} when $q=p$.  

\begin{theorem}
\label{T:ImagesOfiS}
Suppose that $0 \leq \ell \leq m$, $1 \leq p_1 \leq \dots \leq p_m \leq n-3$, and 
$0 \leq i < 1-p_1+\sum_{k=m-\ell+1}^m (n-p_k-2)$.
Then every class in $\pi_i \mathcal{L}_{p_1, \dots, p_m}^n$ is in 
$\displaystyle \sum_{S\subset \{1,\dots,m\}, |S|\leq \ell} \im (\iota_S)_*$, where $\iota_S$ is as in Definition \ref{D:Inclusion}.
\end{theorem}

\begin{proof}
We proceed by induction on $m$, with $m=\ell$ as the basis case.  Suppose we know the theorem to be true for links with $m-1$ components.  
Consider the restriction fibration $\rho_S$, $S=\{2,\dots,m\}$:
\[
\Emb_c\left(\R^{p_1}, \ \R^n -  \coprod_{k=2}^m \R^{p_k} \right)
\to \mathcal{L}_{p_1, \dots, p_m}^n 
\xrightarrow{\rho_S} \mathcal{L}_{p_2, \dots, p_m}^n.
\]
By the direct-sum decomposition from Corollary \ref{C:IntKerPlusSumIm} (with $\mathcal{S}=\{S\}$) and the induction hypothesis, it suffices to prove the analogue of the theorem where $\mathcal{L}^n_{p_1,\dots, p_m}$ is replaced by $\Emb_c\left(\R^{p_1}, \ \R^n -  \coprod_{k=2}^m \R^{p_k} \right)$.  

Our strategy for proving the latter statement is to first separate knotting and braiding phenomena, and then translate braiding phenomena into Whitehead products on wedges of spheres.   
From Lemma \ref{L:HtpyRetract}, we have a retraction 
\begin{equation*}
r:\Emb_c\left(\R^{p_1}, \ \R^n -  \coprod_{k=2}^m \R^{p_k} \right) \to
\Omega^{p_1}\left(\R^{n-p_1} -  \coprod_{k=2}^m \R^{p_k-p_1} \right).
\end{equation*}
Because $r$ and its section induce a splitting of homotopy groups, it suffices to show that classes in the codomain and homotopy fiber of $r$ are in $\im \left(\bigoplus_{|S|\leq \ell} (\iota_S)_*\right) = \sum_{|S|\leq \ell} \im (\iota_S)_*$.

Let $F$ be the homotopy fiber of $r$, which consists of pairs $(f,\gamma)$ such that $f$ is in the domain of $r$ and $\gamma$ is a path from $r(f)$ to the standard embedding $e_{p_1}$.  Similarly define $F'$ as the homotopy fiber of the map 
\[
r':\Emb_c\left(\R^{p_1}, \ \R^n \right) \to
\Omega^{p_1}\left(\R^{n-p_1} \right).
\]
There is a homotopy equivalence $F\to F'$ induced by the inclusion $\R^n - \coprod_{k=2}^m \R^{p_k} \to \R^n$ and the affine-linear map that sends $(t_1^*,0,\dots,0)$ to $0^n$.  Its inverse is induced by an affine-linear inclusion $\R^n \to \R^n -  \coprod_{k=2}^m \R^{p_k}$ whose image lies in $(-\infty, t_2^*) \x \R^{n-1}$.  (Recall that $t_1^*, \dots, t_m^*$ are the first coordinates of the components $e_1,\dots, e_m$ of the standard long link.)  Since the codomain of $r'$ is contractible, we have homotopy equivalences $\Emb_c\left(\R^{p_1}, \ \R^n \right) \to F' \to F$.  In particular, any element in $\pi_*F$ comes from knots, i.e., is in the image of $\iota_1$.

The codomain of $r$ is homotopy equivalent to $\Omega^{p_1}\left(\bigvee_{k=2}^m S^{n-p_k-1}\right)$.  
By the Hilton--Milnor theorem \cite{Hilton:1955} (which we may apply since $n-p_k\geq 3$) and the assumed upper bound on $i$, any element of $\pi_i$ of the latter space is in the image of the map 
\[
\bigoplus_{S\subset \{2,\dots, m\}, |S|\leq \ell-1}  \pi_{i+p_1} \left( \bigvee_{k \in S} S^{n-p_k-1}\right) 
\to \pi_{i+p_1} \left(\bigvee_{k=2}^m S^{n-p_k-1}\right);
\]
indeed, if a Whitehead product of classes from $\ell$ different summands lies in $\pi_{i+p_1}$, then 
$i +p_1\geq 1+ \sum_{k=m-\ell+1}^m (n-p_k-2)$.  A wedge-sum of $\ell-1$ spheres corresponds to a link with $\ell$ components, so the proof is complete.
\end{proof}

\subsection{Spherical links and further calculations for long links}
\label{S:SphericalLinks}
Theorem \ref{T:GraphingDim>1}, the last result of this Section, extends the range of Theorem \ref{T:pi_iLinksIso}, provided we start with 2-dimensional links.  The graphing sequence ends with isotopy classes of long links, so we first digress to relate those to isotopy classes of spherical links, which are known in the range we will consider.

\begin{lemma}
\label{L:SphericalLinks}
If $1\leq p \leq n-3$, then the closure map 
\[
\ \widehat{\cdot} \ : \mathcal{L}_{m \cdot p}^n \to 
\Emb \left( \coprod_1^m S^p, \, S^n \right)
\]
induces a bijection on path components.
\end{lemma}

\begin{proof}
We first check that closure induces a surjection on path components.  Given a spherical embedding $f$, let $f_+$ be its restriction to $\coprod_1^m D^p$, where $D^p$ is the upper hemisphere of $S^p$.  Find an isotopy $F_+=F_+(x,t)$ of $f_+$ to the embedding used to construct the closure of a long link (called $g=(g_1,\dots,g_m)$ in Definition \ref{D:Closure}).  By the isotopy extension theorem, $F_+$ extends to an isotopy $F$ of spherical embeddings starting at $f$.  
Then $F(-,1)$ is in the image of the closure map and is isotopic to $f$, so $\ \widehat{\cdot} \ $ is surjective on $\pi_0$.

For injectivity, we will use an intermediate space, which we will show has the same path components as $\mathcal{L}_{m \cdot p}^n$.
Fix a diffeomorphism $\phi:\coprod_1^m D^p \to  \left( \coprod_{i=1}^m \{(t_i^*,0^{n-p-1})\} \x D^p\right)$, where each summand of $D^p$ in the domain is the upper hemisphere of $S^p$.  
We define the space $\mathrm{DbEmb}\left( \coprod_1^m S^p, \, S^n \right)$ 
of {\bf disk-based embeddings} as the space of pairs $(f,g)$ where 
$f \in \Emb \left( \coprod_1^m S^p, \, S^n \right)$ and $g: D^n \to S^n$ is a smooth, orientation-preserving embedding such that
\begin{itemize}
\item
 $g \phi$ is the restriction of $f$ to $\coprod_1^m D^p$, and
\item 
$g(D^n) \cap \im f = f \left( \coprod_1^m D^p_i \right)$.
\end{itemize}
We topologize it as a subspace of $\mathrm{Emb}\left( \coprod_1^m S^p, \, S^n \right) \x \Emb(D^n, S^n)$.  
This definition is motivated by the $d$-based links used by Habegger and Lin \cite[Definition 2.1]{Habegger-Lin} for link maps with $p=1$ and $n=3$.

We start by establishing a bijection on isotopy classes between disk-based embeddings and long links.
The projection $(f,g)\mapsto g$ to the second factor is a fibration
\begin{equation}
\label{Eq:DbToDisks}
\mathcal{L}_{m\cdot p}^n  \to
\mathrm{DbEmb}\left( \coprod_1^m S^p, \, S^n \right) \to \Emb^+(D^n, S^n)
\end{equation}
where $\Emb^+(D^n, S^n)$ is the space of orientation-preserving embeddings $D^n \to S^n$.  
(The fact that this projection and \eqref{Eq:DbToSphericalLinks} below satisfy the homotopy lifting property can be seen by lifting a homotopy to the space of embeddings of the $n$-manifold with boundary obtained from a neighborhood of $g(D^n) \cup f\left(\coprod_1^m S^p\right)$.)
The inclusion of the fiber over $g$ is $f \mapsto (\widehat{f}, g)$, where $\widehat{f}$ is the closure obtained via the images of $D^p$ under $g$.
The base space is homotopy equivalent to $SO(n+1)$ by shrinking and linearizing.  
In particular, it is path-connected.

It now suffices to check that the boundary map $\pi_1SO(n+1) \to \pi_0 \mathcal{L}_{p\cdot m}^n$ in the long exact sequence of the fibration \eqref{Eq:DbToDisks} is trivial.  
Find a neighborhood of $g(D^n)$ homeomorphic to $D^{n-1} \x I$ such that the induced inclusion $D^n \to D^{n-1} \x I$ is the standard one.  
The generator of $\pi_1SO(n+1)$ in this context can be represented by a loop of rotations of the $D^{n-1}$ factor through an angle of $2\pi$.  The effect on a long link $f$ in the fiber is to send it to (the isotopy class of) the long link obtained by rotating the factor of $I^n$ in $I^n = I^{n-1} \x I$ by $2\pi$, that is, to $f$ itself.  Thus the action is trivial on $\pi_0$, so the inclusion of the fiber in \eqref{Eq:DbToDisks} is a bijection on $\pi_0$.

We now connect to isotopy classes of spherical links. 
The projection $(f,g)\mapsto f$ of a disk-based embedding to the first factor is a fibration
\begin{equation}
\label{Eq:DbToSphericalLinks}
\mathrm{DbEmb}\left( \coprod_1^m S^p, \, S^n \right) \to \Emb\left( \coprod_1^m S^p, \, S^n \right).
\end{equation} 
We will show that it is injective on $\pi_0$ by showing that its fibers are path connected.
Since the closure map factors through the map \eqref{Eq:DbToSphericalLinks}, this will establish injectivity on $\pi_0$ of the closure map.

Because $n-p \geq 3$,  $\pi_0$ of the base space in \eqref{Eq:DbToSphericalLinks} is a group,
so its path components are homotopy equivalent to each other.  Therefore it suffices to consider the fiber $F$ over a standard trivial link.  
By considering a small neighborhood of the interval $[t_1^*, t_m^*] \x 0^{n-1} = \bigcup_{i=2}^m [t_{i-1}^*, t_{i}^*] \x 0^{n-1}$ in $D^n$, we see that $F$ is homotopy equivalent to the space of framed embeddings of $m-1$ intervals with interiors lying in $M:=D^n - \left( \coprod_1^m S^p \right)$ and prescribed values and $p$-frames at the endpoints. 
In turn, $F$ fibers over the space $E$ of unframed such embeddings in $M$, with fiber $\widetilde{F}$ given by the space of such framings.  
By a linearization argument again, $\pi_0 \widetilde{F}$ 
is $\left(\pi_1(O(n-1), O(n-p-1)) \right)^m$, which is trivial since $n-p\geq 3$.  

It remains to show that $E$ is path-connected.
Note that $M$ has a handle decomposition with one $n$-disk and $m$ handles of index $n-p-1$.  
Let $A \in E$.  Since $n-p-1\geq 2$, we can perform an isotopy of each sub-arc of $A$ that lies in a handle, fixing its endpoints, so that it ends up in the $n$-disk.  Since $n \geq 4$, we can take each such isotopy to be an isotopy of $A$, i.e., we can avoid self-intersections.  By a similar argument using that $n \geq 4$, we can find a further isotopy to a fixed standard arc in the $n$-disk.  Thus $E$ is path-connected.  Hence so is $F$, so the map \eqref{Eq:DbToSphericalLinks} is injective on $\pi_0$, and the closure map is bijective on $\pi_0$.
\end{proof}

The work of Budney \cite[Proposition 3.9 (3)]{Budney:Family} gives a stronger analogue of Lemma \ref{L:SphericalLinks} for $\K_p^n$, namely that the closure map is $(n-p-2)$-connected.  

\begin{remark}[Closure of classical long links]
For $p=1$ and $n=3$, the closure map is not injective on isotopy classes for $m\geq 2$.  For example, long links differing by conjugation by a pure braid have the same closure.  In the proof above, we would accordingly have $\pi_0 E \neq \{\ast\}$, i.e., multiple isotopy classes of based embeddings of arcs in $S^3 - \left(\coprod_1^m S^1 \right)$, 
including for $m\geq 3$ distinct classes lying in the image of a 2-disk.  
(For such an embedding of arcs $A_1 \cup \dots \cup A_{m-1} =A$ in $D^2 - \{q_1, \dots, q_m\}$ with $m\geq 3$, a corresponding pure braid $\beta$ can be found by ``thickening'' $A_i$ to an unbased loop $L_i$ enclosing both of its endpoints.  Then take $\beta$ to be any preimage in $\Aut(F_m)$ of the element of $\mathrm{Out}(F_m)$ that sends $x_ix_{i+1}$ to $L_i$.)
For $m=2$, conjugation by pure braids is trivial, but the failure of injectivity of closure can be seen from the fact that $\pi_0\widetilde{F}$ would be nontrivial.
One need only look at long links with 2 and 3 crossings to find an example of two non-isotopic long links with isotopic closure that arises in this way.
\end{remark}

\begin{theorem}
\label{T:GraphingDim>1}
Suppose $1 \leq p \leq n-3$.  
\begin{itemize}
\item[(a)]
For 2-component links, consider the sequence of maps induced by graphing:
\begin{equation*}
\pi_{2n-3p-3}\LL{p}{n} \to \pi_{2n-3p-4}\LL{p+1}{n+1} \to \dots \to \pi_0 \LL{2n-2p-3}{3n-3p-3}.
\end{equation*}
If $p\geq 2$, then all the maps are isomorphisms, and these groups are isomorphic to 
\begin{align*}
\begin{array}{ll}
\Z^3 \oplus \pi_{2n-2p-3} S^{n-p-1} & \text{ if $n-p$ is odd}\\
(\Z/2)^3 \oplus \pi_{2n-2p-3} S^{n-p-1} & \text{ if $n-p$ is even},
\end{array}
\end{align*}
with $\im (\iota_j)_* \cong \Z$ (respectively $\Z/2$)for each $j\in \{1,2\}$ if $n-p$ is odd (respectively even).\\
If $p=1$, then the first map is surjective, and the remaining maps are isomorphisms.  
\item[(b)]
For 3-component links, consider the sequence of maps induced by graphing:
\begin{equation*}
\pi_{2n-3p-3}\mathcal{L}_{3\cdot p}^{n} \to \pi_{2n-3p-4}\mathcal{L}_{3\cdot(p+1)}^{n+1} \to \dots \to \pi_0 \mathcal{L}_{3\cdot(2n-2p-3)}^{3n-3p-3}.
\end{equation*}
If $p\geq 3$, then all the maps are isomorphisms, and these groups are isomorphic to 
\begin{align*}
\begin{array}{ll}
\Z^7 \oplus \left(\pi_{2n-2p-3} S^{n-p-1}\right)^3 & \text{ if $n-p$ is odd}\\
\Z \oplus (\Z/2)^6 \oplus \left(\pi_{2n-2p-3} S^{n-p-1}\right)^3 & \text{ if $n-p$ is even}, 
\end{array}
\end{align*}
with $\im (\iota_j)_* \cong \Z$ (respectively $\Z/2$) for each $j\in \{1,2,3\}$  if $n-p$ is odd (respectively even) \\
and $\bigcap_{|S| \leq 2} \ker (\rho_S)_* \cong \Z$.\\
If $p=2$, then the first map is surjective, and the remaining maps are isomorphisms.  
\end{itemize}
\end{theorem}

We will see in Theorem \ref{T:Knots} that the first graphing map is not an isomorphism for $p=1$.  
Conjecture \ref{C:ExtendThmF} suggests how to complete these descriptions for $p=1$ and $p=2$.

\begin{proof}
We use the fibration 
\[
\Omega \LL{p+j-1}{n+j-1} \to \LL{p+j}{n+j} \to \PE{p+j}{n+j}
\]
where $1 \leq j \leq 2n-3p-3$.
By Lemma \ref{L:Connectivity}, the pseudoisotopy embedding space $\PE{p+j}{n+j}$ is $(2n-2p-5)$-connected.  Thus 
\[
\pi_{i} \LL{p+j-1}{n+j-1} \to \pi_{i-1} \LL{p+j}{n+j}
\]
is an isomorphism if $i \leq 2n-2p-5$ and surjective if $i=2n-2p-4$.
The claims about the sequence itself follow because $p\geq 2$ implies $2n-3p-3 \leq 2n-2p-5$, and $p=1$ yields $2n-3p-3=2n-2p-4$.  
By Lemma \ref{L:SphericalLinks}, $\pi_0 \LL{2n-2p-3}{3n-3p-3} \cong \pi_0\Emb\left( \coprod_1^2 S^{2n-2p-3}, \, S^{3n-3p-3} \right)$.
The proof of the statement for 2-component links is completed using the calculation of the latter group from Haefliger's results \cite[Th\'eor\`eme 10.7]{Haefliger:1966CMH} \cite[Corollary 8.14]{Haefliger:1966Annals}.  Each 3-fold direct sum comes from two summands of $\pi_0 \Emb(S^{2n-2p-3}, S^{3n-3p-3}) \cong \pi_0 \K_{2n-2p-3}^{3n-3p-3}$ and one summand of $\pi_{n-p-1}V_M(\R^{M+n-p-1})$ where $M\geq 2$.
This also proves the claim about $\im (\iota_j)_*$.

In the setting of 3 components, we have an analogous fibration 
\[
\Omega \mathcal{L}_{3\cdot(p+j-1)}^{n+j-1} \to \mathcal{L}_{3\cdot(p+j)}^{n+j} \to \mathcal{P}_{3\cdot(p+j)}^{n+j}.
\]
By Lemma \ref{L:3DPECube}, the base space is $(2n-2p-6)$-connected.  
Thus 
\[
\pi_{i} \mathcal{L}_{3\cdot(p+j-1)}^{n+j-1} \to \pi_{i-1} \mathcal{L}_{3\cdot(p+j)}^{n+j}
\]
is an isomorphism if $i \leq 2n-2p-6$ and surjective if $i=2n-2p-5$.  If $p \geq 3$, then 
$2n-3p-3\leq 2n-2p-6$, while if $p=2$, then $2n-3p-3=2n-2p-5$.
By Lemma \ref{L:SphericalLinks}, $\pi_0 \mathcal{L}_{3\cdot(2n-2p-3)}^{3n-3p-3} \cong \pi_0\Emb\left( \coprod_1^3 S^{2n-2p-3}, \, S^{3n-3p-3} \right)$.
Finally, 
since $4(2n-2p-3)<3((3n-3p-3)-2)$, we can apply another result of
Haefliger \cite[Th\'eor\`eme 9.4]{Haefliger:1966CMH}, which says that 
the latter group is given by the direct sum of isotopy classes of knots and links with fewer than 3 components together with one summand of $\pi_{3n-3p-3}S^{3n-3p-3}$.  The remaining 6 summands of $\Z$ or $\Z/2$ come from knots on each of the 3 components and links on each of the 3 pairs of components, and the 3 summands of $\pi_{2n-2p-3} S^{n-p-1}$ come from links on each of the 3 pairs of components.
This yields the claimed groups in both parities as well as the claims about $\im (\iota_j)_*$ and $\bigcap_{|S|<3} \ker (\rho_S)_*$.
\end{proof}

The following rational result of other authors is for comparison with Theorem \ref{T:GraphingDim>1} and for later use in Section \ref{S:Knots} in the case where $p=1$.

\begin{theorem}
\label{T:TripodHedgehog}
\cite{Songhafouo-Turchin:HHA, Fresse-Turchin-Willwacher:Emb}
If $1\leq p \leq n-3$, then 
\[
\pi_{2n-3p-3}\LL{p}{n} \otimes \Q
\cong
\left\{
\begin{array}{ll}
\Q^4 & n-p \text{ odd}, \\
\Q^3 &  n-p \text{ even and } p=1\\
0 & \text{ else}\\ 
\end{array}
\right. 
\]
\end{theorem}

\begin{proof}
This follows from a result of Fresse, Turchin, and Willwacher \cite{Fresse-Turchin-Willwacher:Emb},
a special case of which was obtained by Songhafouo Tsopm\'en\'e and Turchin \cite{Songhafouo-Turchin:HHA} in higher codimensions, including all cases where $1=p\leq n-3$.
They describe these groups via a complex of graphs with leaves labeled by the link components.  
It is roughly dual to the graph complex described below in Section \ref{S:ConfSpaceInt}, except that there are no link strands, only labels by them, and the graphs are required to be connected.
It generalizes work of Arone and Turchin \cite{Arone-Turchin:Htpy} from one component to multiple components.
If all source components have dimension $p$ and the target has dimension $n$, the degree of a graph is $(n-1)|E| - n |I| - p|L|$, where $|E|$, $|I|$, and $|L|$ are the numbers of edges, internal (i.e., non-leaf) vertices, and leaves respectively.  

From this, one can deduce that only two types of graphs can contribute to $\pi_{2n-3p-3} \mathcal{L}_{p,p}^n$.  One type is the ``tripod'' graph (the trivalent tree with 3 leaves) in degree $2n-3p-3$, with the 4 possible leaf-labelings by $\{1,2\}$, for $n-p$ odd.
The other type is the ``2-hair hedgehog'' graph (the trivalent graph with 2 leaves and one double edge) in degree $2n-2p-4$, with the 3 possible leaf-labelings by $\{1,2\}$, for $n-p$ even and $p=1$.  
 (Each of these graphs is 2-torsion in one parity of $n-p$, due to certain orientation relations.)
 These are the first two graphs shown in \cite[Section 3.2, Table B]{Arone-Turchin:Htpy}.
(In addition, by \cite[Theorem 3.1]{Arone-Turchin:Htpy}, 
these are the only graphs that contribute to $\pi_i \mathcal{L}_{p,p}^n \otimes \Q$, $i\leq 2n-3p-3$, apart from a single edge with distinct leaf-labels in degree $n-2p-1$, where we know $\pi_{n-2p-1}\mathcal{L}_{p,p}^n \cong \Z$ by Corollary \ref{C:DirectSum}.)
\end{proof}

\section{Homotopy classes in spaces of long knots and links from joining pure braids}
\label{S:Knots}



Our main purpose now is to prove the last main result, Theorem \ref{T:E}. where we describe generators for  the homotopy groups of spaces of links that we have calculated.  In it, we also relate the image of graphing from spheres and configuration spaces to the previously known first nontrivial homotopy groups of spaces of long knots. 

In Section \ref{S:ConfSpaceInt}, we briefly review configuration space integrals, which produce cohomology classes in spaces of links from a certain cochain complex of graphs (a.k.a.~diagrams).  
We specify a handful of graph cocycles in Section \ref{S:Cocycles}.
In Section \ref{S:DualHomologyClasses}, we describe how to produce dual homology classes by resolving singular links.  
The proof of Theorem \ref{T:Knots} (a.k.a.~Theorem \ref{T:E}), given in Section \ref{S:KnotsThm}, involves identifying certain such homology classes with the homotopy classes in the Theorem statement, or more precisely their images under the Hurewicz map.  It also relies on known results for long knots of various dimensions and the compatibility of graphing and joining components.
We discuss possible future directions in Section \ref{S:Questions}.

\subsection{Configuration space integrals for cohomology of spaces of long 1-dimensional links}
\label{S:ConfSpaceInt}

The proof of the last result will use several cohomology classes in spaces of long links that we will describe via configuration space integrals.
We now digress for a brief overview of these integrals for long $m$-component, 1-dimensional links. 
Full details are given in our previous joint work \cite{KMV:2013}, which extended work of Cattaneo, Cotta-Ramusino, and Longoni \cite{Cattaneo:2002} and Bott and Taubes \cite{Bott-Taubes:1994} from knots to long links.  

The classes produced by this construction are indexed by linear combinations of {\bf link diagrams}, which are graphs on $m$ oriented line segments, for example as in formulas \eqref{Eq:KappaCocycle} through \eqref{Eq:LkNum3} below.
Such a diagram $\Gamma$ consists of a vertex set $V(\Gamma)$ and an edge set $E(\Gamma)$.  The vertices are partitioned into sets $V_{\mathrm{seg}}(\Gamma)$ and $V_{\mathrm{free}}(\Gamma)$ of {segment vertices} and {free vertices}.  An edge joining two segment vertices is called a chord.  Diagrams where all the edges are chords, which play a somewhat special role, are called a chord diagrams.  A part of one of the $m$ segments between two vertices is called an arc.  All vertices are required to have valence at least 3, where both edges and arcs count towards valence.  For brevity, we will say ``graph'' or ``diagram'' to mean ``link diagram.'' 

For a smooth manifold $X$, the configuration space $\Conf(m,X)$ has a compactification $\Conf[m,X]$ due to Fulton and Macpherson.  It is a smooth manifold with corners homotopy equivalent to $\Conf(m,X)$.  
For any graph $\Gamma$, one defines a bundle $\xi_\Gamma$ over the space $\mathcal{L}_{m \cdot 1}^n$ of $m$-component 1-dimensional links.
To be able to compactify at infinity, one uses a version of $\mathcal{L}_{m \cdot 1}^n$ where the standard link has all $m$ components approach infinity in different directions.  
The total space of $\xi_\Gamma$ is the pullback of the diagram 
\[
\mathcal{L}_{m \cdot 1}^n \x \Conf_0\left[V_{\mathrm{seg}}(\Gamma), \coprod_1^m \R\right] \to 
\Conf\left[V_{\mathrm{seg}}(\Gamma), \R^n\right] \leftarrow 
\Conf\left[V(\Gamma), \R^n\right].
\]
In the left-hand configuration space, the subscript 0 indicates that the segment vertices are required lie on the components and in the order given by the graph $\Gamma$.
Thus the fiber of $\xi_\Gamma$ over a link $f$ is the subspace $\Conf_\Gamma$ of $\Conf[V(\Gamma), \R^n]$ where the points in $V_{\mathrm{seg}}(\Gamma)$  lie in the image of $f$, on certain components and in the order determined by $\Gamma$.  Each pair of vertices $i,j$ determines a map $\phi_{ij}:\Conf(V(\Gamma), \R^n) \to S^{n-1}$ given by the unit vector from point $i$ to point $j$.  Let $\omega$ be a unit volume form on $S^{n-1}$.  One sends a graph $\Gamma$ to a differential form via the map
\[
I: \Gamma \mapsto \int_{\Conf_\Gamma} \bigwedge_{e \in E(\Gamma)} \phi_{ij}^*\omega
\]
where the integration is along the fiber $\Conf_\Gamma$ of the bundle $\xi_\Gamma$ over the space of links.  

To orient the configuration space over which one integrates and determine the sign of the form to be integrated, one needs certain labelings on the graphs $\Gamma$, which depend on the parity of $n$.  For $n$ odd, one needs an ordering of the vertices and an orientation of each edge.  For $n$ even, one needs an ordering of the segment vertices and an ordering of the edges.  
Changing the labeling changes the integral only by a sign, so if two labeled graphs $\Gamma$ and $\Gamma'$ differ by a permutation $\sigma$ of labels and (for $n$ odd) $r$ edge reversals, one sets $\Gamma = (-1)^{\mathrm{sign}(\sigma)+r} \Gamma'$, 
An equivalence class of labeling is viewed as an orientation of a graph.

There is a coboundary operator $\delta$ on such graphs that makes the space of graphs into a cochain complex.
On a graph $\Gamma$, it is defined as the signed sum of edge contractions
\[
\delta \Gamma := \sum_{e} \eps(e) \Gamma/e
\]
over all arcs $e$ and all edges $e$ that are not chords or self-loops.  We first define the sign $\eps(e)$ for $n$ odd.  If $e$ is an edge or arc with endpoints $i < j$, then
\begin{equation*}
\eps(e):= 
\left\{
\begin{array}{ll}
(-1)^j 
& \text{ if } e=(i \to j) \\
(-1)^{j+1}
 & \text{ if } e=(i \leftarrow j)
\end{array}
\right.
\end{equation*}
If $n$ is even and $e$ is an arc with endpoints $i < j$, define $\eps(e)$ as above (where the arc orientation comes from the orientation of the segments).  If $n$ is even and $e$ is an edge, set 
\begin{equation*}
\eps(e) = 
(-1)^{e+1+|V_{seg}(\Gamma)|}
\end{equation*}
where by abuse of notation $e$ also denotes the label on this edge.

With this coboundary operator, the integration map $I$ is a cochain map if $n \geq 4$.  
At the level of cohomology, $I$ is known to be injective on the subspace of trivalent graphs \cite{Cattaneo:2002}.  
It produces all of the real cohomology of the space of braids in dimension at least 4 \cite{KKV:2020}.

\subsection{Some low-dimensional graph cocycles}
\label{S:Cocycles}

We will focus on a handful of graph cocycles, which for $m=1,2,3$ yield cohomology classes in dimension $2n-6$ in the space of $m$-component $1$-dimensional long links in $\R^n$ for any $n \geq 4$.  

The graphs below are oriented using the following shorthand: segment vertices are ordered first by their order on the components (with e.g.~all those on the first component coming first), a free vertex is last, edges are oriented from smaller to larger label, and edges are ordered by the smallest label of an endpoint.  
With this convention, some of the formulas giving cocycles are the same in both parities of $n$ except for some signs; for these classes ($\kappa$, $\mu$, and the $\nu_i$), we list the sign for $n$ odd above the sign for $n$ even.  
For the other classes ($\eta$ and $\lambda$), we separately list the formulas for the two parities (indicated by subscripts), but we later use these symbols without the subscripts to refer to both parities at once.

\begin{align}
\begin{split}
\label{Eq:KappaCocycle}
\kappa &:=  \quad
\raisebox{-4.5pc}{\includegraphics[scale=0.22]{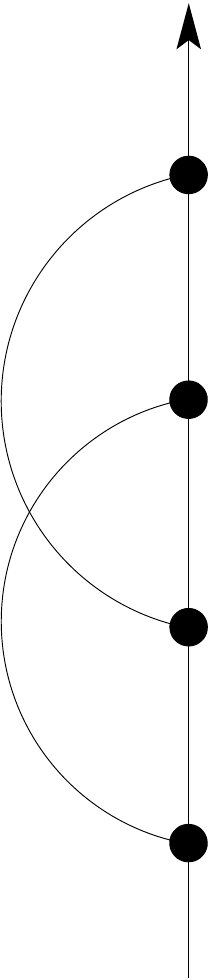}}\quad \mp \quad 
\raisebox{-4.5pc}{\includegraphics[scale=0.22]{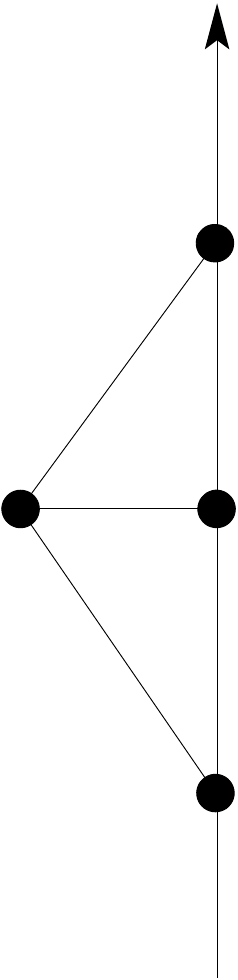}}.\\
\end{split}
\end{align}
\begin{align}
\begin{split}
\label{Eq:EtaCocycleOdd}
\eta_{\mathrm{odd}} &:= \quad
\raisebox{-2.5pc}{\includegraphics[scale=0.3]{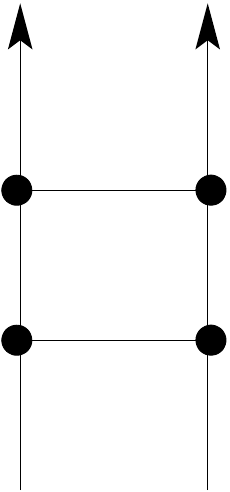}}\quad  - \quad
\raisebox{-2.5pc}{\includegraphics[scale=0.3]{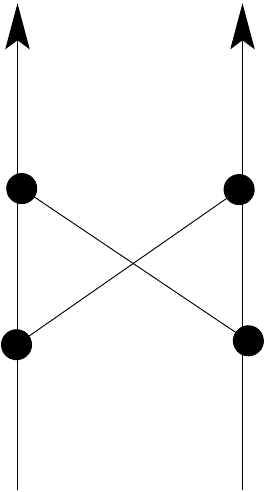}} \\
 \end{split}
 \end{align}
\begin{align}
\label{Eq:LambdaCocycleOdd}
\lambda_{\mathrm{odd}}:= \quad
\raisebox{-2.5pc}{\includegraphics[scale=0.3]{crossing-chords-unlabeled.pdf}} \quad - \quad 
\raisebox{-2.5pc}{\includegraphics[scale=0.3]{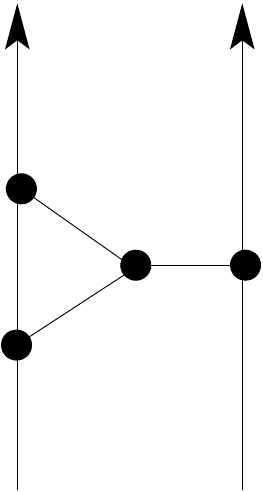}} \quad - \quad 
\raisebox{-2.5pc}{\includegraphics[scale=0.3]{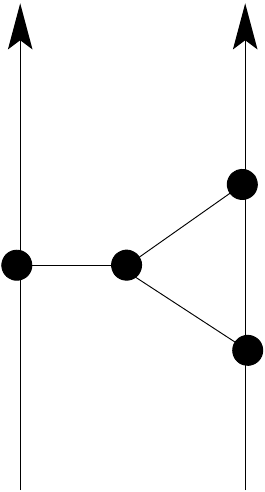}} \quad + \quad 
\raisebox{-2.5pc}{\includegraphics[scale=0.3]{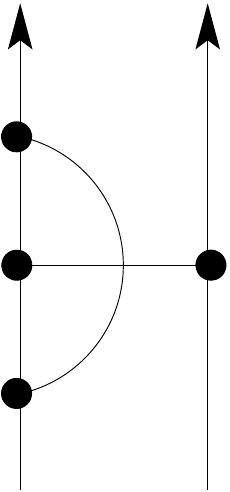}} \quad + \quad 
\raisebox{-2.5pc}{\includegraphics[scale=0.3]{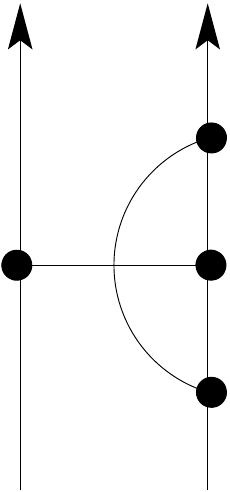}}  
\end{align}
 \begin{align}
\label{Eq:EtaCocycleEven}
\eta_{\mathrm{even}}:= \quad
\raisebox{-2.5pc}{\includegraphics[scale=0.3]{double-chord-unlabeled.pdf}} \quad + \quad 
\raisebox{-2.5pc}{\includegraphics[scale=0.3]{tripod-left-unlabeled.pdf}} \quad + \quad 
\raisebox{-2.5pc}{\includegraphics[scale=0.3]{tripod-right-unlabeled.pdf}} \quad + \quad 
\raisebox{-2.5pc}{\includegraphics[scale=0.3]{psi-left-unlabeled.pdf}} \quad + \quad 
\raisebox{-2.5pc}{\includegraphics[scale=0.3]{psi-right-unlabeled.pdf}}  
\end{align}
\begin{align}
\label{Eq:LambdaCocycleEven}
\lambda_{\mathrm{even}}:= \quad
\raisebox{-2.5pc}{\includegraphics[scale=0.3]{crossing-chords-unlabeled.pdf}} \quad + \quad 
\raisebox{-2.5pc}{\includegraphics[scale=0.3]{tripod-left-unlabeled.pdf}} \quad - \quad 
\raisebox{-2.5pc}{\includegraphics[scale=0.3]{tripod-right-unlabeled.pdf}} \quad + \quad 
\raisebox{-2.5pc}{\includegraphics[scale=0.3]{psi-left-unlabeled.pdf}} \quad - \quad 
\raisebox{-2.5pc}{\includegraphics[scale=0.3]{psi-right-unlabeled.pdf}}  
\end{align}
 \begin{align}
\begin{split}
\label{Eq:MuCocycle}
\mu&:= \ 
\raisebox{-2.5pc}{\includegraphics[scale=0.3]{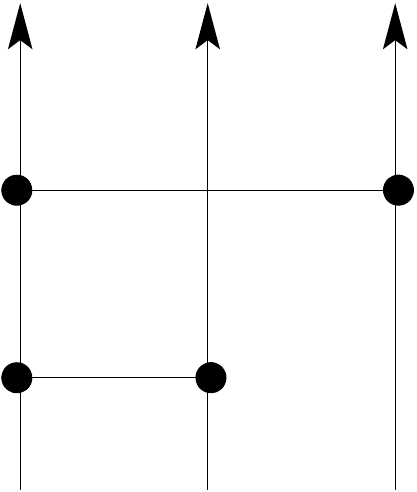}} \  + \
\raisebox{-2.5pc}{\includegraphics[scale=0.3]{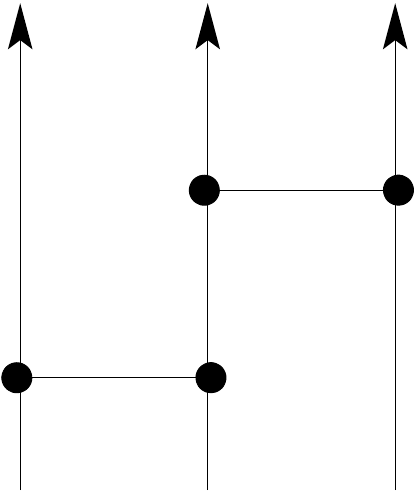}} \ + \ 
\raisebox{-2.5pc}{\includegraphics[scale=0.3]{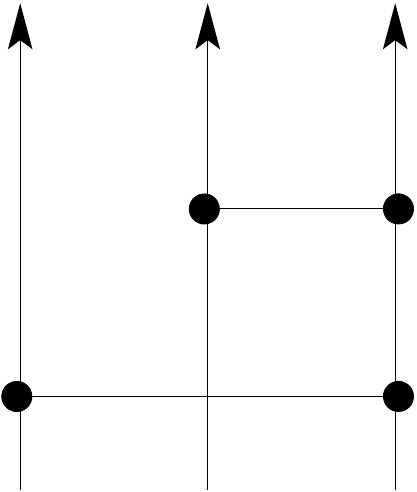}} \ \mp \ 
\raisebox{-2.5pc}{\includegraphics[scale=0.3]{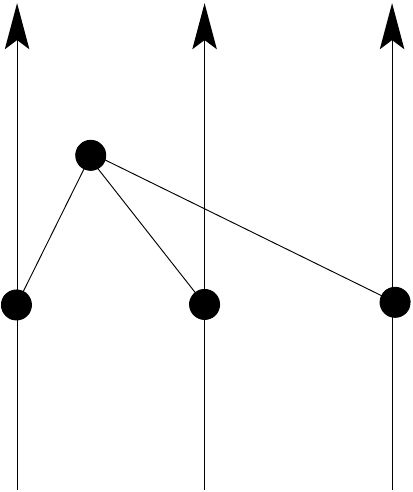}}
\end{split}
\end{align}
\begin{align}
\label{Eq:LkNum12}
\nu_1:= \
\raisebox{-2.5pc}{\includegraphics[scale=0.3]{LM-unlabeled.pdf}} \  \mp \
\raisebox{-2.5pc}{\includegraphics[scale=0.3]{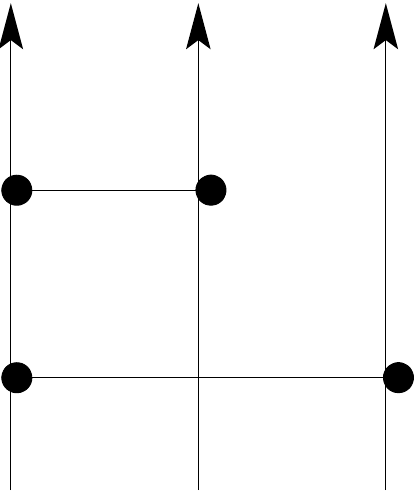}}  \qquad
\nu_2:= \
\raisebox{-2.5pc}{\includegraphics[scale=0.3]{LR-unlabeled.pdf}} \  \mp \
\raisebox{-2.5pc}{\includegraphics[scale=0.3]{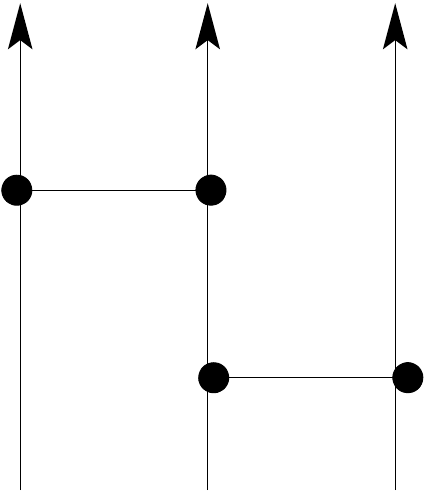}} 
\end{align}
\begin{align}
\label{Eq:LkNum3}
\nu_3:= \
\raisebox{-2.5pc}{\includegraphics[scale=0.3]{MR-unlabeled.pdf}} \  \mp \
\raisebox{-2.5pc}{\includegraphics[scale=0.3]{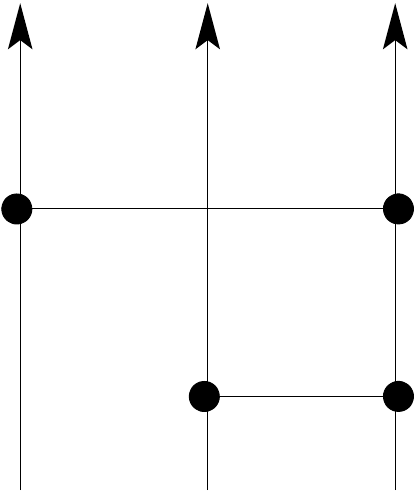}} 
\end{align}

We will use these classes above to identify families of 1- , 2-, and 3-component long links.
Recall the restriction maps $\rho_i: \mathcal{L}_{1,1}^n \to \K_1^n$, $i=1,2$ and $\rho_{ij}: \mathcal{L}_{1,1,1}^n \to \mathcal{L}_{1,1}^n$, $1 \leq i < j \leq 3$ from Definition \ref{D:Restriction}.

\begin{lemma}
\label{L:Cocycles} \  
\begin{itemize}
\item[(a)]
The space $H^{2n-6}(\K_1^n; \, \R)$ is spanned by  $I(\kappa)$.
\item[(b)] 
The class $I(\eta)$ spans a subspace of $H^{2n-6}(\mathcal{L}_{1,1}^n; \, \R)$ 
which pulls back via the graphing map $\Omega S^{n-2} \simeq \Omega\Conf(2,\R^{n-1}) \to \mathcal{L}_{1,1}^n$  
 isomorphically onto $H^{2n-6}(\Omega\Conf(2,\R^{n-1}); \R)$.  
If $n$ is even, the evaluation of $I(\eta)$ on the image of this graphing map (pre-composed with the Hurewicz map) is up to a sign the Hopf invariant of a class in $\pi_{2n-5}S^{n-2}$.  
\item[(c)] 
The set $\{I(\eta), I(\lambda), \rho_1^*I(\kappa), \rho_2^*I(\kappa)\}$ is a basis for $H^{2n-6}(\mathcal{L}_{1,1}^n; \, \R)$.  
\item[(d)]
The set $\{ I(\mu), I(\nu_1), I(\nu_2), I(\nu_3), \rho_{12}^*I(\eta), \rho_{13}^*I(\eta), \rho_{23}^*I(\eta)\}$ is a basis for a subspace of \\ $H^{2n-6}(\mathcal{L}_{1,1,1}^n; \, \R)$ that maps isomorphically onto $H^{2n-6}(\Omega\Conf(3,\R^{n-1}); \, \R)$ via the graphing map $\Omega\Conf(3,\R^{n-1}) \to \mathcal{L}_{1,1,1}^n$.  
\end{itemize}
\end{lemma}

\begin{proof}
The definition of the coboundary operator ensures that with our conventions for graph orientations, these cochains are cocycles.  They lie in a dimension which admits no nonzero coboundary map, so they represent nontrivial and moreover linearly independent cohomology classes.

Part (a) then follows from the fact that  $H^{2n-6}(\K_1^n; \, \Z) \cong \Z$, established in work of Turchin \cite{Turchin:JKTR2007}.

For parts (b) and (d), we use the description of $H_*(\Omega\Conf(m, \R^{n-1}); \R)$, due to Cohen and Gitler \cite[Theorem 2.3]{Cohen-Gitler:TAMS}, as the space of horizontal chord diagrams modulo the 4T relations, where each chord has degree $n-3$.  (The cohomology is thus isomorphic to Vassiliev invariants of pure braids.)  
Write a horizontal chord between strand $i$ and $j$ as $t_{ij}$, where $1\leq i < j \leq m$, and write a horizontal chord diagram as a product of $t_{ij}$.  A basis of horizontal chord diagrams in dimension $d(n-3)$ (before the quotient by 4T) is given by all degree-$d$ monomials in the $t_{ij}$, and we visualize the dual basis of functionals by the same types of diagrams.  The 4T relations are the homogeneous degree-2 elements $[t_{ij}, t_{jk} + t_{ik}]$ with $i,j,k$ distinct.  Thus $H^{2n-6}(\Omega\Conf(m, \R^{n-1}); \R)$ is isomorphic to the subspace of functionals in $\mathrm{span}\{(t_{ij}t_{k\ell})^*\}_{i<j, k<\ell}$ which vanish on all elements $[t_{ij}, t_{jk} + t_{ik}]$.  

By a theorem on graphing braids in previous joint work \cite[Theorem 5.19]{KKV:2020}, the graphing map $G:\Omega \Conf(m, \R^{n-1}) \to \mathcal{L}_{m \cdot 1}^n$ induces a surjection $G^*$ in cohomology.
That theorem also implies that for the link diagrams $\sum c_i \Gamma_i$ in parts (b) and (d), 
$G^* \circ I (\sum c_i \Gamma_i) \in H^{*}(\Omega\Conf(m, \R^{n-1}); \R)$ corresponds to the result of remembering only the (duals of) horizontal chord diagrams.  
Hence for $m=2$, $\eta \mapsto (t_{12} t_{12})^*$, and $(t_{12} t_{12})^*$ forms a basis for $H^{2n-6}(\Omega\Conf(2, \R^{n-1}); \R)$.  
The fact that $I(\eta)$ is the Hopf invariant is established in our joint work \cite[Example B.2]{KKV:Primitive}.  
Thus part (b) is proven.

For $m=3$, there are the diagrams $\mu, \nu_1, \nu_2, \nu_3, \rho_{12}^*\eta, \rho_{13}^*\eta,$ and $\rho_{23}^*\eta$, where $\rho_{ij}^*\eta$ is the result of putting $\eta$ on strands labeled $i$ and $j$ with no edges incident to the remaining third component (so $\rho_{ij}^*I(\eta)=I(\rho_{ij}^*\eta)$).
We have 
$\mu \mapsto (t_{12}t_{13})^* + (t_{12}t_{23})^* + (t_{13}t_{23})^*$,
$\nu_i \mapsto (t_{ij}t_{ik})^* - (t_{ik}t_{ij})^*$ (where $i,j,k$ are distinct), and
$\rho_{ij}^*\eta \mapsto (t_{ij}t_{ij})^*$.
By the theorem on graphing braids \cite[Theorem 5.19]{KKV:2020}, each of these must vanish on the 4T relations.
One can directly check their linearly independence in $\mathrm{span}\{(t_{ij}t_{kl})^*\}_{i<j, k<\ell}$.  
We now check the dimension of the space of horizontal 2-chord diagrams modulo 4T for $m=3$.
There are nine diagrams $t_{ij}t_{k\ell}$, with $1\leq i<j\leq 3$, $1\leq k<\ell \leq 3$.  
The three 4T relations (each determined by the pair $i<j$) span a dimension-2 subspace, so the quotient by them has rank 7.  Thus for $m=3$, our 7 cocycles form a basis for $H^{2n-6}(\Omega\Conf(m, \R^{n-1}); \R)$, and part (d) is proven.
(The class $I(\mu)$ is an analogue of the triple linking number, while $\rho_{ij}^*I(\eta)$ is an analogue of the square the linking number between strands $i$ and $j$, and $I(\nu_i)$ is an analogue of the product of the two different pairwise linking numbers involving strand $i$.)

Finally, for part (c), by Theorem \ref{T:TripodHedgehog} for $p=1$, we have that $\pi_{2n-6} \LL{1}{n}\otimes \Q$ has dimension 3 or 4 according to whether $n$ is odd or even.
Since $\L_{1,1}^n$ is a homotopy-associative H-space, the Milnor--Moore theorem \cite{Milnor-Moore} says that its rational homology is the universal enveloping algebra on its rational homotopy.  
Let $a$ denote a generator of $\pi_{n-3}\L_{1,1}^n \otimes \Q$ as well as its image in rational homology under the Hurewicz map.  Under this map, the Whitehead product in homotopy corresponds, up to a sign, to the graded commutator $[ -,\, -]$ under the Pontryagin product $\cdot$ in homology.
If $n$ is even, then $[a,a]=\pm 2a\cdot a$, while if $n$ is odd, $[a,a]=0$ so $a \cdot a$ is not primitive.
Thus in either parity, the rational homology (or cohomology) in degree $2n-6$ is 4-dimensional.  Hence the linearly independent set in question is a basis.
\end{proof}

\subsection{Homology classes dual to configuration space integrals}
\label{S:DualHomologyClasses}

We now describe a process for constructing a family $F$ of 1-dimensional links out of an immersion with $k$ double-points.  Lemma \ref{L:Pairing} below will show that such a family $F$ is dual to a certain chord diagram $\Gamma$ with $k$ chords.

The family $F$ is constructed as follows.  Let $f$ be an immersion $f: \coprod_1^m \R \to \R^n$ with a finite number $k$ of double-points and no intersections of higher multiplicity.  We require that the tangent vectors are linearly independent at each double-point of $f$.  Suppose the double-points are given by $f(s_1)=f(t_1), \dots, f(s_k)=f(t_k)$.  The $s_i$ may lie in different components, as may the $t_i$, and we make no assumptions on the order in which they lie even if some of them do lie on the same component.   
The unit sphere in the orthogonal complement to the two tangent vectors at each double-point is $S^{n-3}$.  Build a $k(n-3)$-parameter family $F$ of (nonsingular) $m$-component links by varying  the strand on which $t_i$ lies through a normal sphere $S^{n-3}_i$ of a small radius within a small neighborhood $N_i$ in $\R^n$ of the $i$-th double-point for $i=1,\dots, k$.  More precisely, parametrize the family by $(D^{n-3})^{\x k}$ and define it using degree-1 maps $(D^{n-3}, \d D^{n-3}) \to (S^{n-3}_i, \ast)$ around $t_i$ for each $i$.  We can arrange for its support to lie in neighborhoods of $t_1, \dots, t_k$, interpolating between the degree-1 maps and constant maps within these neighborhoods.

\begin{lemma}
\label{L:Pairing}
Let $F$ be the family constructed out of a singular link $f$ as just described, and let $\Gamma$ be a graph.  Then $\langle I(\Gamma), F\rangle$ is nonzero if and only if $\Gamma$ is the chord diagram with chords corresponding exactly to the singularities in $f$.
\end{lemma}

For example, if $f=(f_1,f_2)$ has two double points, one where $f_1(s_1)=f_2(t_1)$ and another where $f_1(s_2)=f_2(t_2)$ , with $s_1<s_2$ and $t_1<t_2$, then $F$ pairs nontrivially with the first graph in formula \eqref{Eq:EtaCocycleOdd}, but trivially with the second graph in that formula.  This $F$ would also pair trivially with all the graphs in formula \eqref{Eq:LambdaCocycleOdd}, including the three chord diagrams in that formula.

\begin{proof}[Sketch of proof of Lemma \ref{L:Pairing}]
The proof is based on that of a result of Cattaneo, Cotta-Ramusino, and Longoni \cite[Theorem 6.1]{Cattaneo:2002}.
First, on the space of configurations where fewer than 2 points lie in some neighborhood $N_i$, any integral vanishes because on that space, the variation of $F$ through the $i$-th copy of $S^{n-3}$ admits a homotopy back to the original singular link, thus yielding a degenerate family.  
This implies that if $\Gamma$ is not a chord diagram, then $\langle I(\Gamma), F \rangle=0$.

It remains to consider the integrals of chord diagrams on the space $\mathcal{C}$ of configurations where each $N_i$ contains exactly two points.  
The configuration points in a chord diagram can be identified by their order on the link components, so there is only one way for each $N_i$ to contain a pair of points, i.e., $\mathcal{C}$ is connected.  
If points $j$ and $k$ do not lie in the same $N_i$, then $\phi_{jk}$ restricted to $\mathcal{C}$ is nullhomotopic.  
Thus for a chord diagram $\Gamma$,  $\langle I(\Gamma), F\rangle \neq 0$ only if the chords of $\Gamma$ correspond exactly to the singularities resolved in $F$.

Finally, consider a factor $\phi_{jk}^*\omega$ in $I(\Gamma)$. Its integral over the configuration space $I \x I$ of two points on segments of the link in some $N_i$ produces an $(n-3)$-form.  The pairing of this form with the $(n-3)$-parameter family coming from the resolution of the $i$-th double point is the same as the linking number of the sphere $S^{n-2}$ obtained from the $(n-3)$-parameter family that varies one segment and (a closure of) the other segment.
Then by construction, if the chords in $\Gamma$ correspond to the resolved singularities of $F$, $\langle I(\Gamma), F\rangle$ is product of factors $\pm 1$, hence is itself $\pm 1$.

This argument can be made more precise by taking the limit as the diameters of the neighborhoods $N_i$ and the spheres $S^{n-3}$ inside them approaches zero.  
\end{proof}

\subsection{Long Borromean rings, the Hopf map, and classes of long links and knots}
\label{S:KnotsThm}

In this section, we will use the map $J$ from Definition \ref{D:Join} that joins the last two components of a long link in $\mathcal{L}_{m\cdot p}^n$.  
For $p=1$, we will use the following alternative description of $J$.
Given $f=(f_1, \dots, f_m) \in \mathcal{L}_{m\cdot 1}^n$, remove $f_{m-1}((2,\infty))$ and $f_m((-\infty, -2))$.
Then join $f_{m-1}(2)$ to $f_m(-2)$ by an arc outside $[-1,1]^n$.  Perform an isotopy of this arc and what remains of the image of $f_m$, and pre-compose and post-compose the resulting embedding with diffeomorphisms of $\R$ and $\R^n$ respectively so that the required conditions for a long link are satisfied and so that the standard $m$-component link is taken to the standard $(m-1)$-component link.  
Although we will consider $n\geq 4$, the arc and isotopy may be taken to lie in a 3-dimensional subspace of $\R^n$.
One can show that this description agrees up to homotopy with Definition \ref{D:Join} for a suitable choice of path $\gamma$: the arc corresponds to $\gamma \x 0^{n-2} \x -1$ in the framed tubular neighborhood of $\gamma$, and the restriction of the isotopy to the point $t_m^* \x 0^{n-2} \x 1$ corresponds to $\gamma \x 0^{n-2} \x 1$.  
Since we will assume $n-p \geq 3$, Remark \ref{R:ClassesOfJoin} implies that it suffices for $\gamma$ to end on the positive side of the $(m-1)$-th component.
We are now ready to prove our last main result, Theorem \ref{T:E}:

\begin{theorem}
\label{T:Knots}
Suppose $1 \leq p\leq n-3$ and $2n-3p-3\geq 0$.
\begin{itemize}
\item[(a)]
If $n-p$ is odd, then the map $\pi_{2n-2p-3} S^{n-p-1}\to \pi_{2n-3p-3}\K_p^n (\cong \Z)$ given by $p$-fold graphing followed by joining the two link components sends the Whitehead square $[\mathbbm{1}_{n-p-1}, \mathbbm{1}_{n-p-1}]$ of the identity to twice a generator.
Thus if $n-p=3$, $5$, or $9$, it sends the Hopf fibration to a generator. 
\item[(b)]
The map $\pi_{2n-3p-3}\Omega^p\Conf(3,\R^{n-p}) \to \pi_{2n-3p-3}\K_p^n (\cong \Z$ or $\Z/2$) induced by the composition of $p$-fold graphing followed by joining the three components together maps the ``parametrized long Borromean rings'' class $[b_{21}, b_{31}]$ to a generator.
\item[(c)] 
For $p=1$, there is a basis for $\pi_{2n-6}\LL{1}{n}$ modulo torsion, 
consisting of the images of a generator of $\pi_{2n-6}\K_1^n$ under the inclusions $\iota_1,\iota_2: \K_1^n \to \LL{1}{n}$;
the result of graphing and then joining two components of $[b_{21},b_{31}]$; and
for $n-p$ odd, the result of graphing $[\mathbbm{1}_{n-2}, \mathbbm{1}_{n-2}]$.

If $p\geq 2$, then 
$\pi_{2n-3p-3}\LL{p}{n}$ is generated by the two inclusions of a generator of $\pi_{2n-3p-3}\K_p^n$; the result of graphing and then joining two components of $[b_{21},b_{31}]$; and the image of a generating set of $\pi_{2n-2p-3}S^{n-p-1}$ under graphing.

If $p\geq 3$ and $m\geq 3$, then  $\pi_{2n-3p-3}\mathcal{L}_{m\cdot p}^{n}$ is generated by the $m$ inclusions of a generator of $\pi_{2n-3p-3}\K_p^n$; the result of graphing and then joining two components of $[b_{21},b_{31}]$ for every pair of components $(i,j)$ with $1\leq i < j \leq m$; the image under graphing of a generating set of $\pi_{2n-2p-3}S^{n-p-1}$ for every $(i,j)$ with $1\leq i < j \leq m$; and the result of graphing $[b_{21},b_{31}]$ for every $(i,j,k)$ with $1 \leq i < j < k \leq m$.
\end{itemize}
\end{theorem}

\begin{proof}
The overall strategy is to first prove all the statements for $p=1$ and then use diagram \eqref{Eq:BigDiagram} to prove the statements for larger values of $p$.  
The work of Budney will imply that the bottom row, which involves only spaces of knots, consists almost entirely of isomorphisms; Theorems \ref{T:Lambda} and \ref{T:GraphingDim>1} give crucial information about the rows of graphing maps involving multiple components.

To prove the statements for $p=1$, we will define families of 1-dimensional $m$-component long links, for $m=1,2,3$, by resolving singular links and joining components.  
These families will be called $F$ and $F'$ (for $m=3$), $H$, $L$, and $L'$ (for $m=2$), and $K$ (for $m=1$).  
Via maps $J$ that join components, we will have $(F\pm F') \mapsto (L\pm L') \mapsto K$ and $H\mapsto K$ by construction.
Our first main task is to identify the homology classes of these four cycles using Lemmas \ref{L:Cocycles} and \ref{L:Pairing}, which includes showing that $F\pm F'$ is (up to a sign) the image of $[b_{21}, b_{31}]$, $H$ is (up to a sign) the image of the Hopf map, and $K$ is a generator of $H_{2n-6}\K_1^n$.

\textbf{Families $F$ and $F'$ of 3-component 1-dimensional pure braids}:
Take a singular 3-strand braid $f=(f_1,f_2,f_3): \coprod_1^3 \R \incl \R^n$ with two double-points, given by $f_1(s_1)=f_2(t)$ and $f_1(s_2)=f_3(u)$ with $s_1<s_2$.  
One can construct $f$ by starting with a singular braid in $\R^3$, for example as in Figure \ref{F:Singular3Braid} and taking its image under an inclusion of coordinates $i:\R^3 \incl \R^n$.  
Then construct a $(2n-6)$-parameter family $F$ by resolving the singularities, as described before Lemma \ref{L:Pairing}.
\begin{figure}[h!]
(a) \includegraphics[scale=0.16]{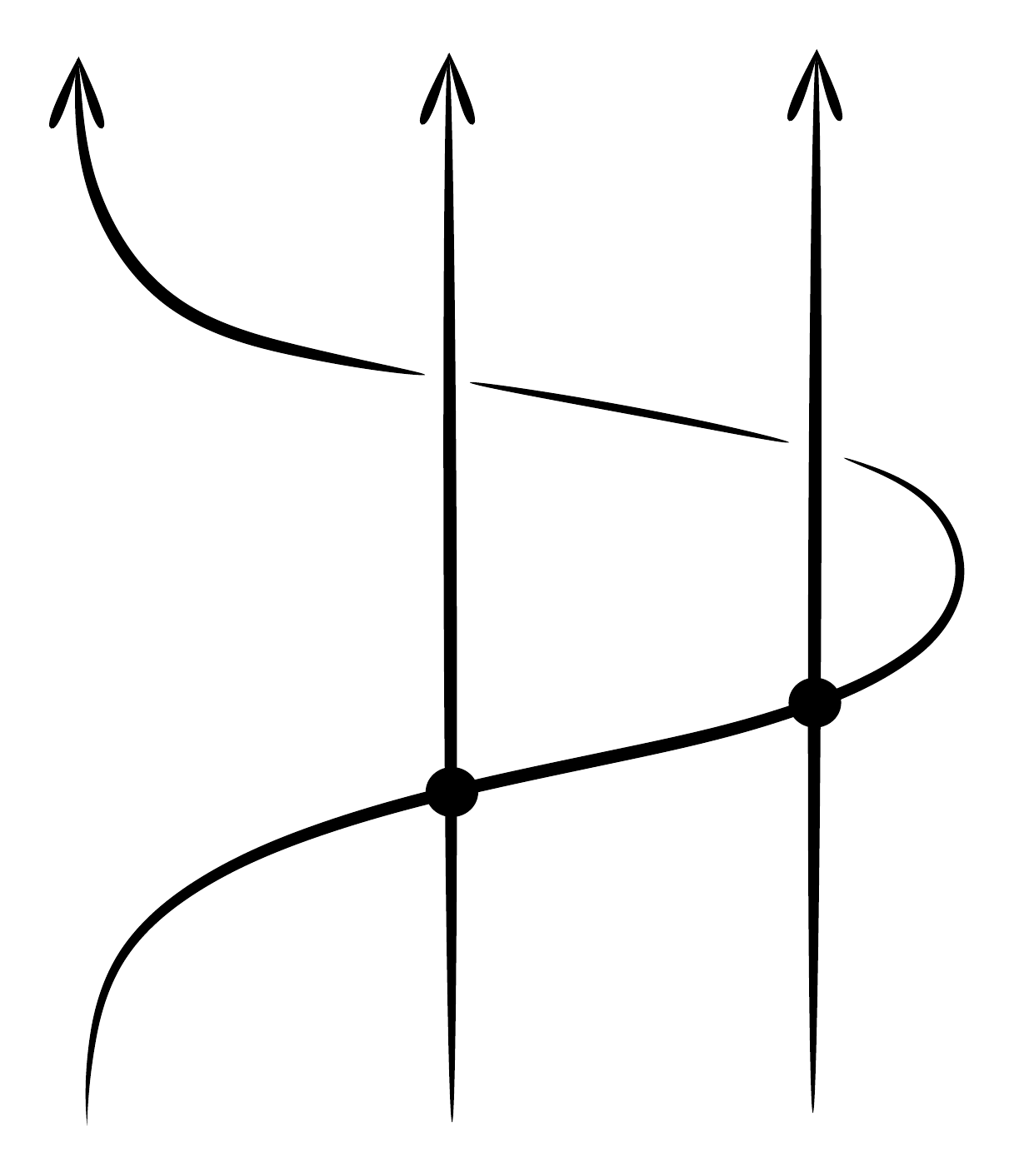} \qquad \qquad \qquad
(b) \includegraphics[scale=0.16]{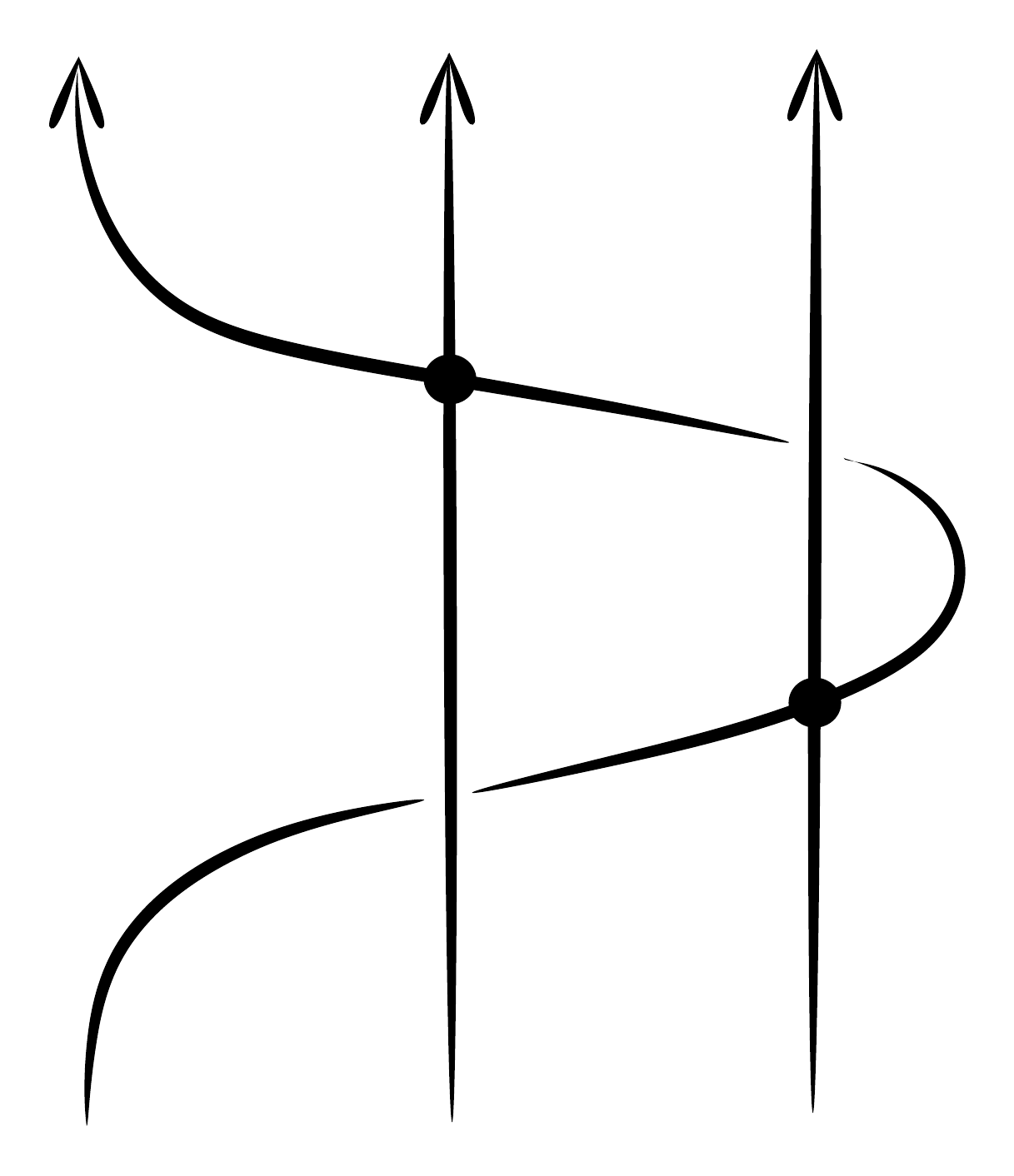}
\caption{ 
(a) The singular 3-strand braid $f=(f_1,f_2,f_3)$ used to build the family $F$.
(b) The singular 3-strand braid $f'=(f_1',f_2',f_3')$ used to build the family $F'$.
}
\label{F:Singular3Braid}
\end{figure}

Next repeat this process, but with a singular 3-strand braid $f'=(f_1', f_2', f_3')$ with two double-points, given by $f_1'(s_1)=f_3'(u)$ and $f_1'(s_2)=f_2'(t)$ with $s_1<s_2$.  
Let $F'$ be the $(2n-6)$-parameter family obtained by resolving the singularities of $f'$.  
We orient it by ordering the parameters as $s_1, s_2, t, u$ and ordering the spheres of resolution in the order that the double-points are listed.

Both $F$ and $F'$ can be taken to be families of pure braids, rather than arbitrary long links.
As a result, we can use Lemma \ref{L:Cocycles} (d).  There, we noted that the graphing map $\Omega \Conf(3,\R^{n-1}) \to \mathcal{L}_{1,1,1}^{n}$ induces an injection in homology \cite{KKV:2020}, and we gave a basis for the dual of its image.
We apply Lemma \ref{L:Pairing} to the cycles $F$ and $F'$  and to the seven graph cocycles $\gamma$ in that basis.  
Only the first graph in the expressions for $\mu$ and $\nu_1$ pairs nontrivially with $F$, and the same is true for the second graph in the expression for $\nu_1$ and $F'$.
Diagrammatically, $F$ is dual to 
$\includegraphics[scale=0.07]{LM-unlabeled.pdf}$ 
and 
$F'$ is dual to 
$\includegraphics[scale=0.07]{ML-unlabeled.pdf}$.
That is, we have that
\begin{align}
\label{Eq:FPairing}
\begin{split}
\langle I(\mu), F \rangle &=1 \\
\langle I(\nu_1), F \rangle &= 1\\ 
\langle I(\nu_i), F \rangle &= 0 \text{ for } i=2,3 \\
\langle \rho_{ij}^*I(\eta), F \rangle &= 0 \text{ for } 1\leq i<j\leq 3
\end{split}
\begin{split}
\langle I(\mu), F' \rangle &= 0 \\
\langle I(\nu_1), F' \rangle &= \mp 1 \\
\langle I(\nu_i), F' \rangle &= 0  \text{ for } i=2,3\\
\langle \rho_{ij}^*I(\eta), F' \rangle &= 0  \text{ for } 1\leq i<j\leq 3
\end{split}
\end{align}
where as before, the sign for $n$ odd is shown above the sign for $n$ even.  
The class $[b_{21}, b_{31}]$ is known to generate $\pi_{2n-6}\Omega\Conf(3,\R^{n-1})\otimes \Q$.  
From our previous joint work \cite[Example 4.8]{KKV:Primitive}, we have $\langle I(\mu), [b_{21}, b_{31}] \rangle = 1$ (possibly up to a sign) and $\langle I(\gamma), [b_{21}, b_{31}] \rangle = 0$ for all the other graph cocycles $\gamma$ appearing in \eqref{Eq:FPairing}.
Hence $F\pm F'$ (for $n$ odd, respectively even) is homologous to $[b_{21}, b_{31}]$, as desired.

\textbf{Families $L$ and $L'$ of 2-component 1-dimensional long links by joining components}:
We next define $(2n-6)$-dimensional families $L$ and $L'$ of 2-component long links 
as the images of $F$ and $F'$ under the map that joins the second and third components, as described just before the Theorem statement.  
We can also describe $L$ and $L'$ as the resolutions of the singular links $\ell$ and $\ell'$ obtained by joining components of the singular links $f$ and $f'$.  (That is, joining components commutes with resolving singularities.)
These singular links are shown in Figures \ref{F:Singular2Braid} and \ref{F:Singular2Link}.
\begin{figure}[h!]
\raisebox{-5pc}{\includegraphics[scale=0.15]{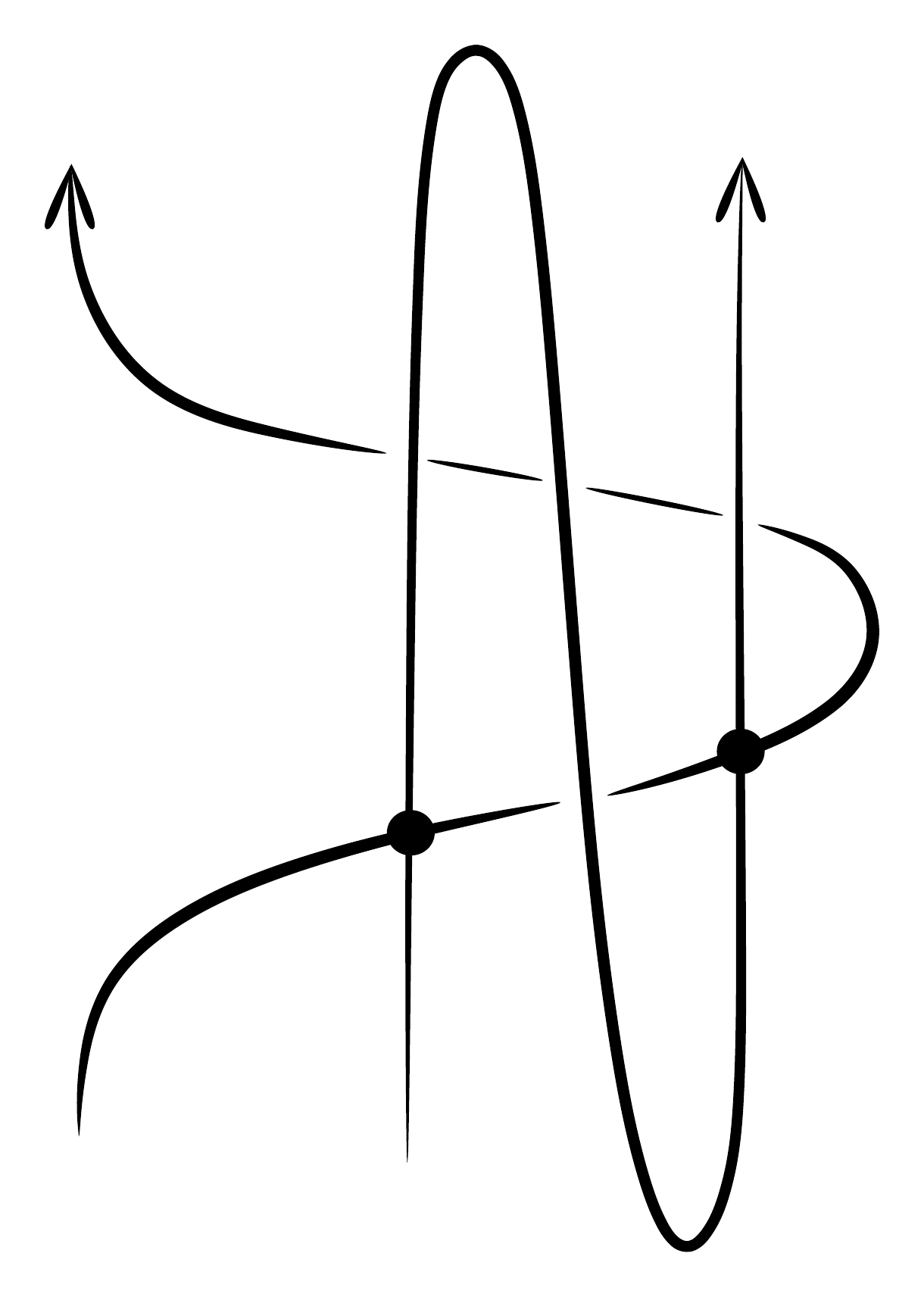}} \quad $\sim$ \quad
 \raisebox{-4.36pc}{\includegraphics[scale=0.15]{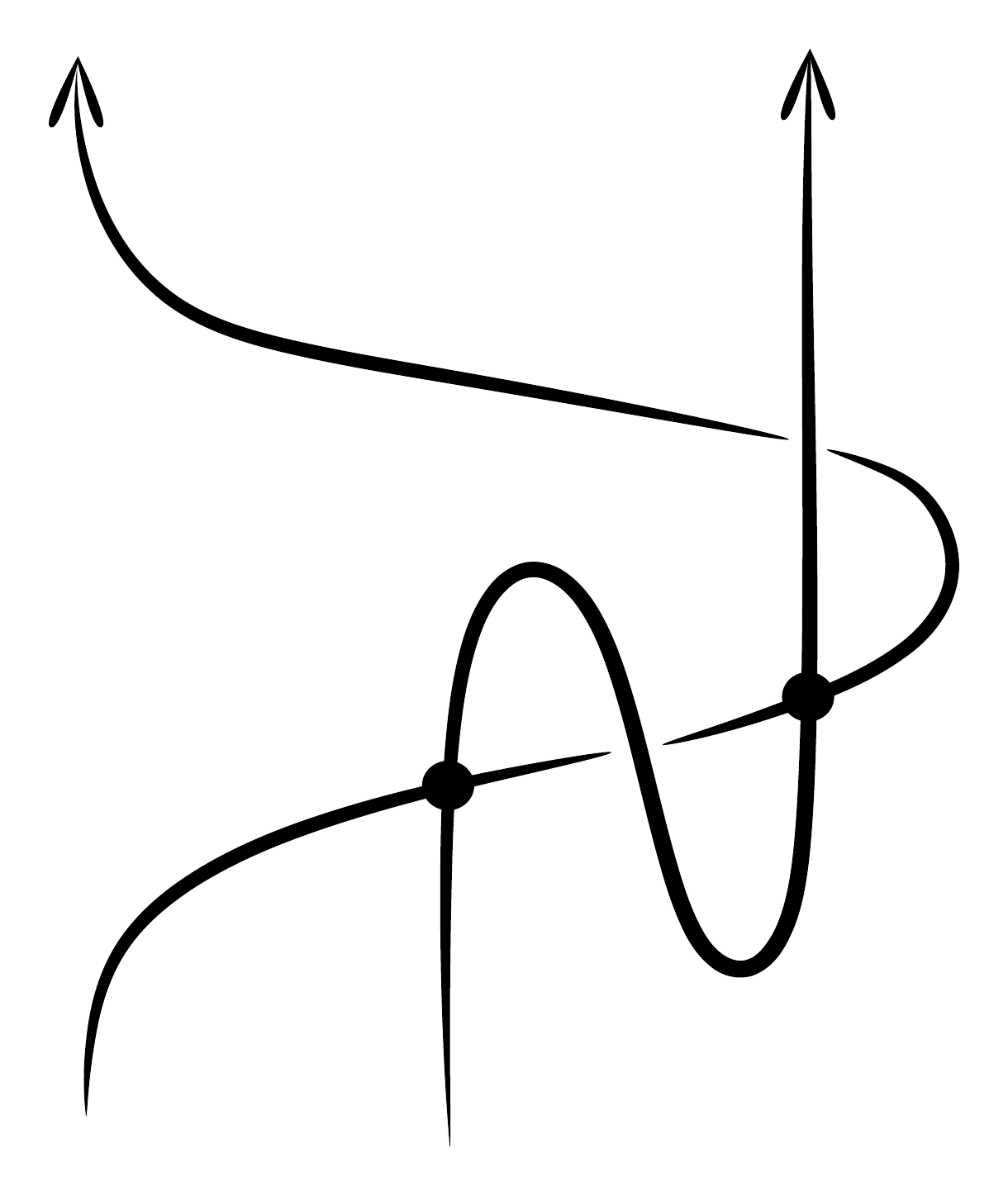}} \quad $\sim$ \quad
 \raisebox{-3.7pc}{\includegraphics[width=5pc, height=8pc]
{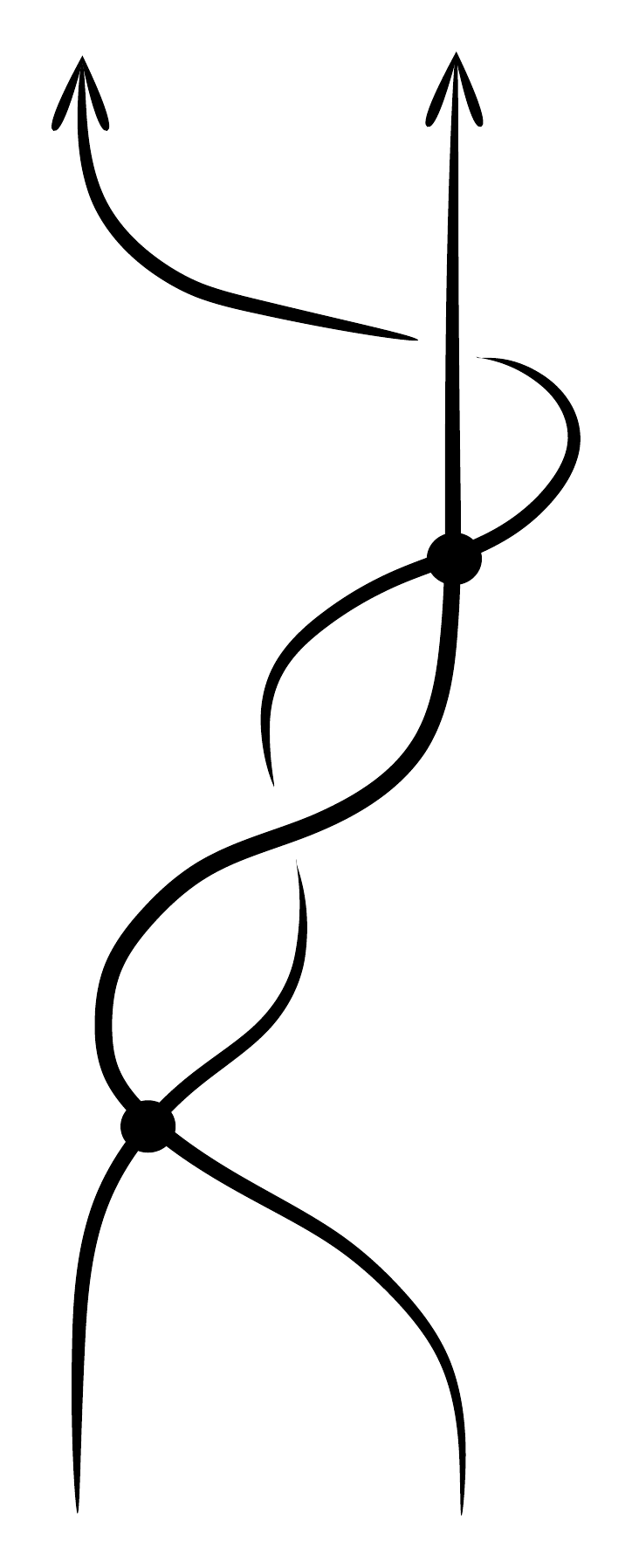}} 
\caption{The singular 2-component long link $\ell=(\ell_1, \ell_2)$obtained by joining strand 2 to strand 3 in the singular braid $f=(f_1,f_2,f_3)$.  An isotopy takes this singular long link to a singular 2-strand pure braid.}
\label{F:Singular2Braid}
\end{figure}

\begin{figure}[h!]
\raisebox{-5pc}{\includegraphics[scale=0.15]{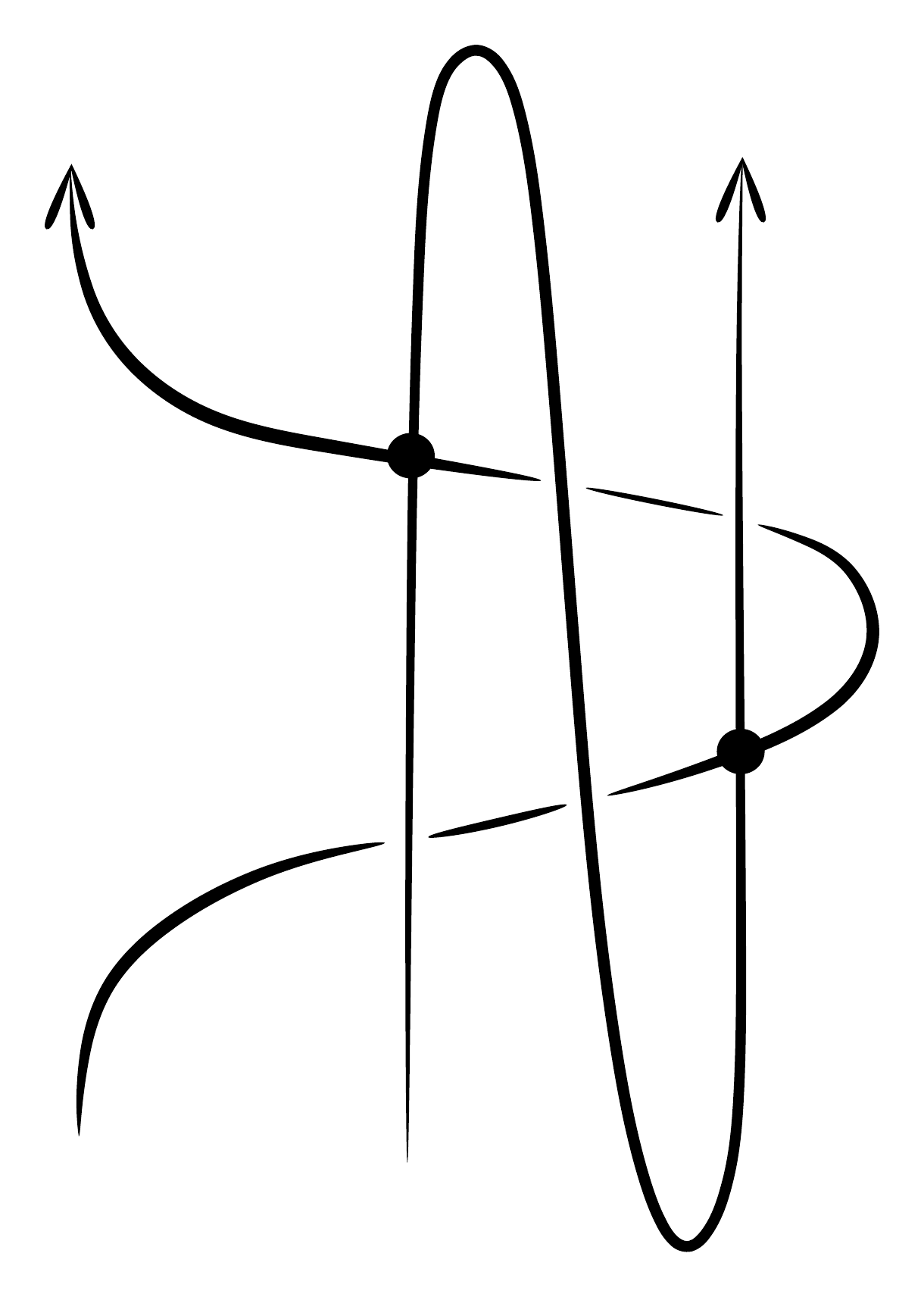}} \quad $\sim$ \quad
\raisebox{-5pc}{\includegraphics[scale=0.15]{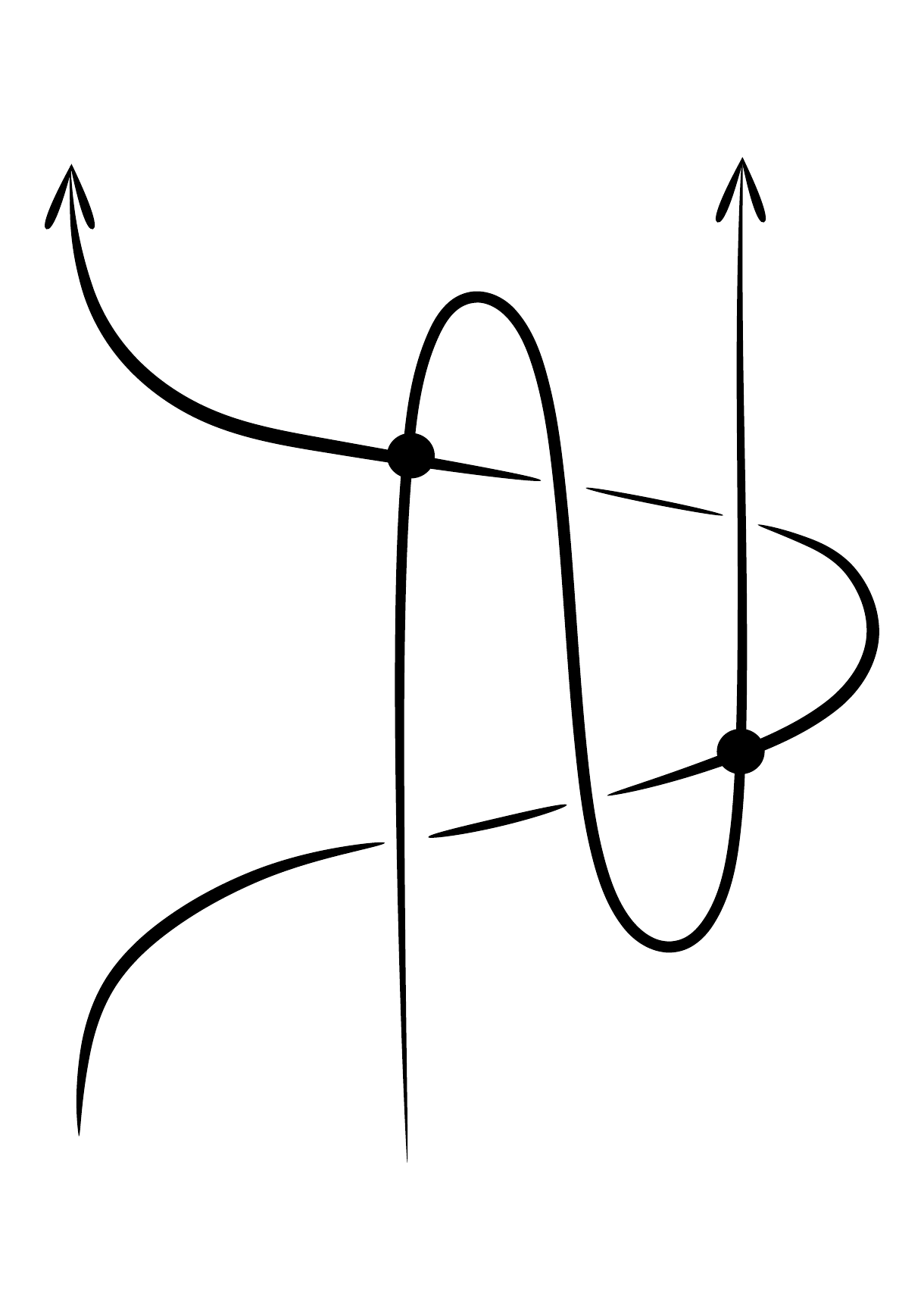}}
\caption{The singular 2-component long link $\ell'=(\ell_1', \ell_2')$ obtained by joining strand 2 to strand 3 in the singular braid $f'=(f'_1,f'_2,f'_3)$.}
\label{F:Singular2Link}
\end{figure}

By Lemma \ref{L:Pairing}
\begin{align*}
\begin{split}
\langle I(\eta), L \rangle &= 1 \\
\langle I(\lambda), L \rangle &= 0 \\
\langle \rho_i^*I(\kappa), L \rangle &= 0 \text{ for } i=1,2 
\end{split}
\begin{split}
\langle I(\eta), L' \rangle &= \left\{\begin{array}{ll} -1 & \text{ if $n$ is odd}\\ 0 &\text{ if $n$ is even}\end{array}\right. \\
\langle I(\lambda), L' \rangle &= 1 \\
\langle \rho_i^*I(\kappa), L' \rangle &= 0 \text{ for } i=1,2. 
\end{split}
\end{align*}
Indeed, $L$ pairs nontrivially with only the first diagram in formulas \eqref{Eq:EtaCocycleOdd} and \eqref{Eq:EtaCocycleEven} for $\eta$, and $L'$ pairs nontrivially with only the first diagram in formulas \eqref{Eq:LambdaCocycleOdd} and \eqref{Eq:LambdaCocycleEven} for $\lambda$.  
Diagrammatically, $L$ is dual to 
$\includegraphics[scale=0.07]{double-chord-unlabeled.pdf}$ 
and 
$L'$ is dual to 
$\includegraphics[scale=0.07]{crossing-chords-unlabeled.pdf}$.

\textbf{A family $H$ of 2-component 1-dimensional pure braids}:
Next, we construct a $(2n-6)$-parameter family $H$ of 2-component 1-dimensional long links.  Start with an immersion $h=(h_1, h_2)$ which is an embedding except at two double-points, where the double-points are given by $h_1(s_1)=h_2(t_1)$ and $h_1(s_2)=h_2(t_2)$ where $s_1 < s_2$ and $t_1 < t_2$.  
See Figure \ref{F:SingularHopfBraid}.  
Resolve the singularities to obtain a $(2n-6)$-parameter family $H$ out of $h$.
We orient it by ordering the parameters as $s_1, s_2, t_1, t_2$ and ordering the resolution spheres in the order that the double-points are listed.
By Lemma \ref{L:Pairing}
\begin{align}
\label{Eq:HPairings}
\begin{split}
\langle I(\eta), H \rangle &=1 \\
\langle I(\lambda), H \rangle &= 0 \\
\langle \rho_i^*I(\kappa), H \rangle &= 0 \text{ for } i=1,2.
\end{split}
\end{align}
Diagrammatically, $H$ is dual to 
$\includegraphics[scale=0.07]{double-chord-unlabeled.pdf}$ 
and 
$L'$ is dual to 
$\includegraphics[scale=0.07]{crossing-chords-unlabeled.pdf}$.
Thus $L\pm L'$ and $H$ represent linearly independent elements of $H_{2n-6}(\LL{1}{n}; \, \Z)$, though $L$ and $H$ represent the same element.

\begin{figure}
\includegraphics[scale=0.16]{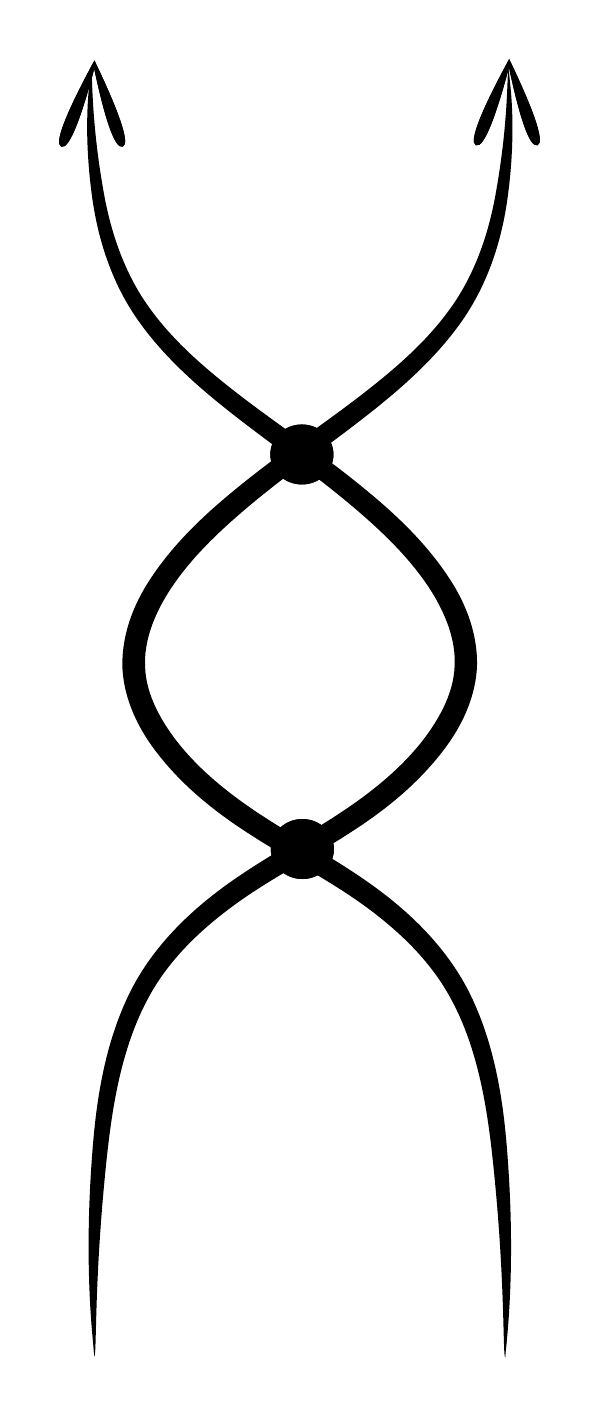} \qquad
\caption{
The singular 2-component braid $h=(h_1,h_2)$ used to construct the family $H$.
}
\label{F:SingularHopfBraid}
\end{figure}

We now identify $H$ with the image of a better known element in the homotopy groups of spheres.
Recall that if $n$ is even, then the free part of $\pi_{2n-5}(S^{n-2})$ maps isomorphically onto $H_{2n-6}(\Omega S^{n-2})$, which via graphing injects into $H_{2n-6}(\LL{1}{n})$.
Lemma \ref{L:Cocycles} (b) and the pairings in \eqref{Eq:HPairings} show that if $n-2=2,4,$ or $8$, then $H$ represents (up to a sign) the class of the pure braid coming from the Hopf fibration.  
More generally, this shows that $2H$ and $\pm [\mathbbm{1}_{n-2}, \mathbbm{1}_{n-2}]$ represent the same class, since the Whitehead square of an even-dimensional identity map has Hopf invariant $\pm 2$.

\textbf{A family $K$ of 1-dimensional long knots}:
Finally, we define a $(2n-6)$-dimensional family $K$ of long knots in $\R^n$ as the image of the family $H$ of 2-component pure braids under the map $J:\LL{1}{n}\to \K_1^n$ described just before the Theorem statement.
As noted earlier, joining commutes with resolving singularities.  Thus we can alternatively define $K$ by applying an obvious extension of the map $J$ to the singular link $h$ to get a singular knot $k$, as shown in Figure \ref{F:KnotFromHopf}(b), and then resolving the singularities of $k$.
One could also obtain $K$ without reference to multiple-component links, defining $k$ directly as a singular knot prescribed by (or dual to) the first diagram 
$\includegraphics[scale=0.07]{CD2-long-knots-unlabeled.pdf}$ 
 in formula \eqref{Eq:KappaCocycle}.
Thus $K$ generates $\pi_{2n-6}\K_1^n$ because the latter definition matches the construction of a generator of $\pi_{2n-6}\K_1^n$ in Budney's work \cite{Budney:Family}, where a quadrisecant counting argument is used to prove that it is a generator.
(By Lemma \ref{L:Pairing}, $\langle I(\kappa), K \rangle = 1$, but this implies only that $K$ is nonzero in $H_{2n-6}(\K_1^n; \, \Z) \cong \pi_{2n-6}\K_1^n$, so we need to use Budney's result.)


\begin{figure}[h!]
\raisebox{-7.6pc}{\includegraphics[scale=0.15]{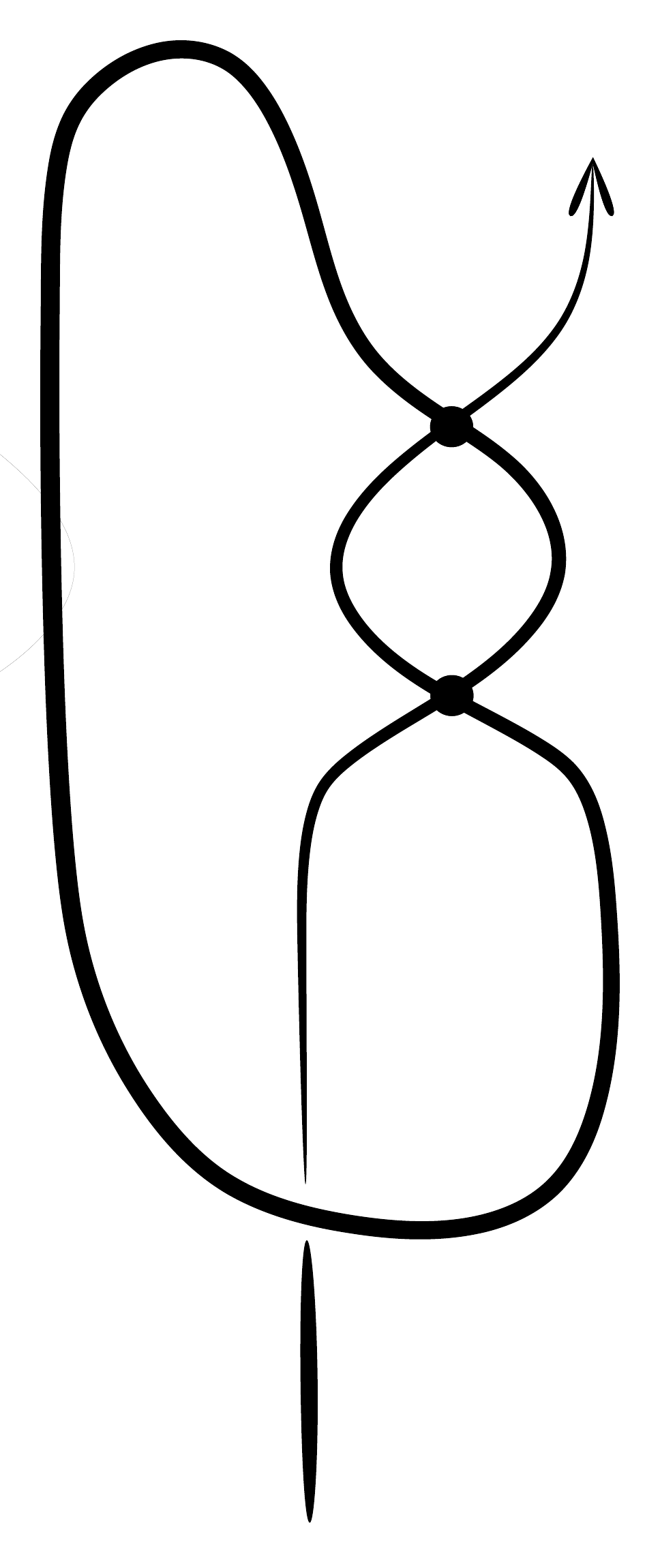}} \quad $\sim$ \quad
 \raisebox{-6.5pc}{\includegraphics[scale=0.15]{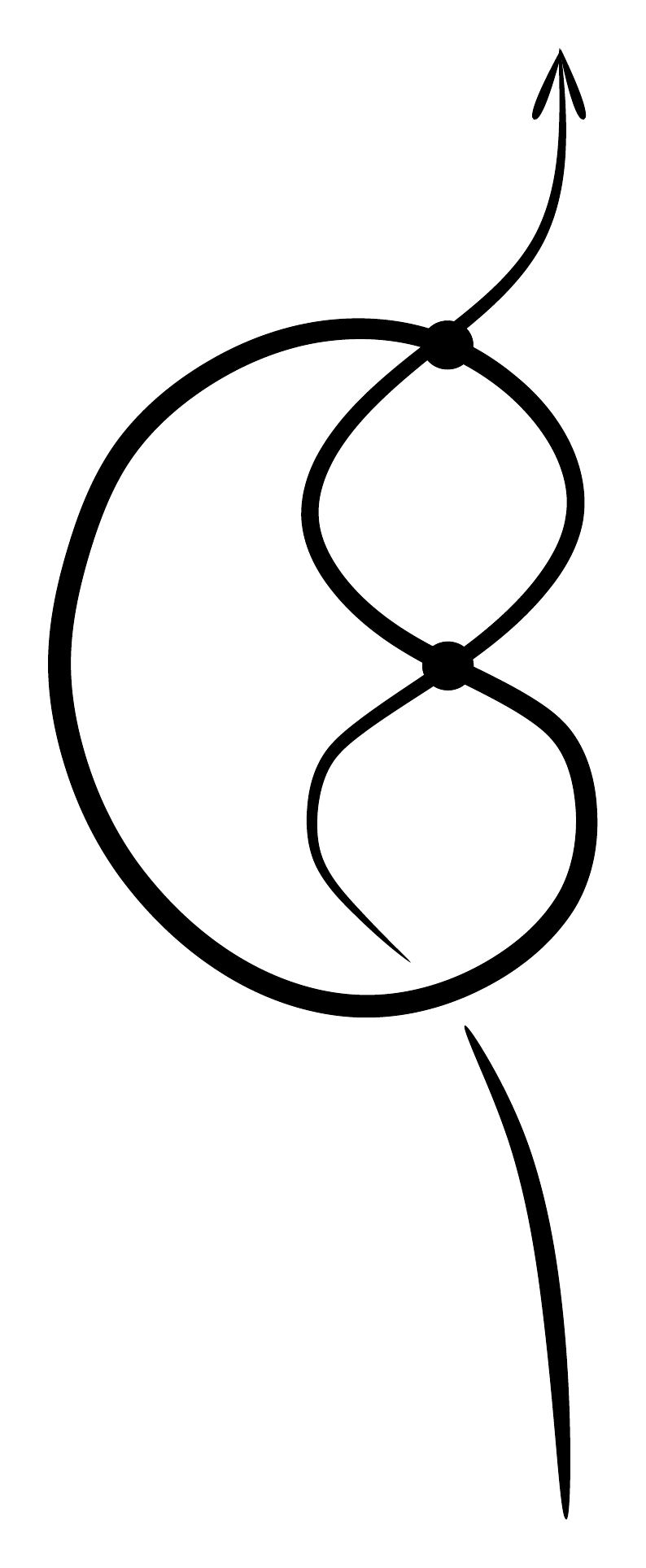}}
\caption{The singular knot $k$ obtained by joining the two components of the singular 2-component braid $h$.}
\label{F:KnotFromHopf}
\end{figure}

\medskip

Since $H \mapsto K$ and $2H$ was identified with $\pm [\mathbbm{1}_{n-2}, \mathbbm{1}_{n-2}]$, part (a) is proven in the case that $p=1$.  In particular, if the ambient dimension $n$ is $4$,$6$, or $10$, graphing the pure braid associated to the Hopf fibration gives a generator of the first nontrivial homotopy group of the space of long knots in $\R^n$.

We next prove part (b) for $p=1$.  We just have to consider the result of joining the two components in the cycle $L\pm L'$.  We saw that $L$ is homologous to $H$ and hence maps to $K$, so it remains to check that the result $K'$ of joining the components of $L'$ pairs trivially with $\kappa$.  Indeed, the singularities of $K'$ do not coincide with the chords of $\kappa$, so by Lemma \ref{L:Cocycles} (a) and Lemma \ref{L:Pairing}, $K'$ is zero in homology, and part (b) is proven for $p=1$.  A picture of $K'$ is shown in Figure \ref{F:SingularKnot}.

\begin{figure}[h!]
  \raisebox{-7pc}{\includegraphics[scale=0.15]{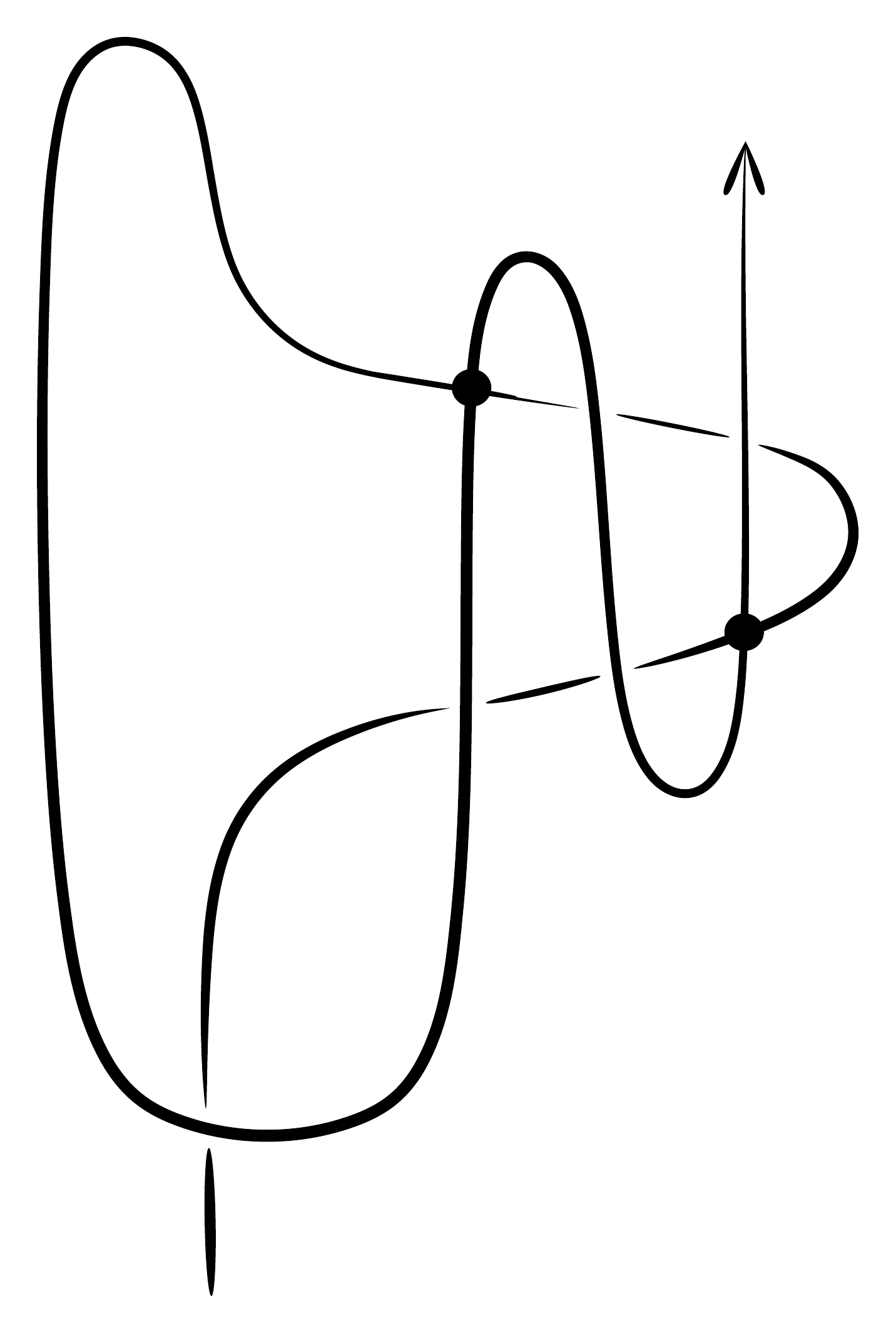}} \quad $\sim$ \quad
 \raisebox{-5pc}{\includegraphics[scale=0.15]{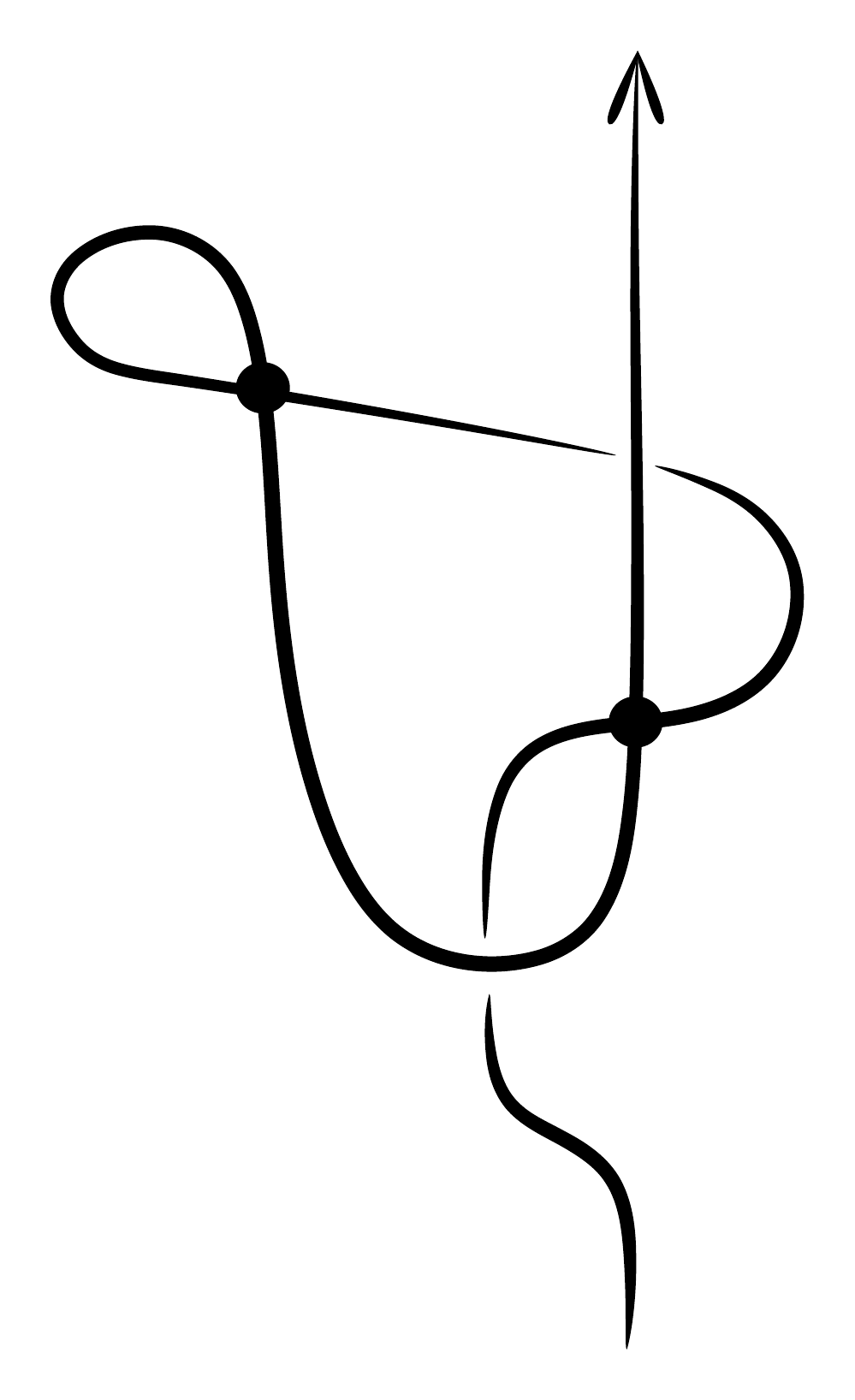}}
\caption{
The singular knot obtained by joining the two components of the singular link $\ell' =(\ell_1', \ell_2')$.
}
\label{F:SingularKnot}
\end{figure}

We now prove the statement for $p=1$ in part (c).  Recall that for $j=1,2$, we have sections $\iota_j: \K_1^n \to \LL{1}{n}$, as in Definition \ref{D:Inclusion}, of the restrictions $\rho_j: \LL{1}{n} \to \K_1^n$ from Definition \ref{D:Restriction}.
The classes $(\iota_1)_*(K)$, $(\iota_2)_*(K)$, and $L\pm L'$ are all in the image of the Hurewicz map because they are  the images of $[b_{21},b_{31}]$ under maps induced by maps of spaces.  
For $n$ even, we likewise identified $H$ with the image of an element in a homotopy group.
Next, it is easy to see that 
\begin{align*}
\begin{split}
\langle \rho_j^*I(\kappa), (\iota_i)_*(K)\rangle &= \delta_{ij} \text{ for } 1 \leq i,j\leq 2\\
\langle I(\eta), (\iota_i)_*(K)\rangle &= 0\\
\langle I(\lambda), (\iota_i)_*(K)\rangle &= 0.
\end{split}
\end{align*}
Thus for $n$ odd (respectively even) these 3 (respectively 4) classes are linearly independent.
Recall from Theorem \ref{T:TripodHedgehog} $\dim \pi_*(\LL{1}{n}) \otimes \Q$ is 3 (respectively 4) 
for $n$ odd (respectively even).  
Furthermore,  $(\iota_1)_*(K)$, $(\iota_2)_*(K)$, and $L\pm L'$ cannot be proper multiples of other classes because each one maps to a generator of $\pi_{2n-6}\K_1^n$.  The same is true for $H$ because by Theorem \ref{T:Lambda}, graphing from a sphere is a split injection.  So part (c)  is proven for $p=1$.

To prove parts (a), (b), and (c) in the case of arbitrary $p$, we use diagram \eqref{Eq:BigDiagram} below.   
Each horizontal arrow is a graphing map, and each vertical arrow is given by joining two components.
By Proposition \ref{P:JG=GJ}, each square commutes.

\begin{equation}
\label{Eq:BigDiagram}
\xymatrix@C1.0pc{
\pi_{2n-3p-3}\Omega^p \Conf(3,\R^{n-p}) \ar[r]^-{G_*}
& \pi_{2n-2p-4} \mathcal{L}_{1,1,1}^{n-p+1} \ar[r]^-{G_*} \ar[d]_-{J_*}
& \pi_{2n-2p-5} \mathcal{L}_{2,2,2}^{n-p+2} \ar[r]^-{G_*} \ar[d]_-{J_*} 
& \dots  \ar[r]^-{G_*} 
& \pi_{2n-3p-3}\mathcal{L}^{n}_{p,p,p} \ar[d]_-{J_*} \\
\pi_{2n-3p-3}\Omega^p \Conf(2,\R^{n-p}) \ar[r]^-{G_*} 
& \pi_{2n-2p-4} \LL{1}{n-p+1} \ar[r]^-{G_*} \ar[d]_-{J_*}
& \pi_{2n-2p-5} \LL{2}{n+2} \ar[r]^-{G_*} \ar[d]_-{J_*}
& \dots  \ar[r]^-{G_*} 
& \pi_{2n-3p-3}\LL{p}{n} \ar[d]_-{J_*} \\
& \pi_{2n-2p-4} \K_1^{n-p+1} \ar@{->>}[r]^-{G_*} 
& \pi_{2n-2p-5}\K_2^{n+2} \ar[r]_-\cong^-{G_*} 
& \dots  \ar[r] _-\cong^-{G_*} 
& \pi_{2n-3p-3}\K_p^n \\
}
\end{equation}

By work of Budney \cite{Budney:Family}, the horizontal maps in the bottom row are isomorphisms except possibly the first one.  
For odd codimension $n-p$, all the groups in the bottom row are $\Z$ and all the maps there are isomorphisms.
For even codimension $n-p$, the first group is $\Z$, while the remaining groups are $\Z/2$, and the first map is surjective while the rest are isomorphisms.  By commutativity, the compositions from the groups 
$\pi_{2n-3p-3}\Omega^p\Conf(m,\R^{n-p})$ (with $m=2,3$) to $\pi_{2n-3p-3} \K_p^{n}$ through $\pi_{2n-2p-4} \K_1^{n-p+1}$ are equal to the compositions through  $\pi_{2n-3p-3}\LL{p}{n}$ and $\pi_{2n-3p-3}\mathcal{L}^{n}_{p,p,p}$ respectively.  Thus we obtain parts (a) and (b) in the case of arbitrary $p$ from the case of $p=1$.

To prove the statement for $p\geq 2$ in part (c), recall that by Theorem \ref{T:GraphingDim>1} (a), 
the graphing maps $\pi_{2n-2p-4}\LL{1}{n-p+1} \to \pi_{2n-3p-3}\L_{p,p}^n$ take the form of surjections
\begin{align}
\Z^4 \oplus T(\pi_{2n-2p-3}S^{n-p-1}) \oplus F &\to \Z^4 \oplus T(\pi_{2n-2p-3}S^{n-p-1}),  & n-p \text{ odd},\\
\Z^3 \oplus \pi_{2n-2p-3}S^{n-p-1} \oplus F' &\to (\Z/2)^3 \oplus \pi_{2n-2p-3}S^{n-p-1}, & n-p \text{ even},
\label{Eq:IteratedGraphing2CompLinksEven}
\end{align}
where $T(A)$ denotes the torsion subgroup of $A$, and where $F$ and $F'$ are possibly nontrivial finite abelian groups.
We just need to check that their restrictions to the summands complementary to $F$ and $F'$ are also surjective, since by diagram \eqref{Eq:BigDiagram} and part (c) for $p=1$, the claimed generators are the images of generators for the domains of these restrictions.

By Theorem \ref{T:Lambda}, the restriction of graphing to $\pi_{2n-2p-3}S^{n-p-1}$ is injective, as is the restriction to its torsion subgroup.  From this fact, it is easy to deduce the desired surjectivity for $n-p$ odd.

For $n-p$ even, note that $T(\pi_{2n-2p-4}\LL{1}{n-p+1}) \cong \pi_{2n-2p-3}S^{n-p-1}$ maps to zero in $\pi_{2n-3p-3}\K_p^n$ because it factors through $\pi_{2n-2p-4}\K_1^{n-p+1} \cong \Z$.  By part (b), each of $[b_{21}, b_{31}]$, $(\iota_1)_*(K)$, and $(\iota_2)_*(K)$ maps downward to the generator of $\pi_{2n-3p-3}\K_p^n \cong \Z/2$ in diagram \eqref{Eq:BigDiagram}.  Thus the image of the $\Z^3$ summand under the map \eqref{Eq:IteratedGraphing2CompLinksEven} is complementary to the image of $\pi_{2n-2p-3}S^{n-p-1}$.  
The 3 generators under consideration are linearly independent because they lie in distinct summands in the decomposition 
\eqref{Eq:pi_iLpqNDecomp2}.
Considering the minimal cardinality of the image of their span gives the desired surjectivity.

We now prove the statement for $p\geq 3$ in part (c).  By Theorem \ref{T:ImagesOfiS}, any class in $\pi_{2n-3p-3}\mathcal{L}_{m\cdot p}^n$ is in $\sum_{|S|\leq 3} \im (\iota_S)_*$.
By Corollary \ref{C:IntKerPlusSumIm}, the subgroup $\bigcap_{|S|\leq 2} \ker (\rho_S)_*$ of $\pi_{2n-3p-3}\mathcal{L}_{3\cdot p}^n$ is complementary to $\sum_{|S|\leq 2} \im (\iota_S)_*$ (which we just described for $p\geq 2$), and by Theorem \ref{T:GraphingDim>1} (b), $\bigcap_{|S|\leq 2} \ker (\rho_S)_*\cong \Z$.  
The image of $[b_{21}, b_{31}]$ under the map $\pi_{2n-3p-3}\Omega^p\Conf(3,\R^{n-p}) \to \pi_{2n-3p-3}\mathcal{L}_{3\cdot p}^{n}$ lies in this subgroup.
By part (b), it maps to a generator of $\pi_{2n-3p-3}\K_p^n$, so it must generate this subgroup. 
\end{proof}

\begin{remark}[Sharpness of Theorem \ref{T:Knots} (a)]
We cannot extend part (a) of Theorem \ref{T:Knots} on graphing spheres, joining 2-component links, and knots to the case where $n-p$ is even.  In that case, although the target is $\Z/2$ for $p\geq 2$, we deduce the failure of this extension from the facts (used in the proof of Theorem \ref{T:Knots} (c)) that $\pi_{2n-2p-3}S^{n-p-1}$ is finite, $\pi_{2n-2p-4}\K_1^{n-p+1} \cong \Z$, and diagram \eqref{Eq:BigDiagram} commutes.
\end{remark}

\begin{remark}[Sharpness of Theorem \ref{T:pi_iLinksIso} for $p=q$]
\label{R:SharpnessBijGraphing}
Theorem \ref{T:Knots} (c) shows that Theorem \ref{T:pi_iLinksIso} does not extend to $\pi_{2n-6} \LL{1}n$, even modulo knotting.  Indeed, the class $L\pm L'$, obtained by joining two components of the parametrized long Borromean rings, comes neither from $\pi_{2n-5}S^{n-2}$ nor from knots.  
Relatedly, Budney's description of the bottom row of diagram \eqref{Eq:BigDiagram} and the splitting 
\eqref{Eq:LinksToKnotsFibn2Comp} shows that the graphing map from $\pi_{2n-6}\LL[0]{1}{n}$ is not an isomorphism, and hence Theorem \ref{T:pi_iLinksIso} is also sharp when $p=q-1$.
\end{remark}

\begin{remark}[Other ways of joining long Borromean rings]
We may join a pair of the three components of $F\pm F'=[b_{21}, b_{31}] \in \pi_{2n-6}\Omega\Conf(3,\R^{n-1})$  in any other of the six possible ways and ask whether we get different classes for $L\pm L' \in \pi_{2n-6} \LL{1}{n}$.
Arguments as in the proof of Theorem \ref{T:Knots} show that in even codimension ($n$ odd), we get the same class for $L+L'$.
In odd codimension ($n$ even), the resulting classes $L-L'$ all yield $\pm 1$ when paired with $I(\eta)$ or $I(\lambda)$.  Since $H$ yields 1 and 0 respectively when paired with these classes, Theorem \ref{T:Knots} holds for any of these joining maps, though one may obtain different bases.
\end{remark}

The next result has been independently obtained in forthcoming work of Gauniyal and Turchin.

\begin{theorem}
\label{T:Borr}
For any $k\geq 2$, let $B$ be the high-dimensional spherical Borromean rings $\coprod_1^3 S^{2k-1} \incl S^{3k}$, defined by 
\begin{align*}
\vec{x}=\vec{0} \quad \text{and} \quad |\vec{y}|^2/4 + |\vec{z}|^2&=1, \\
\vec{y}=\vec{0} \quad \text{and} \quad |\vec{z}|^2/4 + |\vec{x}|^2&=1, \\
\text{and} \quad \vec{z}=\vec{0} \quad \text{and} \quad |\vec{x}|^2/4 + |\vec{y}|^2&=1,
\end{align*}
where $(\vec{x}, \vec{y}, \vec{z})=(x_1, \dots, x_k, y_1, \dots, y_k, z_1, \dots, z_k)$ are coordinates on $\R^{3k}$ and where $S^{3k}$ is its one-point compactification.
Up to a sign, $B$ is isotopic to the image of (a representative of) $[b_{21}, b_{31}]$ under the composition of the $(2k-1)$-fold graphing map $\Omega^{2k-1}\Conf(3,\R^{k+1})\to  \mathcal{L}_{3\cdot(2k-1)}^{3k}$ followed by the closure map.
\end{theorem}

\begin{proof}
Put $p=3$ and $n=k+4$.
By Theorem \ref{T:GraphingDim>1} (b),
graphing maps 
$\pi_{2n-3p-3}\mathcal{L}_{p,p,p}^n$ isomorphically onto $\pi_0\mathcal{L}_{3\cdot(2k-1)}^{3k}$.  By Lemma \ref{L:SphericalLinks}, closure takes the latter group isomorphically onto $\pi_0\Emb\left(\coprod_1^3 S^{2k-1}, \, S^{3k}\right)$.  
As noted in the last paragraph of the proof of Theorem \ref{T:Knots}, iterated graphing maps $[b_{21}, b_{31}]$ to a generator of the subgroup 
given by the intersections of the kernels of the restrictions $\rho_S$ to fewer than 3 components.
Since $B$  generates this subgroup of $\pi_0\Emb\left(\coprod_1^3 S^{2k-1}, \, S^{3k}\right)$ \cite{Massey:1969}, the result follows.
\end{proof}

It is known that for $k$ even, $\pi_0 \Emb(S^{2k-1}, S^{3k})$ is generated by the Haefliger trefoil, which is defined as the result of joining all three components of the high-dimensional Borromean rings $B$ \cite{Haefliger:1962Annals}.  
For $k$ odd, where this group is $\Z/2$, the Manifold Atlas Project website \cite[\S 2]{Skopenkov:Knots} suggests that the analogous question is open.  Theorem \ref{T:Borr} allows us to answer it affirmatively:

\begin{corollary}
\label{C:HaefligerTrefoilEvenCodim}
For $k$ odd, the Haefliger trefoil generates the group $\pi_0 \Emb(S^{2k-1}, S^{3k}) \cong \Z/2$.
\end{corollary}
\begin{proof}
By Theorem \ref{T:Knots} (b),  $J_*J_* G^{2k-1}_* [b_{21}, b_{31}]$ represents a generator of 
$\pi_0 \K_{2k-1}^{3k} := \pi_0 \Emb_c (\R^{2k-1}, \R^{3k})$, where $J$ and $G$ are joining and graphing maps.
By Lemma \ref{L:SphericalLinks}, the closure of this long link thus represents a generator of $\pi_0 \Emb(S^{2k-1}, S^{3k})$.
Joining components commutes with closure, 
so this class is the result of joining components in the closure of $G^{2k-1}_* [b_{21}, b_{31}]$.  By Theorem \ref{T:Borr}, that closed link is isotopic to the result of joining together all three components of $B$.
\end{proof}

\begin{remark}[Knots, links, and braids in $\R^3$ and Vassiliev invariants]
\label{R:ClassicalAnalogue}
The constructions in the proof of Theorem \ref{T:Knots} can be applied when $n=3$.  
In fact, for the specific cocycles considered, the configuration space integral map $I$ gives link invariants in this setting.
The dual families are parametrized by $(S^0)^2$, so 
each of the classes $F + F'$, $L + L'$, $H$, and $K$ in $H_0 \mathcal{L}_{m \cdot 1}^3$ is an alternating sum of four links.  More generally, one can view a link with $i$ double points as an alternating sum of $2^i$ links. 
Now $H_0(\mathcal{L}_{m\cdot 1}^3; \, \R)=\R[\pi_0\mathcal{L}_{m\cdot 1}^3]$ fits into the Vassiliev--Goussarov filtration
\[
\R[\pi_0\mathcal{L}_{m\cdot 1}^3] = J_0 \supset J_1 \supset \dots \supset J_i \supset \dots
\] 
where $J_i$ is the ideal generated by singular $m$-component links with at least $i$ double-points.
A Vassiliev invariant of type $r$ is a link invariant which vanishes on $J_{r+1}$.  
A similar setup applies to pure braids, in which case $J_i$ is the $i$-th power of the augmentation ideal in the group ring on the pure braid group.

The only invariants of type 1 are linking numbers $\mu_{ij}$ of strands $i$ and $j$.  
There is one knot invariant of type 2, called the Casson knot invariant (corresponding to $\kappa$).
The only invariants of type 2 for pure braids are triple linking numbers (corresponding to $\mu$) and 2-fold products of pairwise linking numbers (corresponding to $\eta$ and the $\nu_i$).  
There is an invariant of 2-component long links of type 2 corresponding to $\lambda$.
For any $m$ the dual space $(\R[\pi_0 \mathcal{L}_{m\cdot 1}^3]/J_3)^*$ can be given a basis consisting of the appropriate collection of these invariants.

One extends link  invariants linearly to linear combinations of links.  
Then for any $m$, one can obtain a basis for $\R[\pi_0 \mathcal{L}_{m\cdot 1}^3]/J_3$ by choosing the appropriate elements out of $F + F'$, $L + L'$, $H$, and $K$ and their images under including components.
The classes $F + F'$, $L + L'$ and $K$ are equivalent to classes in $\pi_0\mathcal{L}_{m \cdot 1}^3 \subset H_0\mathcal{L}_{m \cdot 1}^3$ in this stage of the Goussarov--Vassiliev filtration.  (We neglect signs in what follows.)  
A long trefoil is equivalent to $K$, as it has Casson invariant 1.  
A nontrivial 3-component pure braid commutator, whose closure will be the Borromean rings, is equivalent to $F+F'$, since it has the same type-2 invariants.
A long Whitehead link is equivalent to $L+L'$.
The class of $H \in H_0(\Omega S^1)$ is not primitive because there is of course no 2-strand pure braid (or long link) which has zero pairwise linking number but a nonzero square of its pairwise linking number.  This mirrors what happens for other odd values of $n$.
\end{remark}

\subsection{Conjectured extensions and questions}
\label{S:Questions}

We now turn to some conjectures and questions.  We begin with extensions of Theorem \ref{T:Knots} (c).  Recall Theorem \ref{T:TripodHedgehog} and the description of generators in its proof as ``tripod'' and ``2-hair hedgehog'' graphs.

\begin{conjecture}
\label{C:graph-complex} \ 
\begin{itemize}
\item[(a)]
For $n-p$ even and $p=1$, the class in $\pi_{2n-3p-3}\LL{p}{n}$ obtained by graphing and joining 2 components of $[b_{21}, b_{31}]$ corresponds, up to sign, to the 2-hair hedgehog graph with leaves labeled $1,2$.
\item[(b)]
For $n-p$ odd, the class in $\pi_{2n-3p-3}\LL{p}{n}$ obtained by graphing and joining 2 components of $[b_{21}, b_{31}]$ corresponds, up to a sign, to the difference of two tripods with leaves labeled $1,1,2$ and $1,2,2$ respectively.  
The class in $\pi_{2n-3p-3}\LL{p}{n}$ obtained by graphing $[\mathbbm{1}_{n-p-1}, \mathbbm{1}_{n-p-1}]$ corresponds, up to a sign, to the sum of two tripods with leaves labeled $1,1,2$ and $1,2,2$ respectively.  
\end{itemize}
\end{conjecture}

The two classes  in part (a) agree up to a scalar multiple, since they both generate the $1$-dimensional $\Q$-vector space $\ker (\rho_1)_* \cap \ker (\rho_2)_*$.
Similarly, in part (b), the two tripods and the two classes coming from graphing generate the $2$-dimensional $\Q$-vector space $\ker (\rho_1)_* \cap \ker (\rho_2)_*$.

To calculate $\pi_{2n-3p-3} \LL{p}{n}$ for all $p$ and $n$ with $p \leq n-3$, one would just need to calculate the torsion subgroup for $p=1$.  
To calculate $\pi_{2n-3p-3}\mathcal{L}_{m\cdot p}^{n}$ for all $p$ and $n$ with $p \leq n-3$, it suffices by Theorem \ref{T:ImagesOfiS} to calculate $\pi_{2n-9} \mathcal{L}_{2,2,2}^{n}$, and by Theorem \ref{T:TripodHedgehog} it remains only to calculate the torsion subgroup.  
(Indeed, that result and the proof of Theorem \ref{T:Knots} (c) imply that the non-torsion generators that we described for $\pi_{2n-3p-3}\mathcal{L}_{m\cdot p}^{n}$, $p\geq 3$, generate the corresponding rational homotopy group for all $p \geq 1$.)
Theorem \ref{T:Lambda} says that the torsion subgroup contains at least $\pi_{2n-2p-3}S^{n-p-1}$ in both cases, and we suspect that this lower bound is sharp:  

\begin{conjecture}
\label{C:ExtendThmF} 
If $m\geq 2$ and $1 \leq p\leq n-3$, then a minimal generating set for $\pi_{2n-3p-3}\mathcal{L}_{m\cdot p}^{n}$ is given by the $m$ inclusions of a generator of $\pi_{2n-3p-3}\K_p^n$; the result of graphing and then joining two components of $[b_{21},b_{31}]$ for every pair of components $(i,j)$ with $1\leq i < j \leq m$; the image under graphing of a minimal generating set of $\pi_{2n-2p-3}S^{n-p-1}$ for every pair $(i,j)$ with $1\leq i < j \leq m$; and if $m\geq 3$, the result of graphing $[b_{21},b_{31}]$ for every triple $(i,j,k)$ with $1 \leq i < j < k \leq m$.
\end{conjecture}


Potential future work includes studying the realization of classes of long links and long knots by graphing braids and joining components.

\begin{question}
\label{Q:PB}
(a)  Can all classes in $\pi_i\mathcal{L}_{m \cdot p}^n$ be obtained by graphing and then joining components of pure braids, i.e., representatives of elements in $\pi_i\Omega^p \Conf(q, \R^{n-p})$ for some $q\geq m$?  
(b)  If so, how does one find such a representative?
\end{question}

 
An affirmative answer to part (a) would be a higher-dimensional analogue of a theorem of Alexander \cite{Alexander:1923, Cohen-Gitler:LCBaT} that any knot can be obtained by joining components of a pure braid.
It would also be an analogue of a result of Stanford \cite{Stanford:arXiv1998}.  
His result says that knots are equivalent in the Vassiliev filtration quotient $J_0/J_r$ (described in Remark \ref{R:ClassicalAnalogue}) if and only if they differ by a pure braid in the $r$-th stage of the lower central series of the pure braid group.  

In higher dimensions, results on embedding spaces and graph complexes  \cite{Fresse-Turchin-Willwacher:Emb} suggest that the rational homology groups of spaces of links with codimension at least 3 can be viewed as generalizations of terms $J_{r-1}/J_r$,
with larger values of $r$ corresponding to homology groups in higher degrees and with the homotopy groups consisting of the primitive elements.  
The following facts also support this analogy: for a space of embeddings of codimension at least 3, the Taylor tower converges to that space, while for classical knots the Taylor tower is conjectured to be a universal Vassiliev invariant over $\Z$ \cite{BCSS:2005}.  Indeed, after some progress in previous joint work \cite{BCKS:2017}, the conjecture was proven over $\Q$ by Kosanovi\'c \cite{Kosanovic:arXiv2020} and up to $p$-torsion for (roughly) small primes $p$ by Boavida de Brito and Horel \cite{Boavida-Horel}.

However, a negative answer to part (a) seems likely.  
Elements with cycles in the graph complex computing rational homotopy (such as the ``2-hair hedgehog'' mentioned in 
the proof of Theorem \ref{T:TripodHedgehog}) lie in degrees incompatible with the graphing map.  
Indeed, Turchin conjectures that the graphing map is rationally zero on graphs with cycles.  
(Theorem \ref{T:Knots} (c) implies that graphing sends the hedgehog class to a nontrivial class, but it is torsion.)
Together with Proposition \ref{P:JG=GJ}, this suggests a negative answer to part (a).

On the other hand, one could restrict Question \ref{Q:PB} to $i=0$.
Over $\Q$, isotopy classes are given by trivalent trees modulo the Jacobi relation \cite{Songhafouo-Turchin:Forum}.  
Our previous joint work \cite{KKV:Primitive} shows that graphing realizes a proper subspace of these groups, though we did not check that our identification in terms of trees agrees with that of \cite{Songhafouo-Turchin:Forum}.
Here we conjecture that all elements of $\pi_0\mathcal{L}_{m \cdot p}^n \otimes \Q$ can be obtained by graphing 
elements of $\pi_i\Omega^p \Conf(q, \R^{n-p})$ and joining components.
Over $\Z$, isotopy classes are less tractable, but perhaps our question can be answered without a full description of these groups.
It would also be interesting to compare any progress on this question to the work of Kosanovi\'c on realizing isotopy classes of knots via
 graspers \cite{Kosanovic:Graspers}.

\medskip

In contrast to joining components of braids, one can also build $m$-component links out of links with fewer than $m$ components.  
One type of such map is given by the inclusion maps $\iota_S$ from Definition  \ref{D:Inclusion}.  In the classical setting where $p_1= \dots =p_m=1$ and $n=3$, an element in $\im (\iota_S\oplus \iota_T)$ for some proper, nonempty $S\subset \{1,\dots,m\}$ with $T := \{1,\dots, m\}-S$, is called a split link.
For $m=2$ or $3$, Theorems \ref{T:GraphingDim>1} and \ref{T:Knots} (c) implement the decompositions \eqref{Eq:IntKerPlusSumIm} and \eqref{Eq:pi_iLpqNDecomp2} in terms of split links, and the latter further decomposes the first summand in terms of braids and another type of link.  Our methods may extend to arbitrary $p_1, \dots, p_m$, though we focused on equidimensional source manifolds for simplicity. 

Another map that increases the number of components is cabling, where parallel components are added to a knot.  
We omit a precise space-level definition of a cabling map.  Clearly such a map, like the inclusion maps, cannot produce braids out of knots, and in the classical case, there are links that are neither split, nor cables, nor braids.  Thus it is too much to ask for a decomposition in terms of just split links and cables.
Batelier and Ducoulombier \cite{Batelier-Ducoulombier} used the Swiss cheese operad to 
decmopose the space of 2-component long links in $\R^3$ in terms of cables, split links, braids, and a remaining subspace, extending our previous work on 2-component long links \cite{Burke-K:2015, BBK:2015}.
In higher dimensions, one could consider families of links obtained by split links, cables, and pure braids, where these classes no longer form distinct components but instead subspaces.

\begin{question}
\label{Q:CablesSplitLinks}
Is there a systematic decomposition of the space of links $\mathcal{L}_{p_1, \dots, p_m}^n$, $n \geq 4$, in terms of 
the following: pure braids, split links, cables, and/or the results of graphing them and (in the case of pure braids) joining their components?
\end{question}

Such a decomposition could lead to constructions of explicit geometric representatives of homotopy and homology classes in those spaces of links.  The example of the 2-hair hedgehog class, which suggested a negative answer to Question \ref{Q:PB}, does not immediately suggest a negative answer to Question \ref{Q:CablesSplitLinks}; from a purely dimensional viewpoint, it could be the cabling a knot class corresponding to the same graph but with the same label on both leaves.

\begin{remark}
Preliminary computations similar to the ones in the proof of Theorem \ref{T:Knots} suggest that in the bases given in part (c), one can use a cable of the generator in $\pi_{2n-6} \K_1^n$ to replace $J_*[b_{21}, b_{31}]$ for $n$ odd.  However, for $n$ even, this cable spans the same subgroup as the image under graphing of $[\mathbbm{1}_{n-2}, \mathbbm{1}_{n-2}]$ so in this parity, $J_*[b_{21}, b_{31}]$ is neither split, nor a cable, nor a braid.  
\end{remark}


\appendix
\section{Injectivity of graphing for other spaces of 2-component links}
\label{S:Appendix}

We now cover some variants of Theorem \ref{T:A}, the injectivity of graphing for spaces of 2-component links.
In Section \ref{S:Spherical}, we prove Theorem \ref{T:LambdaSpherical}, an analogue of Theorem \ref{T:A} for spaces of spherical links in a sphere.  In Section \ref{S:Alpha}, we prove Theorem \ref{T:Alpha}, which is the analogue for spaces of long link maps, using the $\alpha$-invariant.  Injectivity in the latter setting applies only to homotopy groups in a range.  

\subsection{Injectivity of graphing for spherical links}
\label{S:Spherical}

For Theorem \ref{T:LambdaSpherical}, the analogue of Theorem \ref{T:Lambda} for spherical links in a sphere, we need to define the following subspace of links.  This is because for spherical links, the restriction to one component does not in general admit a section.  (For example, for $p=q=1$ and $n=3$, the component of the Hopf link is homotopy equivalent to $SO(4)$ \cite[Corollary 4.4(e)]{Havens-K}, while the component of the unknot is homotopy equivalent to $SO(4)/SO(2)$ \cite{Hatcher:SmaleConj}.  The long exact sequence in homotopy at $\pi_2$ and $\pi_1$ shows that the projection $SO(4)\to SO(4)/SO(2)$ admits no section.)

\begin{definition}
\label{D:Brunnian}
Define the space of \textbf{Brunnian long links} 
$\mathrm{BrEmb}_c\left(\coprod_1^m \R^{p_i}, \, \R^n\right)$
as the the subspace of links $f$ in 
$\Emb_c\left(\coprod_1^m \R^{p_i}, \, \R^n\right)$ 
such that the restriction of $f$ to any $m-1$ of its components is isotopic to a standard link $(e_1, \dots, e_{m-1})$.  
Define the space $\mathrm{BrEmb}\left(\coprod_1^m S^{p_i}, \, S^n \right)$ of {\bf Brunnian spherical links in a sphere} as the subspace of spherical links $f$ in $\Emb\left(\coprod_1^m S^{p_i}, \, S^n \right)$ such that $f$ is isotopic to a trivial link, where a spherical link is trivial if each component $S^{p_i}$ bounds a disk $D^{p_i+1}$ in $S^n$.
\end{definition}

Thus $\mathrm{BrEmb}_c\left(\coprod_1^m \R^{p_i}, \, \R^n \right)$ is a union of path components of  
$\Emb_c\left(\coprod_1^m \R^{p_i}, \, \R^n \right)$.
In many of the cases we consider, the latter space is path-connected, in which case all embeddings are Brunnian.  
The graphing and closure maps preserve the Brunnian property.

A version of the next map, sometimes called the linking coefficient, appears at least as far back as work of Zeeman \cite{Zeeman:LinkingSpheres, Zeeman:Isotopies} and Haefliger \cite{Haefliger:1962Top}.
In fact, the section given by graphing is sometimes called the Zeeman map.

\begin{definition}
\label{D:Lambda}
For $p \leq q <n$, define a map
\[ 
\pi_0\BrEmb(S^p \sqcup S^q, \ S^n)   \overset{\lambda}{\longrightarrow} \pi_p S^{n-q-1}
\]
as follows.  Denote an element of the domain by $b=(b_1,b_2)$.
Take any isotopy of $b_2$ to the standard embedding $\widehat{e}_2$, and extend it to an isotopy from $b$ to a link $b' = (b'_1, \widehat{e}_2)$.  
Define $\lambda(b)$ as the composite
\[
S^p \overset{b_1'}{\longrightarrow} S^n - S^q  
\longrightarrow S^{n-q-1} \x \R^{q+1}
\longrightarrow S^{n-q-1}
\]
where the second map is the diffeomorphism \eqref{Eq:SpheresComplementDiffeo} and the third map is projection onto the first factor.
\end{definition}

\begin{lemma}\cite[Th\'eor\`eme 7.1]{Haefliger:1966CMH}
The map $\lambda$ is well defined.  
\end{lemma}

\begin{proof}
We need to show that $\lambda(b)$ is independent of the choice of isotopy from $b_2$ to $\widehat{e}_2$ and the choice of an extension of it to the first component.  Suppose $b'$ and $b''$ are the endpoints of two such isotopies starting at $b$.  Then $b'$ and $b''$ are isotopic links, each restricting to $\widehat{e}_2$ on $S^q$ (though the links throughout the isotopy need not satisfy this property).  

An isotopy from $b'$ to $b''$ can be extended to a path of diffeomorphisms of $S^n$ starting at the identity.  Restrict this path to a disk $D^n$ which contains $b_1'(S^p)$ and which intersects $\widehat{e}_2(S^q)$ in a disk $D^q \subset D^n$.
The endpoints $g_0$ and $g_1$ of this path lie in the space $\Emb^+(D^n, S^n; \, D^q)$ of orientation-preserving embeddings $D^n \incl S^n$ which agree with the fixed inclusion on $D^q$.  
By a shrinking and linearization argument, $\Emb^+(D^n, S^n; \, D^q) \simeq SO(n-q)$.  In particular, it is path-connected, so there is a path in $\Emb^+(D^n, S^n; \, D^q)$ from $g_0$ to $g_1$.  

Lift this path to the space $\Diff^+(S^n; \, \widehat{e}_2(S^q))$ of diffeomorphisms of $S^n$ which fix $\widehat{e}_2(S^q)$ pointwise.  
We can do so by the isotopy extension theorem \cite[Theorem 8.1.3]{Hirsch:1976} or more generally the fact that the restriction $\Diff^+(S^n; \widehat{e}_2(S^q)) \to\Emb^+(D^n, S^n; D^q)$ is a fibration \cite{Palais:LocalTriv, Lima:1964}.
Restricting the latter path to the image of $b'$ gives an isotopy from $b'$ to $b''$ which leaves the second component fixed.
\end{proof}

\begin{theorem}
\label{T:LambdaSpherical}
Let $1 \leq p\leq q \leq n-2$, and let $i \geq 0$.  Then $\pi_i \Emb(S^p \scu S^q, \, S^n)$ has a direct summand of $\pi_{i+p} S^{n-q-1}$, 
the inclusion of which is induced by a based homotopy equivalence 
$j:S^{n-q-1} \to R^{n-p} - \R^{q-p}$, 
the inclusion $\Emb_c(\ast, \, \R^{n-p} - \R^{q-p}) \incl \Emb_c(\ast \scu \R^{q-p}, \, \R^{n-p})$
and the composition 
 $\Omega^p \Emb_c(\ast \scu \R^{q-p}, \, \R^{n-p}) \to \Emb_c(\R^p \scu \R^q, \, \R^n) \to \Emb_c(S^p \scu S^q, \, S^n)$
of the graphing map followed by the closure map.
\end{theorem}

Theorem \ref{T:LambdaSpherical} is valid for $q=n-1$ but trivial because we assume $p\geq 1$.  We impose that hypothesis because we defined the closure map only for codimension at least 2.

\begin{proof} 
To prove the theorem, it suffices to prove the analogous statement where $\Emb(S^p \scu S^q, \, S^n)$ is replaced by $\BrEmb(S^p \scu S^q, \, S^n)$ because the former space is a union of path components of the latter space.
We will show that the following composition of maps is the identity:
\begin{align}
\label{Eq:LambdaCompositionSpherical}
\begin{split}
\pi_{i+p} S^{n-q-1}  
\underset{\Omega^p(j)_*}{\overset{\cong}{\longrightarrow}} \pi_{i} \Omega^p (\R^{n-p} - 0^{n-q} \x\R^{q-p}) 
&\overset{\cong}{\longrightarrow} \pi_{i}  \Omega^p \BrEmb_c(\{\ast\} \sqcup \R^{q-p}, \, \R^{n-p})  \\
& \overset{G^p_*}{\longrightarrow} \pi_{i} \BrEmb_c(\R^p \sqcup \R^q, \, \R^n) \\
& \overset{\widehat{\cdot}}{\longrightarrow} \pi_{i} \BrEmb(S^p \sqcup S^{q}, \, S^n) \\
& \overset{G^i_*}{\longrightarrow} \pi_{0} \BrEmb(S^{i+p} \sqcup S^{i+q}, \, S^{i+n}) \\
&  \overset{\lambda}{\longrightarrow} \pi_{i+p}S^{n-q-1}
\end{split}
\end{align}
The first map comes from a homotopy equivalence as in formula \eqref{Eq:ComplementHtpyEqv}. 
The arrow just above $G_p^*$ involves an affine-linear map to obtain the required behavior at the boundary of $(-1,1)^n$. 
The map denoted $G^i$ (by abusing notation already used for long links) is the composition of $i$ maps, starting with 
\[
\Omega \BrEmb(S^p \sqcup S^q, \, S^n)  \to 
\BrEmb(S^{p+1} \sqcup S^{q+1}, \, S^{n+1}). 
\]
To obtain this first map, start by graphing a loop restricted to $I$ to get an embedding $S^p \x I \scu S^q \x I \to S^n \x I$.  
Then attach two disks $D^{n+1}$ to $S^n \x I$ using diffeomorphisms $\d D^{n+1} \to S^n \x \{\pm 1\}$.  In each copy of $D^{n+1}$, fix smooth proper embeddings of $D^{p+1}$ and $D^{q+1}$ and glue two copies of each of these to the embeddings of $S^p \x I$ and $S^q \x I$ to obtain smooth embeddings of $S^{p+1}$ and $S^{q+1}$.  We can use fixed embeddings of $D^{p+1}$ and $D^{q+1}$ independent of the embeddings of $S^p$ and $S^q$ because the loop of embeddings is based at the standard trivial link $\widehat{e}$.
The remaining $i-1$ maps needed to obtain $G^i$ are defined similarly.

We next apply the Pontryagin--Thom correspondence to $\pi_{i+p}S^{n-q-1}$.  A homotopy class is identified with a bordism class of framed manifolds by taking the preimage of a regular value of a smooth representative $\R^{i+p} \to S^{n-q-1}$.
Recall from Definition \ref{D:Closure} that the closure map uses disks of radii $\sqrt{n}$ contained in $t_i^*\x \R^{n-1}$ and that $t_2^* - t_1^*=2/3$.
We fix a regular value $y \in S^{n-q-1}$ far from the first coordinate axis and more specifically in the neighborhood of vectors whose angle with $(0,1,0,\dots,0)$ is less than $\arctan(1/(3\sqrt{n}))$.
This guarantees that the disks used in the closure do not contribute to the bordism class, or more precisely that the framed submanifold representing a bordism class in the target copy of $\pi_{i+p}S^{n-q-1}$ lies in $\R^i \x D^p \subset \R^{i+p}$.  

By this choice of regular value, we need only consider the restrictions $D^{i+p} \scu D^{i+q} \to D^{i+n}$ of elements of $\BrEmb(S^{i+p} \sqcup S^{i+q}, \, S^{i+n})$ to evaluate their image under $\lambda$.  
This allows us to essentially reduce the composite \eqref{Eq:LambdaCompositionSpherical} to a composite of maps
involving long links.  
In this setting, $\lambda$ is homotopic to the composite of the homotopy left-inverse $r$ to $G^{i+p}$ from formula \eqref{Eq:EmbToMap} followed by an iterated looping of the homotopy inverse to $j$.  Thus the composite \eqref{Eq:LambdaCompositionSpherical} is indeed the identity.
\end{proof}

\subsection{Injectivity of graphing for link maps in a range}
\label{S:Alpha}
Our next main result is proven using the $\alpha$-invariant of spherical links in Euclidean space.  
It is essentially a Gauss map and appears at least as early as the work of Massey and Rolfsen \cite{Massey-Rolfsen}.  Further developments on it include work of Koschorke \cite{Koschorke:1988, Koschorke:1997} and A.~Skopenkov \cite{Skopenkov:2000}.  

\begin{definition}
\label{D:Alpha}
We define $\alpha$ as the map
\begin{align*}
\alpha: \Link_*(S^p \scu S^q, \, \R^n) &\longrightarrow \Map_*(S^p \x S^q, \Conf(2,\R^n))\\
g=(g_1,g_2) &\longmapsto 
((t,u) \mapsto (g_1(t), g_2(u))).
\end{align*}
\end{definition}

Now $\Conf(2,\R^n) \simeq S^{n-1}$, so if $p,q\leq n-2$, then the restriction of any map 
$S^p \x S^q \to \Conf(2,\R^n)$ to $S^p\vee S^q$ is nullhomotopic.
Thus the map induced by $\alpha$ on $\pi_i$ for any $i \geq 0$ can be written as 
\begin{align*}
\alpha_*: \pi_i\Link_*(S^p \scu S^q, \, \R^n) &\to \pi_i\Map_*((S^p \x S^q)/(S^p \vee S^q), \, \Conf(2,\R^n))\\
&\cong  \pi_{i+p+q}\Conf(2, \R^n)\\
&\cong  \pi_{i+p+q}S^{n-1}
\end{align*}

\begin{theorem}
\label{T:Alpha}
Suppose $1 \leq p \leq q \leq n-2$ and $i\geq 0$. 
Consider the map $\pi_{i+p}S^{n-q-1} \to \pi_i \Link_c(\R^p \scu \R^q, \, \R^n)$ induced by 
a based homotopy equivalence $j:S^{n-q-1} \to R^{n-p} - \R^{q-p}$, 
the inclusion $\epsilon: \Link_c(\ast, \, \R^{n-p} - \R^{q-p}) \incl \Link_c(\ast \scu \R^{q-p}, \, \R^{n-p})$, and 
graphing $G^p:  \Link_c(\ast \scu \R^{q-p}, \, \R^{n-p}) \to \Link_c (\R^p \scu \R^q, \ \R^n)$.
Let $E^q: \pi_{i+p} S^{n-q-1} \to \pi_{i+p+q} S^{n-1}$ be the $q$-fold suspension map.
Then $(G^p \circ \epsilon \circ \Omega^p(j))_*$ is injective, nonzero, or split injective if the $q$-fold suspension map is respectively injective, nonzero, or an isomorphism.  
\end{theorem}

Putting $i=0$ above gives a result of Scott \cite{Scott:1968}, which was improved upon by Massey and Rolfsen \cite{Massey-Rolfsen}, Koschorke \cite{Koschorke:1990}, and Habegger and Kaiser \cite{Habegger-Kaiser}.
Like Theorem \ref{T:LambdaSpherical}, Theorem \ref{T:Alpha} is valid for $q=n-1$ but trivial because we assume $p\geq 1$.  Again, we impose that hypothesis because we defined the closure map only for codimension at least 2.

\begin{proof}
The main idea is to check that the composite $\pi_{i+p} S^{n-q-1} \to \pi_{i+p+q} S^{n-1}$, shown below in the composite \eqref{Eq:AlphaComposition}, agrees with the $q$-fold iterated suspension $E^q$.  
For notational ease, we will suppress $\Omega^p(j)_*$ and $\epsilon_*$ from the notation and write $\widehat{G}^p$ to denote the $p$-fold graphing $G^p$ followed by the closure map $\ \widehat{\cdot}\ $.  Thus we write the composite below as $[f] \mapsto \alpha_*\widehat{G}^p_*[f]$.
\begin{align}
\label{Eq:AlphaComposition}
\begin{split}
\pi_{i+p} S^{n-q-1}  
\overset{\cong}{\longrightarrow}  \pi_i \Omega^p (\R^{n-p} - \R^{q-p}) 
& \longrightarrow \pi_i \Omega^p \Link_c(\{\ast\} \scu \R^{q-p}, \, \R^{n-p})\\
& \xrightarrow{G^p_*} \pi_i \Link_c(\R^p \sqcup \R^q, \, \R^n) \\
& \overset{(\ \widehat{\cdot} \ )_*}{\longrightarrow} \pi_i \Link_*(S^p \sqcup S^q, \,\R^n) \\
& \overset{\alpha_*}{\longrightarrow} \pi_{i+p+q} S^{n-1}
\end{split}
\end{align}
All of the above maps are homomorphisms, even if $i=0$, because $G^p$ is a map of H-spaces.

We will use the Pontryagin--Thom theorem to identify elements of the first and last homotopy groups in the composite \eqref{Eq:AlphaComposition} with bordism classes of manifolds.  
In particular, we will identify $[f]$ with a bordism class of framed submanifold of $\R^{i+p}$ of dimension $i+p+q-n+1$, and we will identify $\alpha_*\widehat{G}^p_*[f]$ with a class of framed submanifold of $\R^{i+p+q}$, also of dimension $i+p+q-n+1$.  
As in the proof of Theorem \ref{T:LambdaSpherical}, let $y\in S^{n-q-1}$ be a regular value of $f$ that lies in the neighborhood of vectors making an angle of less than $\arctan(1/(3\sqrt{n}))$ with $(0,1,0,\dots,0)$.  
This will simplify a consideration related to the closure map.
Let $X:= f^{-1}\{y\}$ be the framed bordism class corresponding to $[f]$.

The result of applying the first map in the composite \eqref{Eq:AlphaComposition} to a map $f: \R^{i+p} \to S^{n-q-1}$ is more explicitly written as the composition
\begin{align*}
\xymatrix@C2pc@R1pc{
\R^i \x \R^p \ar[r]^-f & S^{n-q-1}  \ar@{^(->}[r]^-i_-\simeq & S^{n-q-1} \x \R^{q-p+1} \ar[r]^-{h^{-1}}_-\cong  & \R^{n-p} - \R^{q-p}
}
\end{align*}
where $i$ maps $S^{n-q-1}$ to $S^{n-q-1} \x 0^{q-p+1}$.
Above $h$ is the diffeomorphism \eqref{Eq:ComplementDiffeo}.
Note that $h^{-1}i(S^{n-q-1})$ is the unit sphere in $\R^{n-q} \x 0^{q-p} \subset \R^{n-p}$.

The second arrow in the composite \eqref{Eq:AlphaComposition} comes from mapping
$\R^{n-p} - \R^{q-p}$ 
homeomorphically onto the subspace of
$\Link_c(\{\ast\} \scu \R^{q-p}, \, \R^{n-p})$ where the $\R^{q-p}$ component is standard.
Thus the restriction to the first component of the image of $[f]$ in $\pi_i \Link_c(\R^p \sqcup \R^q, \, \R^n)$ is represented by the map
\begin{align*}
\R^i \x \R^p &\longrightarrow (\R^{n-p}  - \R^{q-p}) \x\R^p \\
(s,t) &\longmapsto (h^{-1} i f (s,t), \, t).
\end{align*}
Hence on the subspace $D^p \x D^q \subset S^p \x S^q$, the class $\alpha_*(\widehat{G}^p_*[f])$ is represented by the composite 
\begin{equation}
\label{Eq:AlphaGComposite}
\xymatrix@R=0.2pc{
\R^i \x D^p \x D^q \ar[r] & \Conf(2, \R^n) \ar[r]^-\simeq &  S^{n-1}\\
(s,t,u)  \ar@{|->}[r]& \left( (h^{-1} i f (s,t),t), \, (0^{n-q}, u) \right) &
}
\end{equation}
%
One can check that $(y,0^q) \in S^{n-1}\subset \R^n$ is a regular value of the composite \eqref{Eq:AlphaGComposite}.  
Because of Definition \ref{D:Closure} of the closure map and our choice of regular value $y$, the preimage of $(y,0^q)$ under $\alpha \widehat{G}^p f$ is contained in $\R^i \x D^p \x D^q \subset \R^i \x S^q \x S^q$.
Let $Z:=(\alpha \widehat{G}^p f)^{-1}\{(y,0^q)\}$.  
From the formula in  \eqref{Eq:AlphaGComposite} and the fact that 
$h^{-1}i(S^{n-q-1})$ is the unit sphere in $\R^{n-q} \x 0^{q-p} \subset \R^{n-p}$, we deduce that
\begin{align*}
Z&=\{(s,t,0^{q-p},t) : h^{-1} i f (s,t) = h^{-1}i (y)\}  \\
&=\{(s,t,0^{q-p}, t) : f(s,t)=y \}
\end{align*}
where the second equality comes from the fact that $h^{-1}i$ is injective.  

The iterated suspension $E^q [f]$ corresponds to the image of $X$ under the inclusion $\iota_{\R^{i+p}, \R^{i+p+q}}:\R^{i+p} \incl \R^{i+p+q}$ of the first $i+p$ coordinates.  
Now
\[
\iota_{\R^{i+p}, \R^{i+p+q}} (X) = \{ (s,t, 0^q): f(s,t)=y\}.
\]
This submanifold (with its induced framing) is bordant to $Z$ (with its induced framing) via the ambient isotopy 
\begin{align*}
\R^{i+p+q} \x [0,1] &\to \R^{i+p+q} \\
((t_1, \dots, t_{i+p+q}),r) &\mapsto (t_1, \dots, t_{i+q}, 
 t_{i+q+1} + r t_{i+1}, \dots t_{i+q+p}+rt_{i+p}) .
\end{align*}

Thus the composition $\pi_{i+p}S^{n-q-1} \to \pi_{i+p+q}S^{n-1}$ shown in \eqref{Eq:AlphaComposition} is the $q$-fold suspension $E^q$, so the map $G^p_*: \pi_{i+p}S^{n-q-1} \to \pi_{i} \Link_c\left(\R^p \sqcup \R^q, \, \R^n \right)$ is injective or nonzero if $E^q$ is, and if $E^q$ is an isomorphism, $G^p_*$ maps isomorphically onto a direct summand.  
\end{proof}

One can immediately adapt Example \ref{Ex:LkNumExamples} on linking number classes from Theorem \ref{T:Lambda} to Theorem \ref{T:Alpha}, i.e., from embeddings to link maps.

\begin{remark}[Spaces of based spherical link maps in Euclidean space]
\label{R:BasedSphericalLinkMaps}
An analogue of Theorem \ref{T:Alpha} also holds for the space  $\Link_*\left(S^p \sqcup S^q, \ \R^n\right)$.  Indeed, the proof applies to $\widehat{G}^p_*$ just as well as it does to $G^p_*$.  
\end{remark}

\begin{remark}[Spaces of embeddings and the $\alpha$-invariant]
\label{R:AlphaEmbeddings}
Theorem \ref{T:Alpha} and its proof are equally valid for spaces of embeddings rather than link maps.  
However, for spaces of embeddings, Theorem \ref{T:Lambda} already gives a stronger result.
\end{remark}

     \bibliographystyle{alpha}    
    \bibliography{refs}

\end{document}